\documentclass[DIV=13]{scrartcl}
\usepackage[sb]{libertine} %
\usepackage[T1]{fontenc}
\usepackage{textcomp}
\usepackage[varqu,varl]{zi4}%
\usepackage[amsthm,upint]{libertinust1math} %
\usepackage[scr=boondoxo,bb=boondox]{mathalpha} %
\usepackage{bm}
\usepackage{makecell}
\usepackage{mathtools}
\usepackage{todonotes}
\usepackage{booktabs}
\usepackage{enumitem}
\usepackage{placeins}
\usepackage{lineno}
\usepackage{colortbl}
\usepackage{subcaption}
\usepackage{tikz}
\usetikzlibrary{arrows.meta,positioning,calc,decorations.pathreplacing}
\usepackage{pgfplots}
\pgfplotsset{compat=newest}%
\usetikzlibrary{external,fit,backgrounds, matrix, shapes.misc,decorations.pathreplacing,decorations.markings,positioning,calc,3d,arrows.meta,fadings,patterns}
\usepgfplotslibrary{groupplots,statistics,colorbrewer,fillbetween,patchplots}
\usepackage{pgfplotstable}
\pgfmathdeclarefunction{lg10}{1}{%
  \pgfmathparse{ln(#1)/ln(10)}%
}
\def\XMIN{-0.2}
\def\XMAX{10.2}

\pgfplotsset{
  mysharedx/.style={
    xmin=\XMIN, xmax=\XMAX,
  }
}
\usepackage[backend=biber]{biblatex}
\usepackage{bm}
\usepackage{xcolor}
\usepackage{float}
\RequirePackage{siunitx}
\RequirePackage[algoruled,norelsize, linesnumbered, lined, commentsnumbered]{algorithm2e}

\newcommand\restrict[1]{\raisebox{-.5ex}{$|$}_{#1}}
\usepackage[format=plain,labelfont=bf]{caption}
\usepackage{amsthm}
\newtheoremstyle{mystyle}%
{}%
{}%
{}%
{}%
{\bfseries}%
{.\quad}%
{ }%
{\thmname{#1}\thmnumber{ #2}\thmnote{ #3 }}%
\theoremstyle{mystyle} %
\newtheorem{theorem}{Theorem}[section]

\theoremstyle{mystyle} %

\theoremstyle{mystyle} %
\newtheorem{remark}[theorem]{Remark}

\usepackage[pdftex,colorlinks=true,linkcolor={blue!50!black},citecolor={red!50!black},urlcolor=blue]{hyperref}
\usepackage{doi}
\newcommand{\V}{\mathcal V}
\newcommand{\Vh}{\mathcal V_h}
\newcommand{\vh}{\mathbf v_h}

\newcommand{\Pt}{\mathbf{P}}

\newcommand{\vt}{\mathbf{v}}
\newcommand{\gt}{\mathbf{g}}
\newcommand{\xt}{\mathbf{x}}

\newcommand{\rt}{\mathbf{r}}

\newcommand{\wt}{\mathbf{w}}
\newcommand{\nt}{\mathbf{n}}

\newcommand{\taut}{\pmb{\tau}}

\newcommand{\sigmat}{\pmb{\sigma}}
\newcommand{\phit}{\pmb{\phi}}
\newcommand{\epsilont}{\pmb{\epsilon}}

\newcommand{\te}{\pmb{e}}
\renewcommand{\div}{\operatorname{div}}

\definecolor{myblue}{RGB}{0 83 139}
\definecolor{myred}{RGB}{114 16 69}
\definecolor{mygreen}{RGB}{0 94 0}
\title{A Hybrid Neural Network-Finite Element Method for
  the Viscous-Plastic Sea-Ice Model}
\author{Nils Margenberg\thanks{
    University of Magdeburg,
    Institute for Analysis and Numerics,
    Universit\"atsplatz 2,
    39104 Magdeburg,
    Germany,
    \texttt{\{nils.margenberg,carolin.mehlmann\}@ovgu.de}
  }
  \and Carolin Mehlmann\footnotemark[1]
}

\begin{document}
\maketitle%
\begin{abstract}
We present an efficient hybrid Neural Network–Finite Element Method (NN-FEM) for solving the viscous–plastic (VP) sea-ice model. The VP model is widely used in climate simulations to represent large-scale sea-ice dynamics. However, the strong nonlinearity introduced by the material law makes VP solvers computationally expensive, with the cost per degree of freedom increasing rapidly under mesh refinement. High spatial resolution is particularly required to capture narrow deformation bands known as linear kinematic features  in viscous-plastic models. To improve computational efficiency in simulating such fine-scale deformation features, we propose to enrich coarse-mesh finite element approximations with fine-scale corrections predicted by neural networks trained with high-resolution simulations. The neural network operates locally on small patches of grid
  elements, which is efficient due to its relatively small size and
  parallel applicability across grid patches. An advantage of this local approach is that it generalizes well to different
  right-hand sides and computational domains, since the network operates on small subregions rather than learning details tied to a specific choice of boundary conditions, forcing, or geometry. The numerical examples quantify the runtime and  evaluate the error
  for this hybrid approach with respect to the simulation of sea-ice deformations.
Applying the learned network correction enables coarser-grid simulations to achieve qualitatively similar accuracy at approximately 11 times lower computational cost relative to the high-resolution reference simulations. Moreover, the learned correction accelerates the Newton solver by up to \SI{10}{\percent} compared to runs without the correction at the same mesh resolution.

\end{abstract}
\section{Introduction}

The sea-ice cover in both hemispheres is characterized by narrow, elongated deformation zones known as linear kinematic features (LKFs)~\cite{Kwok2008}. These structures, which correspond to leads and pressure ridges, are important for the heat exchange between ocean, and atmosphere~\cite{Maykut1978}.
Accurately capturing LKFs in sea-ice models critically depends on the choice of sea-ice rheology~\cite{Hutter2022}. Currently, most climate models apply the viscous–plastic (VP) constitutive relation~\cite{Blockley2020}, although several alternative rheological formulations have been proposed, e.g.~\cite{Rampal2016}. Recent studies have shown that VP-based models can reproduce key statistical and scaling properties of observed LKFs~\cite{Hutter2020}.

The nonlinearity introduced by the viscous–plastic material law makes solving the equations computationally expensive, with the computational cost per degree of freedom increasing rapidly as the mesh is refined~\cite{Koldunov, StadlerMehlmann2023}. At high horizontal mesh resolution the computational cost of converged numerical sea-ice approximations can become as expensive as running the global ocean model~\cite{Koldunov}.
Since mesh resolution is a key factor determining both the accuracy and computational cost of simulating linear kinematic features (LKFs), we investigate whether a neural network can help to achieve qualitatively similar LKF patterns within the viscous–plastic sea-ice model on coarser horizontal grids, thereby reducing the overall computational cost of sea-ice–ocean simulations run at high spatial mesh resolution.

In addition to grid resolution, the simulation of LKFs in the viscous-plastic model depends on several factors, including solver convergence~\cite{Koldunov, Lemieux2012}, rheological parameters~\cite{Hutter2020}, wind forcing~\cite{Hutchings2005}, and numerical spatial discretization~\cite{Mehlmannetal2021}.
Most sea-ice models use structured quadrilateral meshes, such as the \emph{Community Ice Code} (CICE)~\cite{Hunke2015}, where velocities and stresses are discretized with bilinear finite element basis functions, or the \emph{Massachusetts Institute of Technology general circulation model} (MITgcm)~\cite{LOSCH2010}, which employs a finite-volume scheme with central differences. Recent developments introduce unstructured triangular finite element meshes, as in the \emph{Finite-Volume Sea-Ice–Ocean Model} (FESOM)~\cite{Danilov2015} and the sea-ice module of the \emph{Icosahedral Nonhydrostatic Weather and Climate Model}~\cite{Mehlmann2021}, or quasi-structured finite-volume approaches, as used in the \emph{Finite-Volume Community Ocean Model} (FVCOM) and the \emph{Model for Prediction Across Scales} (MPAS)~\cite{Gao2011, Petersen2019}.

 \paragraph*{Approach and Related Work} We formulate the method within a finite element discretization framework. The proposed hybrid neural network–finite element approach enriches robust, converged coarse-mesh solutions with fine-scale corrections predicted by a neural network trained on high horizontal resolution data. The core idea is to project the coarse-scale solution onto a refined mesh obtained through grid refinement. On this fine mesh, the neural network refines the solution and generates an updated right-hand side for the sea-ice momentum equation. This updated field is then projected back onto the coarse mesh and used to advance the simulation in the next time step.
 Owing to the locality of the finite element discretization (i.e., its compact stencil), the neural network operates locally, making the method particularly well suited for learning spatially localized structures such as linear kinematic features (LKFs). Furthermore, the localized nature of the approach renders it attractive for large-scale parallel computing architectures used in climate simulations. The concept is closely related to hybrid finite element–neural network methods developed to enhance coarse-scale approximations of the Navier–Stokes equations~\cite{margenberg2022nnmg,margenberg2023hybrid}. A dynamic super-resolution approach was recently proposed for the shallow-water equations, in which a U-Net is used to correct coarse ICON-O simulations~\cite{Witte2025DynamicSR}. This method achieves errors comparable to those of a mesh twice as fine while preserving key flow balances. This concept is closely aligned with our hybrid correction philosophy.
Machine learning has recently been applied to sea-ice modeling to improve
subgrid-scale representation and forecasting. Finn \textit{et\,al.} used
U-Net-based corrections for a Maxwell elasto-brittle sea-ice rheology, reducing forecast
errors by over 75\%~\cite{finn2023subgrid}. Durand \textit{et\,al.} developed
CNN surrogates for long-term sea-ice thickness prediction, outperforming
persistence and climatology baselines~\cite{durand2024cnn}.

\paragraph*{Contributions and Limitations}
Our contributions are threefold. First, we explore the concept of the multiscale hybrid finite element-neural network in the context of the highly nonlinear viscous-plastic sea-ice model and evaluate its performance with respect to the simulation of LKFs and the computational cost. The results show that our method achieves high accuracy while reducing time-to-solution relative to high resolution reference simulations.
 Second, we perform detailed ablation studies that investigate the
importance of the network's input quantities. Third, we generalize the use of
neural networks beyond mesh cells, allowing them to operate locally on
small patches of mesh elements. We provide computational studies on the
impact of patch size and number of predicted levels.

The presented hybrid neural network-finite element method (NN-FEM) strongly depends on the hierarchical mesh structure and the locality of the finite element approach to relate the coarse-scale approximation with the fine-scale counterpart.
Nevertheless, it is flexible regarding the numerical discretization details: the coarse-grid problem can be solved with
various standard solvers, and the local enrichment via neural networks could be
done without constructing a global fine grid. Other local discretization methods such as finite differences or finite volumes might also be a suitable setup. The NN-FEM is presented on quadrilateral mesh in the context of the viscous-plastic rheology. We are confident that the method can be applied to nested triangular meshes and be used with different sea-ice rheologies.

\section{\label{sec:problem} Governing Equations}
Let $\Omega\subset\mathbb R^2$ be a bounded Lipschitz domain and $I=(0,T] \subset \mathbb{R}$ be the time interval of interest. The sea-ice dynamics in the classical viscous-plastic sea-ice model of Hibler~\cite{Hibler1979} is
prescribed by three variables: the
horizontal velocity $\vt\colon\Omega\times I\to\mathbb R^2$, the thickness
$H\colon\Omega\times I\to[0,\infty)$ and the concentration $A\colon\Omega\times I\to[0,1]$. The equations of motion are given by the following system.
\begin{subequations}\label{eq:model}
 \begin{align}
  \rho_{\mathrm{ice}}H\bigl(\partial_t \vt + f_c\te_z\times ( \vt -\vt_w )\bigr)
  &= \div \sigmat + \taut(\vt), \label{eq:momentum}\\
  \partial_t A + \div (A\vt) &= 0, \text{in }\Omega\times I, \label{eq:concentration}\\
  \partial_t H + \div (H\vt) &= 0, \text{in }\Omega\times I, \label{eq:thickness}
 \end{align}
\end{subequations}
where $\rho_{\mathrm{ice}}$ is the sea-ice density, $f_c$ the Coriolis parameter,
$\te_z$ the out-of-plane unit vector ($\te_z\times(v_1,v_2)=(-v_2,v_1)$) and $\vt_\mathrm w$ the near surface ocean velocity.
The stress
$\sigmat$ is given by the VP rheology \eqref{eq:stress}.
We use $\div$ for the divergence on vectors and the row-wise divergence on
second-order tensors; spatial derivatives are with respect to $x\in\Omega$. We neglect the thermodynamics. Therefore, there are no source and sink terms in the balance laws (\ref{eq:concentration}) and (\ref{eq:thickness}).

\paragraph*{Drag Forces.}
The oceanic and atmospheric drag forces $\taut(\vt)$ (force per unit horizontal area) are given by
\begin{equation}
\taut(\vt)
= \taut_w(\vt)+\taut_a(\vt),
\label{eq:forcing}
\end{equation}
with $\taut_w(\vt)=C_{\mathrm w}\rho_{\mathrm w}\|\vt_{\mathrm w}-\vt\|_2(\vt_{\mathrm w}-\vt)$ and $\taut_a=C_{\mathrm a}\rho_{\mathrm a}\|\vt_{\mathrm a}\|_2\vt_{\mathrm a}$. Here,
 $\rho_{\mathrm w},\rho_{\mathrm a}>0$ are water and air densities,
$C_{\mathrm w},C_{\mathrm a}>0$ dimensionless drag coefficients, and
$\vt_{\mathrm w},\vt_{\mathrm a}:\Omega\times I\to\mathbb R^2$ prescribed the near surface atmospheric and oceanic flows.

\paragraph*{VP Rheology.}
The viscous-plastic rheology relates internal stresses $\sigmat$ to strain rates
\begin{equation}
\dot\epsilont \coloneq \tfrac12\bigl(\nabla\vt + \nabla\vt^\top\bigr)
\qquad
\dot\epsilont' \coloneq \dot\epsilont - \tfrac12\operatorname{tr}(\dot\epsilont)I,
\label{eq:strain}
\end{equation}
where $\operatorname{tr}$ denotes the trace, $\dot\epsilont'$ is the deviatoric part and $I$ the identity matrix.
The viscous--plastic
constitutive relation reads as follows.
\begin{equation}
\sigmat = 2\eta\dot\epsilont' + \zeta\operatorname{tr}(\dot\epsilont)I - \tfrac{P}{2}I.
\label{eq:stress}
\end{equation}
 The viscosities, $\eta, \zeta$, and sea-ice strength, $P$, are given as
\begin{equation}
\eta = e^{-2}\zeta,
\qquad
\zeta = \frac{P}{2\Delta(\dot\epsilont)}, \qquad P(H,A) = P^\star H \exp\bigl(-C(1-A)\bigr).
\label{eq:viscosities}
\end{equation}
Here, $e=2$ denotes the eccentricity of the elliptical yield curve and $P^\star, C$ are the ice strength parameters.
To achieve a smooth transition between the viscous and the plastic regime, we follow~\cite{Kreyscher2000, Mehlmann2017} and choose the following regularization:
\begin{equation}
  \Delta(\dot{\boldsymbol{\epsilon}})=\sqrt{\Delta_P(\dot{\boldsymbol{\epsilon}})^2+ \Delta_\text{min}^2}.
  \label{eq:delta}
\end{equation}
In case of the plastic regime, $\Delta_P(\dot{\boldsymbol{\epsilon}})$ is defined as
\begin{equation}
\Delta_P(\dot{\boldsymbol{\epsilon}})=\sqrt{ \frac{2}{e^2}\dot\epsilont':\dot\epsilont' + \bigl(\operatorname{tr}(\dot\epsilont)\bigr)^2 },
  \label{eq:delta_p}
\end{equation}
whereas the {viscous regime} is given as
\begin{equation}
 \Delta_\text{min}(\dot{\boldsymbol{\epsilon}})= 2\times 10^{-9}.
\end{equation}

\paragraph*{Initial and Boundary Data.}
The system is completed by prescribing the following initial and boundary data.
\begin{alignat*}{2}
  \vt&=\vt_0&&\quad\text{ on } \{t=0\}\times \Omega,\\
  A = A_0,\quad H &= H_0&&\quad\text{ on } \{t=0\}\times \Omega,\\
  \vt&= 0 &&\quad\text{ on }I\times \partial\Omega,\\
  A = A^{\text{in}},\quad H &= H^{\text{in}} &&\quad\text{ on } I\times\Gamma^{\text{in}},
\end{alignat*}
where $\Gamma^{\text{in}}:=\{x\in
\partial\Omega~\vert~\nt\cdot \vt < 0\}$ and $\nt$ is the outward unit vector.

\paragraph{Variational Formulation.}
In order to apply a finite element formulation to approximate the coupled system (\ref{eq:momentum})-(\ref{eq:thickness}), we state its variational formulation.
By multiplying the system (\ref{eq:momentum})-(\ref{eq:thickness}) with test functions and applying partial integration to the stress term, we obtain
\begin{subequations}\label{eq:model-var}
 \begin{align}
  (\rho_{\mathrm{ice}}H\bigl(\partial_t \vt+ f_c\te_z\times ( \vt -\vt_w )\bigr),\phit)_\Omega
  &= (\sigmat, \nabla \phit)_\Omega + (\taut(\vt),\phit)_\Omega \label{eq:momentumvar}\\
  (\partial_t A + \div (A\vt),\phi_A) &= 0, \label{eq:concentrationvar}\\
  (\partial_t H + \div (H\vt),\phi_H) &= 0,  \label{eq:thicknessvar}
 \end{align}
\end{subequations}

where $(\cdot, \cdot)$ is the $L^2$ scalar product on $\Omega$.

\section{Discretization in Space and Time.}\label{sec:discretization}
The viscous–plastic sea-ice model is implemented in the software library Gascoigne~\cite{beckerFiniteElementToolkit}.
To solve the coupled system (\ref{eq:momentum})–(\ref{eq:thickness}), we follow the standard procedure in the literature and apply a splitting approach in time. First, we solve the balance laws (\ref{eq:concentration})–(\ref{eq:thickness}) and then, using the updated sea-ice concentration and thickness, we compute the solution of the momentum equation.
Let \(0 = t^0 < t^1 < \dots < t^N = T\) be a uniform partition of \([0,T]\) with step size
\[
k \coloneq t^{n+1} - t^{n}.
\]
We set \(\vt^n = \vt(t^n)\), \(H^n = H(t^n)\), and \(A^n = A(t^n)\). Then the splitting in time reads as:
\begin{enumerate}
\item \emph{Transport equations:} Given \(\vt^{n}\), update \((A^{n+1}, H^{n+1})\) by
 solving the advection balances (\ref{eq:concentration})–(\ref{eq:thickness}) with fixed velocity at \(\vt^{n}\).
\item \emph{Momentum equation:} With \((A^{n+1}, H^{n+1})\) fixed, compute
 \(\vt^{n+1}\) from the momentum equation (\ref{eq:momentum}) at time \(t^{n+1}\).
\end{enumerate}
In the momentum solve, the external fields \(\vt_{\mathrm a}\), \(\vt_{\mathrm w}\) and coefficients depending on
\((A,H)\) are evaluated at \(t^{n+1}\). The transport equations are discretized implicitly in time using streamline diffusion. The choice is made for convenience based on the academic software library. Any other transport discretization would be possible. As the focus of the manuscript lies on the development of an efficient solving procedure for the sea-ice momentum equation, we outline the discretization of the sea-ice momentum in detail:

\paragraph{Discretization of the Sea-Ice Momentum Equation}
Owing to the stiffness of the momentum equation~(\ref{eq:momentum}), explicit time-stepping methods would necessitate prohibitively small time steps. Therefore, we adopt an implicit time-stepping strategy, as is standard in the literature, specifically the implicit Euler method. For the spatial discretization of the momentum equation we use a finite element discretization. Let $\{\mathcal T_h\}$ be a shape-regular,
quasi-uniform triangulation of $\Omega$ into convex quadrilaterals $K$, with
mesh size $h=\max_{K\in\mathcal T_h}\operatorname{diam}(K)$. Using the mesh
$\{\mathcal T_h\}$, we construct our finite element space using standard
isoparametric biquadratic Lagrange elements on each quadrilateral
$K\in\mathcal T_h$. The conforming spaces are
\begin{equation}\label{eq:FEspace}
 S_h\coloneq\{\phi_h\in H^1(\Omega)\;:\; \phi_h|_K\in\mathbb Q_2(K)\;\forall K\in\mathcal T_h\},\qquad
 \Vh\coloneq[S_h]^2,\qquad
 \Vh^0\coloneq\{\vh\in\Vh:\ \vh|_{\partial\Omega}=\mathbf 0\}\,,
\end{equation}
where $H^1(\Omega)$ consists of functions, defined on $\Omega$ with square-integrable derivatives.
In every time step $t^{n+1}$, we separate the nonlinear momentum
equation (\ref{eq:momentumvar}) into an operator on the unknown, $\mathcal{A}(\vt_h^{n+1}, {\phit_h})$, and a right-hand side, $\mathcal{F}({\phit_h})$. We seek a $\vt^{n+1}_h$ such that

\begin{equation}\label{eq:Newvar}
\mathcal{R}(\vt^{n+1}_h,\phit_h):= \mathcal{A}(\vt_h^{n+1}, {\phit_h}) -\mathcal{F}({\phit_h})=0, \, \forall \phit_h \in \Vh^0,
\end{equation}

where
\begin{align}
\mathcal{A}({\vt_h}^{n+1}, {\phit_h}) &:= (\rho_\text{ice}H_h^{n+1}\vt_h^{n+1},\phit_h)
+ k ( f_c \mathbf{e}_z \times \vt_h^{n+1}, {\phit_h})
+ k ({\sigmat}^{n+1}, \nabla {\phit}_h)
- k ({\taut}_w(t^{n+1}, \vt_h^{n+1}, {\phit_h}),\label{eq:A}\\
\mathcal{F}({\phit_h}) &:= (\rho^n \vt_h^{n},{\phit_h})
+ k ({\taut}_a(t^n), {\phit_h})
+ k ( f_c \mathbf{e}_z \times {\vt}_{\mathrm{w}}(t^n), {\phit_h}).
\end{align}
We solve the finite-dimensional nonlinear system of equations with a Newton-Krylov method, which is briefly outlined in the following.
Let \(\{\phit_h^i\}_{i=1}^{N}\) be a basis of \(\Vh^0\). Then any function $\vt_h \in \Vh^0$
can be represented as
\[
 \vt_h=\sum_{j=1}^{N} v_j\,\phit_h^j,
 \qquad \mathbf x=(v_1,\dots,v_{N})^\top\in\mathbb R^{N}.
\]
We introduce the discrete residual vector \(\rt\in\mathbb{R}^{N}\) by
\(r_i\coloneq\mathcal R(\vt^ {(n+1)}_h,\,\phit_h^i)\) for \(i=1,\dots,N\) and the right-hand side vector
$\mathbf{f} \in \mathbb{R}^{N}$ by \(f_i \coloneq \mathcal{F}(\phit_h^i)\) for \(i=1,\dots,N\).
Analogous to (\ref{eq:Newvar}) we seek \(\mathbf x^{n+1}\) such that
\(
 \rt(\mathbf x^{n+1}) = \mathbf 0.
\)
For a given iterate $\vt_h^{(k)}$, $k = 0, 1, 2, \dots$, with coefficient vector $\mathbf{x}^{(k)}$,
the Newton iteration based on the Jacobian $ \mathbf{J}$ reads
\begin{equation}
 \label{eq:newton-iterate}
 \mathbf{J}(\mathbf{x}^{(k)}) \, \delta \mathbf{x}^{(k)} = -\mathbf{r}(\mathbf{x}^{(k)}),
\end{equation}
followed by the update
\[
\mathbf{x}^{(k+1)} = \mathbf{x}^{(k)} + \alpha^{(k)} \, \delta \mathbf{x}^{(k)}.
\]
Here $\alpha^{(k)}$ is the line search parameter. Details on the assembly of the Jacobian by the bilinear form $\mathcal{A}(\cdot, \cdot) $ (\ref{eq:A}) are given in~\cite{Mehlmann2017}. The arising linear systems are solved by a preconditioned Krylov
method (e.\,g.\ GMRES preconditioned by geometric multigrid). The procedure is described in~\cite{Mehlmann2017}.
To conclude this section, we state the solution of the coupled problem with the
Newton-Krylov solve for the momentum equations in the $n$-th
timestep.

\begin{algorithm}[H]
 \caption{\label{alg:nk}Solution process of the momentum equation in one time step, $t_n \to t_{n+1}$.}
  \DontPrintSemicolon
  \textbf{Momentum equations:} Set initial guess $\vt_h^{(0)} \leftarrow \vt_h^{n}$ (use the solution of the previous time step $\vt_h^n$ as initial guess for the first step) \; %
  \For{$k=0,1,2,\dots$ until convergence}{
   \label{alg:nk:S1}Assemble residual vector \(\rt\) with entries \(r_i\coloneq\mathcal R(\vt_h,\,\phit_h^i)\) for \(i=1,\dots,N\).\;
   \label{alg:nk:S2}Assemble the Jacobian $\mathbf J^{(k)}$.\;
   \label{alg:nk:S3}Solve $\mathbf J(\mathbf x^{(k)})\delta\mathbf x^{(k)} = -\rt(\mathbf x^{(k)})$ by GMRES with geometric multigrid preconditioning (cf.\ \eqref{eq:newton-iterate}).\;
   \label{alg:nk:S4}Apply a line search and update $\vt_h^{(k+1)} \leftarrow \vt_h^{(k)} + \alpha^{(k)} \delta\vt_h^{(k)}$.\;
   \label{alg:nk:S}\uIf{$\|\rt^{(k+1)}\|\le \varepsilon_{\mathrm{nl}}\|\rt^{(0)}\|$}{\textbf{break}}
  }
  Accept $\vt_h^{n+1} \leftarrow \vt_h^{(k+1)}$.\;
\end{algorithm}
The solving procedure described in Algorithm~\ref{alg:nk} provides the numerical solver part
of the hybrid neural network-finite element method presented in this paper. The deep neural network in the hybrid
method is used to improve the numerical solution after Newton convergence. We
describe this process in detail in the next section.

\begin{remark}
 In our implementation, we employ a \emph{modified Newton solver} according to~\cite{Mehlmann2017}. There, globalization of Newton's
 method is achieved by splitting the analytical Jacobian into a symmetric
 positive definite part and a remainder. By adaptively damping the remainder,
 robust multigrid preconditioning is preserved while ensuring convergence. For
 a detailed formulation we refer to~\cite{Mehlmann2017}.
\end{remark}

\section{Hybrid Neural Network-Finite Element Method}\label{sec:dnnmg}
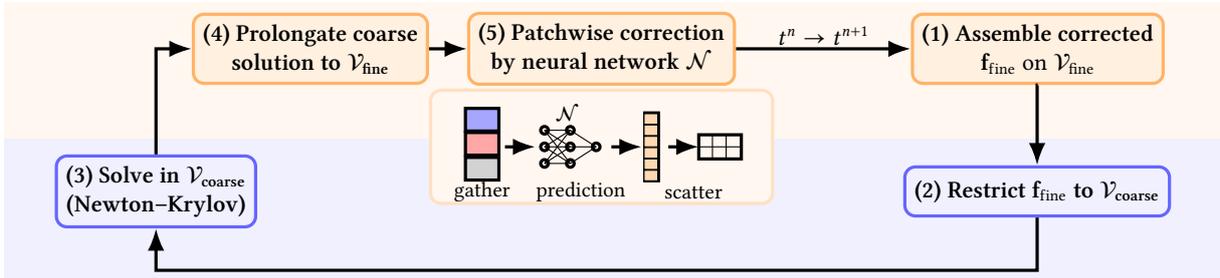
\begin{figure}[htb]
\centering
\begin{tikzpicture}[
 font=\footnotesize, >=Latex,
 node distance=5mm and 5mm,
 boxL/.style={draw=blue!60,rounded corners,very thick,fill=blue!6,minimum width=2.4cm,minimum height=6.5mm,align=center},
 boxF/.style={draw=orange!60,rounded corners,very thick,fill=orange!12,minimum width=2.3cm,minimum height=6.5mm,align=center},
 arrow/.style={->,very thick}
]
\fill[blue!6] (-2,-1.2) rectangle (14,0.7);
\fill[orange!6] (-2,0.7) rectangle (14,2.5);
\node[boxL] (solve) at (0,0) {\bfseries (3) Solve in \(\V_\text{coarse}\)\\\bfseries (Newton--Krylov)};
\node[boxL,right=85mm of solve] (rest) {\bfseries (2) Restrict $\mathbf f_{\mathrm{fine}}$ to \(\V_\text{coarse}\)};
\node[boxF] (prol) at ($(solve.north east)+(7mm,14mm)$) {\bfseries (4) Prolongate coarse\\
 \bfseries solution to \(\V_\text{fine}\)};
\node[boxF,right=of prol] (patch) {\bfseries (5) Patchwise correction\\\bfseries by neural network $\mathcal N$};
\node[boxF,right=of patch,xshift=18mm] (rhs) {\bfseries (1) Assemble corrected\\$\mathbf f_\mathrm{fine}$ on \(\V_\text{fine}\)};
\draw[arrow] (solve) |- (prol.west);
\draw[arrow] (prol) -- (patch);
\draw[arrow] (patch) -- node[above,yshift=-1mm]{\bfseries \(t^{n}\to
 t^{n+1}\)} (rhs);
\draw[arrow] (rhs.south) -| (rest.north);
\draw[arrow] (rest.south) |- ($(rest.south)+(-7mm,-7mm)$) --
($(rest.south)+(-75mm,-7mm)$) -| (solve.south);
\begin{scope}[shift={($(patch)+(0,-1.3)$)}, x=1cm, y=1cm]
 \node[draw=orange!24,rounded corners,fill=orange!6,very thick,inner
 sep=2pt,minimum width=4.5cm,minimum height=1.5cm] (mini) {};
 \node[font=\scriptsize] at ($(mini.west)+(0.70,-0.6)$) {gather};
 \node[font=\scriptsize] at ($(mini.west)+(2,-0.6)$) {prediction};
 \node[font=\scriptsize] at ($(mini.west)+(3.48,-0.6)$) {scatter};
 \begin{scope}[shift={($(mini.west)+(0.48,0)$)}]
  \draw[very thick,fill=blue!35]  (0.00, 0.21) rectangle +(0.44,0.3);
  \draw[very thick,fill=red!35]  (0.00,-0.11) rectangle +(0.44,0.3);
  \draw[very thick,fill=gray!35]  (0.00,-0.44) rectangle +(0.44,0.3);
 \end{scope}
 \draw[->,very thick] ($(mini.west)+(1,0)$) -- ($(mini.west)+(1.4,0)$);
 \begin{scope}[shift={($(mini.west)+(1.52,0)$)}]
  \foreach \y in {0.22,0.00,-0.22} \filldraw[very thick,fill=white] (0,\y) circle (0.05);
  \foreach \y in {0.22,0.00,-0.22} \filldraw[very thick,fill=white] (0.34,\y) circle (0.05);
  \filldraw[very thick,fill=white] (0.68,0.00) circle (0.05);
  \foreach \yA in {0.22,0.00,-0.22}
   \foreach \yB in {0.22,0.00,-0.22}
    \draw (0,\yA) -- (0.34,\yB);
  \foreach \yB in {0.22,0.00,-0.22} \draw (0.34,\yB) -- (0.68,0.00);
  \node[font=\scriptsize] at (0.3,0.45) {$\mathcal N$};
 \end{scope}
 \draw[->,very thick] ($(mini.west)+(2.4,0)$) -- ($(mini.west)+(2.8,0)$);
 \begin{scope}[shift={($(mini.west)+(3.35,0)$)}]
  \draw[very thick] (-0.33,-0.45) rectangle +(-0.18,0.90);
  \foreach \y in {-0.45,-0.3,-0.15,0.0,0.15,0.3}
   \draw[fill=orange!30] (-0.34,\y) rectangle +(-0.18,0.14);
  \draw[->,very thick] (-0.2,0.00) -- (0.2,0.00);
  \draw[very thick] (0.20,-0.16) rectangle +(0.56,0.32);
  \foreach \x in {0.20+0.18,0.20+0.36} \draw (\x,-0.16) -- (\x,0.16);
  \foreach \y in {-0.32+0.16,-0.32+0.32} \draw (0.20,\y) -- (0.76,\y);
 \end{scope}
\end{scope}
\end{tikzpicture}
\caption{\label{fig:overview-dnn-mg}
 \emph{ Compact overview of one hybrid NN–FEM time step $t^{n}\to t^{n+1}$:
 (1) use the corrected state $\vt_\text{fine}^{n}$ from the previous
 step to assemble the fine-space right-hand side $\mathbf f_\text{fine}^{n+1}$ and
 (2) restrict it to the working space,
 $\mathbf f_\text{coarse}^{n+1} = \mathbf R\mathbf f_\text{fine}^{n+1}$,
 (3) solve the sea-ice momentum equations on $\V_\text{coarse}$ to obtain $\vt_\text{coarse}^{n+1}$,
 (4) prolongate $\vt_\text{coarse}^{n+1}$ to the auxiliary fine space,
 $\vt_\text{fine}^{n+1} = \mathbf P\vt_\text{coarse}^{n+1} \in \V_\text{fine}$,
 (5) apply a patch-wise neural-network correction based on the residual
 and state information to compute a fine-scale increment
 $\delta\vt_\text{fine}^{n+1} \in \V_\text{fine}$, and
  store the corrected fine-grid state
 $\vt_\text{fine}^{n+1} = \vt_\text{fine}^{n+1} + \delta\vt_\text{fine}^{n+1}$ for use
 in assembling $\mathbf f_\text{fine}^{n+2}$ at the beginning of the next time step.
}}
\end{figure}

In the NN-FEM framework, a neural network predicts a fine-scale correction in an auxiliary space, which is applied as a post-processing step to the standard finite element discretization of the momentum equations. Our approach extends the deep neural
network finite element approach proposed in~\cite{margenberg2022nnmg}.
Within each time step, after the Newton method of Algorithm~\ref{alg:nk} has
converged, we predict fine-scale, patch-local corrections in a
richer finite element space (Fig.~\ref{fig:overview-dnn-mg}), i.\,e.~using the same finite element space on a refined spatial mesh. This hybridization
aims to improve the numerical solution without re-linearization and without
solving on the finer mesh. This hybrid method is only applied to the momentum
equations, the transport update remains unchanged.

The hybrid NN-FEM approach yields a separation of roles: the classical finite
element method and Newton solver ensure stability and consistency on the working
mesh, while the neural network provides data-driven subgrid-scale enhancement in
a refined auxiliary space. The following subsections introduce the level
hierarchy, notation, and transfer operators used in our hybrid algorithm.

\subsection{Level Hierarchy and Notation}
Let $\{\mathcal T_0,\dots,\mathcal T_L\}$ be a sequence of
uniformly refined meshes for $\Omega$, where $l=0,\cdots L$ indicates the refinement level. Further, let $\mathcal T_{L+S}$ denote a
uniformly refined \emph{auxiliary} fine mesh obtained by $S\ge1$ refinements of
$\mathcal T_L$. To each refined mesh $\mathcal T_h \in \{\mathcal T_0,\dots, \mathcal T_{L+S}\}$ of we can relate a finite element discretization $\V_h^0$, based on standard isoparametric
biquadratic Lagrange elements, see (\ref{eq:FEspace}). To simplify notation, we write $\mathcal T_\text{coarse}$ and $\mathcal T_\text{fine}$ for the working (coarse‐resolution) and auxiliary (high‐resolution) meshes of levels $\mathcal T_L$ and $\mathcal T_{L+S}$, respectively, and analogously use $\V_\text{coarse}$ and $\V_\text{fine}$ for the corresponding finite element spaces $\V_L$ and $\V_{L+S}$. The mesh transfer from a fine finite element space to a coarser one is
accomplished with $L^2$-projections, known as restrictions, and from coarse to
fine ones with interpolations, known as prolongations. We write
$\mathbf P:\V_\text{coarse}\to \V_\text{fine}$ and
$\mathbf R:\V_f\to \V_{c}$ for prolongation and restriction.

\begin{remark}[On the role of Geometric Multigrid in the hybrid NN-FEM approach]
 The mesh hierarchy and transfer operators, $\mathbf P$ and $\mathbf R$, are readily available within the
 Multigrid framework. However, the Geometric Multigrid method, which used as part of the GMRES to solve the linear problems arising in each Newton iteration (see Section \ref{sec:discretization}) is not a necessary
 part of the hybrid NN-FEM approach. Any other solver can be used to approximate the linear problem.
\end{remark}

\subsection{The Hybrid Algorithm}\label{sec:dnnmg-alg}
\begin{algorithm}[H]
 \caption{\label{alg:dnnmg}Hybrid neural network-finite element algorithm per time step $t^{n}\to t^{n+1}$}
 \DontPrintSemicolon%
 Assemble the fine-space right-hand side $\mathbf f_\text{fine}^{n+1}$ in $\V_\text{fine}$
 using the stored corrected state $\vt_\text{fine}^{n}$ from
 l.~\ref{alg:dnnmg:S7} of the previous step and restrict it to the working
 space:
 $\mathbf f_\text{coarse}^{n+1} \leftarrow \mathbf R\,\mathbf f_\text{fine}^{n+1}$.\;\label{alg:dnnmg:S8}%
 Solve equation~\eqref{eq:Newvar} in $\V_\text{coarse}$ with Newton-Krylov
 (Alg.~\ref{alg:nk}) using $\mathbf f_\text{coarse}^{n+1}$ to obtain
 $\vt_\text{coarse}^{n+1}$.\;\label{alg:dnnmg:S1}%
 Prolongate to the auxiliary mesh:
 $\vt_\text{fine}^{n+1}\leftarrow \mathbf P\,\vt_\text{coarse}^{n+1}$.\;\label{alg:dnnmg:S2}%
 Assemble the fine-space residual
 $\rt_\text{fine}^{n+1}$ on $\V_\text{fine}$ for $\vt_\text{fine}^{n+1}$
 (cf.\ Alg.~\ref{alg:nk}, l.~\ref{alg:nk:S1}).\;\label{alg:dnnmg:S3}%
 Gather patch-wise inputs from $\vt_\text{fine}^{n+1}$, $\rt_\text{fine}^{n+1}$,
 and geometry information, $\omega,$ into a matrix $\mathbf X$.\;\label{alg:dnnmg:S4}%
 Evaluate the network: $\mathbf D\leftarrow \mathcal N(\mathbf X)$
 (one row per patch).\;\label{alg:dnnmg:S5}%
 Scatter patch outputs into the global FE coefficient vector
 $\delta\vt_\text{fine}^{n+1}$ and update the fine-grid state
 $\vt_\text{fine}^{n+1}\leftarrow \vt_\text{fine}^{n+1}+\delta\vt_\text{fine}^{n+1}$.\;\label{alg:dnnmg:S6}%
 Store the corrected state $\vt_\text{fine}^{n+1}$ for assembling $\mathbf
 f_\text{fine}^{n+2}$ in the next time step.\;\label{alg:dnnmg:S7}%
\end{algorithm}
We summarize one time step of the hybrid NN-FEM correction. The overall
procedure is illustrated in Fig.~\ref{fig:overview-dnn-mg} and stated in
Algorithm~\ref{alg:dnnmg}. At the beginning of time step $t^{n+1}$ we assume
that a corrected fine-grid state $\vt_{\text{fine}}^{n} \in \V_\text{fine}$ from the
previous step is available. The state $\vt_{\text{fine}}^{n}$ is used to assemble
the fine-space momentum right-hand side
$\mathbf f_\text{fine}^{n+1}$, which is then restricted to the
working space via
\[
 \mathbf f_\text{coarse}^{n+1}
 \coloneq \mathbf R\mathbf f_\text{fine}^{n+1}
 \quad
 (\text{Alg.~\ref{alg:dnnmg}, l.~\ref{alg:dnnmg:S8}}).
\]
This way, the NN-FEM update $\delta\vt_\text{fine}^n$ from the previous step enters the new time step
as an explicit modification of the working space right-hand side
$\mathbf f_\text{coarse}^{n+1}$. All operators and the Newton-Krylov solve remain
on $\V_\text{coarse}$ (see also Appendix~\ref{sec:technical} for notation).

With the modified working mesh right-hand side $\mathbf f_\text{coarse}^{n+1}$,
we perform a Newton-Krylov solve on the working mesh
(Alg.~\ref{alg:nk}) to obtain the coarse solution
$\vt^{n+1}_\text{coarse} \in \V_\text{coarse}$. This state
$\vt^{n+1}_\text{coarse}$ is prolongated to the auxiliary fine space,
\[
 \vt^{n+1}_\text{fine}
 \coloneq \mathbf P\,\vt^{n+1}_\text{coarse}
 \in \V_\text{fine},
\]
and we assemble the corresponding fine-space momentum residual,
\[
 \mathbf r_\text{fine}^{n+1}
 \coloneq \mathbf r_\text{fine}\bigl(\vt^{n+1}_\text{fine}; \mathbf P A^{n+1}, \mathbf P H^{n+1}\bigr),
\]
with $(A^{n+1}, H^{n+1})$ frozen and prolongated to $\V_\text{fine}$
(Alg.~\ref{alg:dnnmg}, l.~\ref{alg:dnnmg:S1}--\ref{alg:dnnmg:S3}).
Next, we localize the network input data
($\mathbf r_\text{fine}^{n+1}$, $\vt^{n+1}_\text{fine}$, and geometric
features, $\omega$) to overlapping patches, evaluate the neural network patch-wise, and
assemble a global fine-space increment,
\[
 \delta \vt_\text{fine}^{n+1} \in \V_\text{fine},
\]
from the localized, patch-wise network outputs
(Alg.~\ref{alg:dnnmg}, l.~\ref{alg:dnnmg:S4}--\ref{alg:dnnmg:S6}).
This gather-predict-scatter pipeline respects Dirichlet boundary conditions and
averages contributions at overlapping degrees of freedom (DoF). Compared to prior
implementations developed in the context of the Navier-Stokes
equations~\cite{margenberg2023hybrid}, we generalize the approach to larger patch
sizes and arbitrary refinement levels between $\V_\text{coarse}$ and
$\V_\text{fine}$. The corresponding operators are described in detail in
Sec.~\ref{sec:patch-ops}.
Finally, we overwrite $\vt^{n+1}_\text{fine}$ with the corrected fine-grid field
\[
 \vt^{n+1}_\text{fine}
 \coloneq \vt^{n+1}_\text{fine} + \delta \vt_\text{fine}^{n+1}
 \in \V_\text{fine},
\]
and store $\vt^{n+1}_\text{fine}$ for use in the assembly of the next
fine-space right-hand side
$\mathbf f_\text{fine}^{n+1}$ at the beginning of time step $t^{n+2}$
(Alg.~\ref{alg:dnnmg}, l.~\ref{alg:dnnmg:S7}--\ref{alg:dnnmg:S8}).
The restricted load
$\mathbf f_\text{coarse}^{n+2} = \mathbf R\mathbf f_\text{fine}^{n+2}$
then serves as the modified working space right-hand side in the subsequent
Newton-Krylov solve. Thus, the NN-FEM correction
$\delta \vt_\text{fine}^{n+1}$ acts as an explicit post-processing on the
auxiliary fine space, $\V_\text{fine}$, that feeds back into the $\V_\text{coarse}$ problem only through
the right-hand side: it does not alter the Jacobian of the Newton-Krylov
iteration and does not require a fine-grid solve on $\mathcal T_\text{fine}$.

\paragraph*{Summary of One Time Step.}
Given $(\vt^{n},A^{n},H^{n})$ on the coarse space $\V_\text{coarse}$ we modify the solution
step for the momentum equation as follows:
\begin{enumerate}
\item \emph{Momentum (Newton-Krylov on $\V_\mathrm{coarse}$):} With $(A^{n+1},H^{n+1})$ fixed and external fields
 evaluated at $t^{n+1}$, assemble the right-hand side $\mathbf f_\mathrm{fine}^{n+1}$ using
 $(\mathbf P A^{n+1},\mathbf P H^{n+1}, \vt^{n}_\text{fine})$ on $\V_\text{fine}$ and restrict
 to $\V_\mathrm{coarse}$: $\mathbf f^{n+1}_\text{coarse}=\mathbf R \mathbf{f}_\mathrm{fine}^{n+1}$.
 Compute $\vt^{n+1}_\text{coarse}$ by solving~\eqref{eq:Newvar} (step \ref{alg:dnnmg:S1} of Alg.~\ref{alg:dnnmg}) via Alg.~\ref{alg:nk}.
\item \emph{Correction (hybrid NN-FEM):} Prolongate $\vt_\text{coarse}^{n+1}$ to $\V_\text{fine}$,
 assemble the fine-space residual on $\V_\text{fine}$, and compute a patch-wise correction by
 Alg.~\ref{alg:dnnmg} (steps~\ref{alg:dnnmg:S2}--\ref{alg:dnnmg:S6}). Store the corrected
 state $\vt_\text{fine}^{n+1}$ for assembling $\mathbf f^{n+2}_\text{coarse}$ at the \emph{beginning} of the next time
 step (steps~\ref{alg:dnnmg:S7}--\ref{alg:dnnmg:S8} of  Alg.~\ref{alg:dnnmg}).
\end{enumerate}

\subsection{Neural Network Inputs}\label{sec:dnnnet}
Before introducing the localization and assembly operators
(Sec.~\ref{sec:patch-ops}), we state the information passed to the network and
the form of its output.

\paragraph{Inputs.}
On the auxiliary mesh $\mathcal T_\text{fine}$ (with fine space $\V_\text{fine}$), we construct, for each small local neighborhood (patch, i.e., a set of adjacent mesh elements), a fixed-length feature vector consisting of three blocks:
\begin{enumerate}
\item the prolongated velocity $\mathbf \Pt\vt_\text{coarse}=\vt_\text{fine}$ restricted to the neighborhood (state);
\item the corresponding entries of the momentum residual $\rt_\text{fine}$;
\item geometry descriptors $\boldsymbol\omega$ (edge lengths, aspect ratios,
 interior angles) of the cells in the neighborhood.
\end{enumerate}
The state $\vt_\text{fine}$ encodes the local sea-ice flow, the residual encodes local error $\rt_\text{fine}$, and the
geometry block $\boldsymbol\omega$ adapts the map to the local cell shape.

\paragraph{Output and Assembly.}
From each local input, the network returns a velocity \emph{correction} in the fine space,
$\V_\text{fine}$. All local corrections are combined into a single fine-space
correction $\delta\vt_\text{fine}$ in
$\V_\text{fine}$. Dirichlet values are enforced, and the result is added to $\vt_\text{fine}$
to obtain a corrected field. The concrete localization, averaging, and
boundary-handling operators are detailed in Sec.~\ref{sec:patch-ops}.

\subsection{Patch Localization of Inputs, Batched Prediction, and Global Assembly}\label{sec:patch-ops}
\begin{figure}[ht]
\centering
\begin{subfigure}{0.32\linewidth}
\centering
\begin{tikzpicture}[scale=0.9, line cap=round, line join=round, >=Latex, font=\scriptsize]
 \draw[very thick] (0,0) rectangle (4,4);
 \foreach \x in {1,2,3} \draw[thin,gray!70] (\x,0) -- (\x,4);
 \foreach \y in {1,2,3} \draw[thin,gray!70] (0,\y) -- (4,\y);
 \fill[blue!6] (2,1) rectangle (3,2);
 \draw[ultra thick,blue!70!black] (2,1) rectangle (3,2);
 \draw[densely dashed,blue!70!black] (2.5,1) -- (2.5,2);
 \draw[densely dashed,blue!70!black] (2,1.5) -- (3,1.5);
 \node[anchor=south,blue!70!black] at (2.5,2) {$M_{K} \subset
  \mathcal T_{L+1}$};

 \draw[thin,gray!70] (0.1,3.8) -- (0.6,3.8);
 \node[anchor=west] at (0.65,3.8) {coarse grid $\mathcal T_\text{coarse}=\mathcal T_L $};
 \draw[ultra thick,blue!70!black] (0.1,3.5) -- (0.6,3.5);
 \node[anchor=west] at (0.65,3.5) {coarse cell $K$};
 \draw[densely dashed,blue!70!black] (0.1,3.2) -- (0.6,3.2);
 \node[anchor=west] at (0.65,3.2) {children on $\mathcal T_{L+1}$};
\end{tikzpicture}
\caption{\label{fig:scp-j1}Single-cell patch $M_{K}$ ($S=1$), with patch size $N_M=0$.}
\end{subfigure}
\hfill
\begin{subfigure}{0.32\linewidth}
\centering
\begin{tikzpicture}[scale=0.9, line cap=round, line join=round, >=Latex, font=\scriptsize]
 \draw[very thick] (0,0) rectangle (4,4);
 \foreach \x in {1,2,3} \draw[thin,gray!70] (\x,0) -- (\x,4);
 \foreach \y in {1,2,3} \draw[thin,gray!70] (0,\y) -- (4,\y);
 \fill[green!6] (2,1) rectangle (3,2);
 \draw[ultra thick,green!50!black] (2,1) rectangle (3,2);
 \foreach \x in {2.25,2.5,2.75}
  \draw[densely dashed,green!50!black] (\x,1) -- (\x,2);
 \foreach \y in {1.25,1.5,1.75}
  \draw[densely dashed,green!50!black] (2,\y) -- (3,\y);
 \node[anchor=south,green!50!black] at (2.5,2) {$M_{K} \subset \mathcal T_{L+2}$};
 \draw[thin,gray!70] (0.1,3.8) -- (0.6,3.8);
 \node[anchor=west] at (0.65,3.8) {coarse grid $\mathcal T_\text{coarse}:=\mathcal T_L$};
 \draw[ultra thick,green!50!black] (0.1,3.5) -- (0.6,3.5);
 \node[anchor=west] at (0.65,3.5) {coarse cell $K$};
 \draw[densely dashed,green!50!black] (0.1,3.2) -- (0.6,3.2);
 \node[anchor=west] at (0.65,3.2) {children on $\mathcal{T_\text{fine}}:=\mathcal T_{L+2}$};
\end{tikzpicture}
\caption{\label{fig:scp-j2}Single-cell patch $M_{K}$ ($S=2$), witch patch size $N_M=0$.}
\end{subfigure}
\hfill
\begin{subfigure}{0.32\linewidth}
\centering
\begin{tikzpicture}[scale=0.9, line cap=round, line join=round, >=Latex, font=\scriptsize]
 \draw[very thick] (0,0) rectangle (4,4);
 \foreach \x in {1,2,3} \draw[thin,gray!70] (\x,0) -- (\x,4);
 \foreach \y in {1,2,3} \draw[thin,gray!70] (0,\y) -- (4,\y);
 \fill[orange!10] (2,0) rectangle (4,2);
 \draw[ultra thick,orange!80!black] (2,0) rectangle (4,2);
 \node[orange!80!black] at (3,0.75) {$\mathcal B$};
 \foreach \x in {2.5,3,3.5} \draw[densely dashed,orange!80!black] (\x,0) -- (\x,2);
 \foreach \y in {0.5,1,1.5} \draw[densely dashed,orange!80!black] (2,\y) -- (4,\y);
 \node[anchor=south,orange!80!black] at (2.9,2)
   {$M_{\mathcal B}\coloneq\bigcup_{K\in\mathcal B} M_{\{K\}}$};
 \draw[thin,gray!70] (0.1,3.8) -- (0.6,3.8);
 \node[anchor=west] at (0.65,3.8) {coarse grid $\mathcal T_\text{coarse}=\mathcal T_L$};
 \draw[ultra thick,orange!80!black] (0.1,3.5) -- (0.6,3.5);
 \node[anchor=west] at (0.65,3.5) {coarse block $\mathcal B$};
 \draw[densely dashed,orange!80!black] (0.1,3.2) -- (0.6,3.2);
 \node[anchor=west] at (0.65,3.2) {children on $\mathcal T_\text{fine}$=$\mathcal T_\text{L+1}$};
\end{tikzpicture}
\caption{\label{fig:mcp-j1}Multi-cell patch $M_{\mathcal B}$ (S=1), with patch size $N_M=1$.}
\end{subfigure}
\caption{\label{fig:patch-visualization-annotated}\emph{
Patch geometry and notation. Thin gray
 lines: coarse mesh $\mathcal T_\text{coarse}$. Dashed lines:
 auxiliary refinement on $\mathcal T_\text{fine}$ restricted to a patch. A \emph{single-cell} patch $M_{K}$ on $\mathcal T_\text{coarse}$ collects the $4^{S}$ child patches on $\mathcal T_\text{fine}$ corresponding to a single coarse cell $K \in \mathcal T_\text{coarse}$ (cf. panel (a) and panel (b). A \emph{multi-cell} patch is the union
 over a set of cells $K$, $\mathcal B\subset\mathcal T_\text{coarse}$ (cf. panel (c)).
 }}
\end{figure}
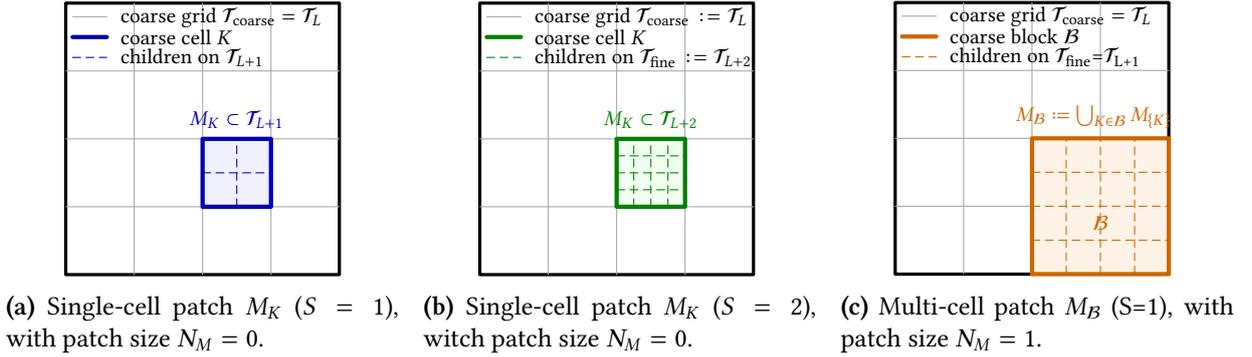

We first localize the prolongated field $\vt_\text{fine}$ and the residual
\(\rt_\text{fine}\) to small patches (cf.~Fig.~\ref{fig:patch-visualization-annotated}). These localized quantities are used
together with geometric information as the input to the network. All patches (collections of neighboring cells) are
processed together in a single batched network inference to predict local
velocity corrections. These local updates are then assembled into one coherent
fine-grid increment, $\delta \vt_\text{fine}$. In our setting, patches overlap,
i.\,e.\ some degrees of freedom receive multiple predictions from the network.
These duplicate predictions are averaged such that the global correction $\delta
\vt_\text{fine}$ lies in the fine space $\V_\text{fine}$. Boundary conditions
are enforced for $\delta \vt_\text{fine}$.

For each coarse cell \(K\in\mathcal T_\text{coarse}\), define its \emph{patch} \(M_K \subset \mathcal T_{L+S}\), $S \geq 1$,
as the union of its \(4^{S}\) child-patches on the auxiliary mesh
\(\mathcal T_\text{fine}\).  The patch size \(N_M\) denotes the number of uniform coarsening steps applied
to the working mesh $\mathcal{V}_{coarse}$ in order to construct the patch. Starting from the working
mesh, the patch is produced by coarsening it \(N_M\) times; hence, larger
values of \(N_M\) yield patches composed of fewer but correspondingly larger
coarse cells. As a consequence, increasing \(N_M\) expands the local context
and raises the number of degrees of freedom per patch (see
Tab.~\ref{tab:io-sizes-side}). The case \(N_M = 0\) represents the situation in
which the patch consists of a single cell of the working mesh
\(\mathcal{V}_{\mathrm{coarse}}\). In general, a patch of size \(N_M\)
comprises \(4^{N_M}\) cells of the working mesh
\(\mathcal{V}_{\mathrm{coarse}}\).
In Fig.~\ref{fig:patch-visualization-annotated} we show
examples for $S=1$  (Fig.~\ref{fig:scp-j1}) and $S=2$ combined with $N_M=0$ (Fig.~\ref{fig:scp-j2}). In Fig.~\ref{fig:scp-j1} and Fig.~\ref{fig:scp-j2} we show patches that consist of a single cell on the working mesh $\V_{\mathrm{coarse}}$ with a patch size of $N_M=1$.
Our implementation also supports multi-cell patches (cf.\ Fig.~\ref{fig:mcp-j1}). However, for
clarity, we present the construction for patches based on the working mesh
\(\mathcal T_\text{coarse}\) (cells \(K\in\mathcal T_\text{coarse}\)). The extension to
multi-cell patches is straightforward. Let
\(\mathbb P\coloneq\{M_K:\ K\in\mathcal T_\text{coarse}\}\) be the set of all
patches. Denote by \(N_{dof}^\text{fine}\coloneq\dim \V_\text{fine}\) the number of fine-grid velocity
DoFs (both components), by \(N_{\mathrm{cells}}^\text{coarse}\coloneq|\mathcal T_\text{coarse}|\)
the number of patches, and by \(N_{\mathrm{dof}}^{M}\) the number of velocity
DoFs supported on a single patch.

\paragraph{Gather.}
The goal of the gather step is to convert global fine-space vectors into a
batched, row-wise layout (one row per patch) for efficient, vectorized network
inference.
For each patch $M_K\in\mathbb P$, let
$\pi_K:\{1,\dots,N_{\mathrm{dof}}^{M}\}\to\{1,\dots,N_{dof}^\text{fine}\}$
be the local-to-global DoF map on $M_K$, and define the
\emph{selection} matrix $G_K\in\{0,1\}^{N_{\mathrm{dof}}^{M}\times N_{dof}^\text{fine}}$
by $(G_K)_{j,i}=1$ iff $i=\pi_K(j)$, else $0$.
Then, for any fine-grid coefficient vector $\mathbf x^\text{fine}\in\mathbb R^{N_\text{fine}}$,
the patch vector is
\[
 q_K = G_K\mathbf x^\text{fine} \in \mathbb R^{N_{\mathrm{dof}}^{M}}\,.
\]
Stacking one row per patch gives the \emph{batched} (row-wise) layout
\begin{equation}
 \label{eq:gather-op}
 \mathcal Q(\mathbf x^\text{fine})
 =
 \begin{bmatrix}
  G_1 \\ \vdots \\ G_{N_{\mathrm{cells}}^{c}}
 \end{bmatrix}\mathbf x^\text{fine}
 \in \mathbb R^{N_{\mathrm{cells}}^\text{coarse}\times N_{\mathrm{dof}}^{M}}.
\end{equation}
The matrices \(G_K\) are
fixed, extremely sparse, and can be stored as index lists.

\begin{figure}[t]
\centering
\begin{tikzpicture}[
 >=Latex, font=\small,
 mat/.style={draw,very thick,rounded corners=2pt,fill=gray!5},
 hdr/.style={align=center,font=\small,inner sep=1pt,fill=white,draw=white,rounded corners=1pt},
 colsep/.style={densely dashed,gray!70,very thick},
 rowlab/.style={anchor=east},
 brace/.style={decorate,decoration={brace,amplitude=5pt}},
 bracedown/.style={decorate,decoration={brace,mirror,amplitude=5pt}}
]
\node[mat,minimum width=6cm,minimum height=4.0cm] (IN) {};
\draw[colsep] ($(IN.north west)!0.375!(IN.north east)$) -- ($(IN.north west)!0.375!(IN.north east)+(0,-4cm)$ |- IN.north) -- cycle;
\draw[colsep] ($(IN.north west)!0.75!(IN.north east)$) -- ($(IN.north west)!0.75!(IN.north east)+(0,-4cm)$ |- IN.north) -- cycle;
\node[hdr] at ($(IN.north west)!0.18!(IN.north east)+(0,0.45)$) {Patch DoFs:\\velocity};
\node[hdr] at ($(IN.north west)!0.56!(IN.north east)+(0,0.45)$) {Patch DoFs:\\residual};
\node[hdr] at ($(IN.north west)!0.875!(IN.north east)+(0,0.45)$) {geometry};
\foreach \i/\y in {1/1.35, 2/0.65, 3/-0.15, 4/-0.95} {
 \node[rowlab] at ($(IN.west)+(-0.25,\y)$) {$M_{\i}$};
 \node[hdr] at ($(IN.west)!0.18!(IN.east)+(0,\y)$) {$\mathcal Q(\vt_\text{fine};M_{\i})$};
 \node[hdr] at ($(IN.west)!0.56!(IN.east)+(0,\y)$) {$\mathcal Q(\rt_\text{fine};M_{\i})$};
 \node[hdr] at ($(IN.west)!0.875!(IN.east)+(0,\y)$) {$\boldsymbol\omega(M_{\i})$};
}
\node[rowlab] at ($(IN.west)+(-0.25,-1.65)$) {$\vdots$};
\node[rowlab] at ($(IN.west)!0.18!(IN.east)+(0,-1.65)$) {$\vdots$};
\node[rowlab] at ($(IN.west)!0.56!(IN.east)+(0,-1.65)$) {$\vdots$};
\node[rowlab] at ($(IN.west)!0.875!(IN.east)+(0,-1.65)$) {$\vdots$};
\node at ($(IN.south)+(0,-0.45)$) {Input: $N_{\mathrm{cells}}^{coarse}\times(2N_{\mathrm{dof}}^{M} + n_{\mathrm{geo}})$};

\node[align=center,draw,very thick,rounded corners=2pt,fill=blue!6,minimum
width=2.5cm,minimum height=1.1cm,right=0.5cm of IN] (NN) {$\mathcal
 N$\\\footnotesize per patch (row-wise)};

\draw[->,thick] (IN.east) -- (NN.west);

\node[mat,minimum width=1.5cm,minimum height=4.0cm,right=0.5cm of NN] (OUT) {};
\node[hdr] at ($(OUT.north)+(0,0.45)$) {Patch DoFs:\\corrections};
\node at ($(OUT.south)+(0,-0.45)$) {Output: $N_{\mathrm{cells}}^{coarse}\times N_{\mathrm{dof}}^{M}$};

\foreach \i/\y in {1/1.35, 2/0.65, 3/-0.15, 4/-0.95} {
 \node[rowlab] at ($(OUT.east)+(0.8,\y)$) {$M_{\i}$};
 \node[hdr] at ($(OUT.west)!0.5!(OUT.east)+(0,\y)$) {$\delta \vt_M^{\i}$};
}
\node[rowlab] at ($(OUT.west)+(-0.25,-1.65)$) {$\vdots$};
\node at ($(OUT.center)+(0,-1.65)$) {$\vdots$};
\draw[->,thick] (NN.east) -- (OUT.west);
\end{tikzpicture}
\caption{\label{fig:dnnmg_batched_io}\emph{Batched input matrix composed of per-patch (one patch = one row) velocity entries \(\mathcal Q(\vt_\text{fine})\), residual entries \(\mathcal Q(\rt_\text{fine})\), and geometric descriptors \(\boldsymbol\omega\); the network \(\mathcal N\) maps each row (patch) to a vector of patch-local corrections \(\delta\vt^{M}_{K}\).}}
\end{figure}
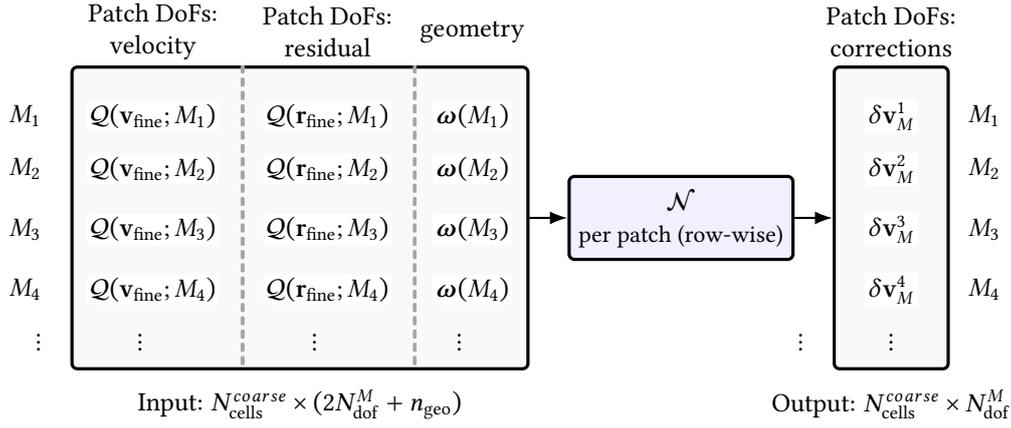
\paragraph{Network Input.}
Let $\vt_\text{fine}=\mathbf P\vt_\text{coarse}$ be the prolongated velocity and
let $\rt_\text{fine}$ be the fine-grid momentum residual. For each patch, we \emph{gather} the local
coefficients of $\vt_\text{fine}$ and $\rt_\text{fine}$ into a row and append
the geometric patch descriptors $\boldsymbol\omega$. Stacking all rows gives
\begin{equation}
 \label{eq:patch-wise-input}
 \mathbf X
 =
 \big[\mathcal Q(\vt_\text{fine})\ \big|\
     \mathcal Q(\rt_\text{fine})\ \big|\
     \boldsymbol\omega\big]
 \in
 \mathbb R^{N_{\mathrm{cells}}^\text{coarse}\times(2N_{\mathrm{dof}}^{M}+n_{\mathrm{geo}})}.
\end{equation}
Here the three blocks are, per row (i.e.\ per patch): local velocity DoFs,
local residual DoFs, and $n_{\mathrm{geo}}$ geometric features. For a given patch $M_{\mathcal{B}}$ with prediction jump $S$,
the flattened row dimension equals the network input size $N_{\mathrm{in}}=2N_{\mathrm{dof}}^{M}+n_{\mathrm{geo}}$, and the
patch correction has length $N_{\mathrm{out}}=N_{\mathrm{dof}}^{M}$; concrete values for
$(N_{\mathrm{in}},N_{\mathrm{out}})$ across $(M,S)$ are listed in
Tab.~\ref{tab:io-sizes-side}.

\paragraph{Patch-wise Prediction.}
The patch network $\mathcal N$ acts \emph{row-wise} on $\mathbf X$,
producing one correction vector per patch:
\[
 \mathbf D = \mathcal N(\mathbf X)
 \in \mathbb R^{N_{\mathrm{cells}}^\text{coarse}\times N_{\mathrm{dof}}^{M}},
\]
where the $K$-th row $\mathbf D_{K,:}$ is the patch-local velocity correction
$\delta\vt^{M_K}\in\mathbb R^{N_{\mathrm{dof}}^{M}}$ for $M_{K}$.

\paragraph{Scatter.}
Patch predictions must be merged into one fine-space update and made compatible
with Dirichlet boundary conditions. Using the gather matrices
$G_K\in\{0,1\}^{N_{\mathrm{dof}}^{M}\times N^\text{fine}_{dof}}$ from above, define the
\emph{scatter} matrices
$S_K\coloneq G_K^\top\in\{0,1\}^{N^\text{fine}_{dof}\times N_{\mathrm{dof}}^{M}}$, which add a
patch vector back into the global layout. Let
$W=\operatorname{diag}(\mu_1,\dots,\mu_{N_\text{fine}})$ be diagonal
\emph{partition-of-unity weights} chosen so that
$W\sum_K S_K G_K = I$ on interior DoFs, and let
$B=\operatorname{diag}(\beta_1,\dots,\beta_{N_\text{fine}})$ with $\beta_i=0$ on
Dirichlet DoFs and $\beta_i=1$ otherwise. The fine-grid increment
assembled from all patches is
\begin{equation}
 \label{eq:scatter-op}
 \delta\vt_\text{fine}
 =
 BW\sum_{K=1}^{N_{\mathrm{cells}}^\text{coarse}} S_K\delta\vt^{M}_K
 =
 BW\sum_{K=1}^{N_{\mathrm{cells}}^\text{coarse}} S_K\mathbf D_{K,:}^{\top}
 \in\mathbb R^{N_\text{fine}}.
\end{equation}
By construction, $\sum_{K} S_K G_K$ is diagonal with entries equal to the local
patch multiplicities; hence on interior DoFs, $W\sum_{K} S_K G_K = I$, and
$BW\sum_{K} S_K G_K = B$ overall. Thus, $W$ averages overlapping predictions of the network,
and $B$ enforces homogeneous Dirichlet values.

\subsection{Neural network}\label{sec:dnntrain}
We train the network $\mathcal N$ in a supervised setting with targets
from fine-space reference simulations, $\vt^\text{fine}_\text{ref}$, with converged Newton-Krylov iterations. For each of these snapshots we build, in
$\V_\text{fine}$, the patch-wise inputs (state, $\vt_\text{fine}$, residual, $\rt_\text{fine}$, geometry, $\omega$) and the corresponding
target, the \emph{corrections} to the coarse solution obtained from a fine reference solve. Concretely, the target on
the patch $M_{\mathcal B}$ is the local difference between the reference velocity and the
prolongated solution $\vt_\text{fine}=\mathbf P\vt_\text{coarse}$, restricted to that patch.

Training samples are the rows of the batched input matrix, collected over all
patches and all snapshots. Mini-batches are formed by randomly sampling rows
across snapshots. The loss is the mean-squared error over patches. We proceed to
a formal specification of the network architecture, targets and loss, together with the optimization
details.
\paragraph{Neural Network Architecture}
In the hybrid NN-FEM approach we use a compact multilayer perceptron (MLP)
that maps the patch-wise input vector
\(\mathbf x \in \mathbb R^{N_{\mathrm{in}}}\) to a local velocity correction
\(\mathcal N(\mathbf x) \in \mathbb R^{N_{\mathrm{out}}}\).
The network has \(\ell\) layers and is defined as follows.
We set
\[
 \mathbf h^{(0)} = \mathbf x\quad\text{and}\quad
 \mathbf h^{(i)} =
 \sigma\bigl(\mathbf W_i \mathbf h^{(i-1)} + \mathbf b_i\bigr),\quad i = 1, \dots, \ell-1\,.
\]
Here \(\mathbf W_i \in \mathbb R^{w \times d_{i-1}}\),
\(\mathbf b_i \in \mathbb R^{w}\), \(d_0 = N_{\mathrm{in}}\),
and \(d_i = w\) for \(i = 1,\dots,\ell-1\).
The output layer is
\[
 \mathcal N(\mathbf x)
 =\mathbf W_\ell \mathbf h^{(\ell-1)} + \mathbf b_\ell,
 \qquad
 \mathbf W_\ell \in \mathbb R^{N_{\mathrm{out}} \times w},\;
 \mathbf b_\ell \in \mathbb R^{N_{\mathrm{out}}}.
\]

Where dimensions match, we use a residual connection
\[
 \mathbf h^{(i)} \leftarrow \mathbf h^{(i)} + \mathbf h^{(i-1)}.
\]
In our implementation we choose the activation function
\(\sigma = \tanh\) and apply layer normalization to the pre-activations,
\[
 \mathbf h^{(i)} =
 \sigma\bigl( \mathrm{LN}_i(W_i \mathbf h^{(i-1)} + b_i) \bigr),
\]
which we have found to yield the most stable results.

\paragraph{Targets and Loss.}
For a fine reference solution $\vt_\text{fine}^{\mathrm{ref}}$, which is a computed on a mesh with higher resolution without deep neural network, and patch $M_K$, we define the
per-patch target increment
$\mathbf y^K_{M} \coloneq G_K\big(\vt_\text{fine}^{\mathrm{ref}}-\mathcal
P\vt_\text{coarse}\big)\in\mathbb R^{N_{\mathrm{dof}}^{M}}$. Let
$\delta\vt^{K}_M\in\mathbb R^{N_{\mathrm{dof}}^{M}}$ denote the network's
prediction for patch $M_K$. A simple mean-squared error with weight decay is
\begin{equation}
 \label{eq:simple-loss}
 \mathcal L(\theta)
 =
 \frac{1}{|\mathcal S|}\sum_{s\in\mathcal S}
 \frac{1}{N_{\mathrm{cells}}^{(s)}}\sum_{K=1}^{N_{\mathrm{cells}}^{(s)}}
 \big\|\delta\vt^{K,(s)}_M-\mathbf y^{K,(s)}_{M}\big\|_2^{2}
 +\alpha\|\theta\|_2^2,
\end{equation}
where $\mathcal S$ indexes training snapshots and $\alpha>0$ is the
Tikhonov/weight-decay parameter.
\paragraph{Optimization Details.}
For the training, we use the AdamW
optimizer~\cite{loshchilovDecoupledWeightDecay2018} (learning rate $10^{-4}$,
batch size $64$), which leads to better stability compared to alternatives, as
we observed in~\cite{margenberg2023hybrid}. We use Tikhonov regularization with
a scaling factor \(\alpha=10^{-3}\). The training is run for 40 epochs. Then,
the model with the lowest validation loss is selected. For stable training, we
use a warm-up with a low initial learning rate and a cyclical scheduler
performing a single cycle over all epochs~\cite{Smith2017}. Network parameters
are initialized with PyTorch's default initialization.

We conclude this section with a few remarks on the network design and
implementation aspects of the neural network.
\begin{remark}[Impact of the correction on Newton's method]\label{rem:dnnmg-newton}
In the sense of inexact
Newton~\cite{DemboEisenstatSteihaug1982,EisenstatWalker1996}, the NN-FEM
correction $\delta\vt_\text{fine}^n$ at $t^{n}$ induces an \emph{additive right-hand-side
 perturbation} of the residual $\mathcal R_{\mathrm{coarse}}^{n+1}$ at $t^{n+1}$, with the same Jacobian since
the perturbation depends on $\vt_\mathrm{fine}^{n}$ and is independent of $\vt_\text{coarse}^{n+1}$. Thus, the local
assumptions for inexact Newton (invertibility and Lipschitz continuity of
$\mathcal J$) can be transferred, and the theory can be applied here.
Admissibility then reduces to a \emph{relative smallness} of the additive load at the
initial guess, as detailed in Remark~\ref{rem:dnnmg-stab}. In particular, the local
Newton convergence region is preserved up to a small factor. A refined treatment and detailed theoretical
investigation is left to future work. Next, we explain how we design the neural
network such that the \emph{relative smallness} assumption can be fulfilled (cf.~\eqref{eq:smallness-assumption-i}). We note
that this reasoning is rather heuristic and will be investigated empirically
through ablation studies in Sec.~\ref{sec:study-ablation}.
\end{remark}
\begin{remark}[Neural network design for stability]\label{rem:nn-design}
To ensure the convergence of the nonlinear solver, the perturbation
of the fine-level right-hand side $\mathbf f_{\text{fine}}^{n+1}$, induced by
$\delta \mathbf v_{\text{fine}}^n$, must be bounded in terms of the unperturbed
residual (cf.~Remark~\ref{rem:dnnmg-stab}). We therefore hypothesize that the neural network should exhibit two key properties:
 \emph{state-awareness}, $\vt_\text{fine}$, and \emph{residual-awareness}, $\rt_\text{fine}$. We achieve this by
 explicitly including both, the current state and residual information, as inputs
 to the network.

 State-awareness requires that the correction
 $\delta\vt_\text{fine}$ aligns with the local solution structure.
 State-awareness is provided by feeding local patches of the velocity field
 $\mathbf P\vt_\mathrm{coarse}$ to the network. This allows the correction
 $\delta\vt_\text{fine}$ to align with the local solution geometry and stress
 state, ensuring consistency with the Jacobian structure
 $\mathcal J(\vt_\mathrm{coarse})$ and preventing the introduction of
 spurious modes.

 Residual-awareness requires that the correction magnitude is appropriately
 scaled relative to the local discretization error, thereby controlling
 the perturbation of $\mathbf f_\mathrm{coarse}$ through the corrected solution $\vt_\mathrm{fine}$
 to satisfy the bound above.
 Residual-awareness is provided by including the local residual
 $\mathcal R(\vt_\mathrm{coarse};\cdot)$ as an additional input feature (cf.~\eqref{eq:Newvar}).
\end{remark}
\begin{remark}[Parallelism, batching, and conservation]
 Gather $\mathcal Q$, prediction $\mathcal N$, and
 scatter operations are embarrassingly parallel across patches:
 patch-wise predictions are mutually independent up to the scatter
 step. In practice, all patch-local data are stacked into a single input array
 so that $\mathcal N$ can be evaluated in one batched call on modern
 accelerators, yielding higher arithmetic intensity than the typically
 memory-bound finite element kernels and thus higher throughput. Boundary
 values are preserved by zeroing Dirichlet DoFs via the scatter mask $B$, and
 partition-of-unity weights $W$ stabilize the assembly by averaging
 contributions at shared DoFs in~\eqref{eq:scatter-op}. Consequently, the global correction avoids
 spurious jumps at patch interfaces.
\end{remark}

\section{Numerical Experiments}
In this section, we describe the numerical experiments performed to generate
training data, train neural networks, and test the NN-FEM algorithm for sea-ice
simulations. The experiments involve several parameter and ablation studies. The
experiments ran on a workstation with two Intel Xeon E5-2640 v4 CPUs, 756 GB
RAM, and an NVIDIA Tesla V100 GPU.
The numerical solver, which is part of Gascoigne3D~\cite{beckerFiniteElementToolkit}, runs on the CPU, while the network is evaluated on the GPU using the libtorch library~\cite{paszkePyTorchImperativeStyle2019}. The neural network architecture is a simple MLP with recurrent connections.
While more advanced architectures are supported, the MLP performed reliably in
previous works~\cite{margenberg2023hybrid}. In Section \ref{sec:bench} we introduce the benchmark setup.  Section \ref{sec:nn-par-and-train} outlines the training data and introduces the metrics for evaluation of the benchmark problem. In the experiments we conduct in Section \ref{sec:study-patchsize-jumplevels}, we
investigate how key parameters of the algorithm, namely the patch size, $N_M$, the
number of predicted levels, $S$, and the choice of network inputs, affect accuracy
and generalization. Rather than optimizing the network architecture, we conduct
a systematic study to understand the influence of the aforementioned parameters
on the performance of the NN-FEM approach.

\subsection{Benchmark Setup}\label{sec:bench}
We solve the sea-ice dynamics~\eqref{eq:momentum}--\eqref{eq:thickness} for the benchmark described in~\cite{Mehlmann2021}.
We consider a square domain $\Omega = (0, \SI{512}{\kilo\meter})^2$, initially covered by a thin, motionless layer of sea-ice. The initial conditions are
\begin{equation*}
\vt(\xt,0) = \mathbf{0},\SI{}{\meter\per\second}, \qquad
A(\xt,0) = 1, \qquad
H(\xt,0) = \SI{0.3}{\meter}.
\end{equation*}
A cyclone propagates from the center of the domain toward the upper-right corner (northwest direction), while the ocean current is prescribed to circulate counterclockwise.
The corresponding wind and ocean velocities in~\eqref{eq:forcing} are specified as
\begin{equation*}
\begin{aligned}
  &\vt_{\mathrm w}(\xt,\,t) = \SI{0.01}{\meter\per\second}
  \begin{pmatrix}
    -1 + 2y/ (\SI{512}{\kilo\meter})\\
    \phantom{-}1 - 2x/ (\SI{512}{\kilo\meter})
  \end{pmatrix},\\
  &\vt_{\mathrm a}(\xt,\,t) =
  \omega(\xt)\,\bar{\vt}_{\mathrm a}^{\max}
  \begin{pmatrix}
    \phantom{-}\cos(\alpha) &\sin(\alpha) \\
    -\sin(\alpha) &\cos(\alpha)
  \end{pmatrix}
  \begin{pmatrix}
    x - m_x(t) \\
    y - m_y(t)
  \end{pmatrix},
\end{aligned}
\end{equation*}
where
\begin{equation*}
\begin{aligned}
  \bar{\vt}_{\mathrm a}^{\max} = \bar{\vt}_{\mathrm a}^{\max}(t) &= \SI{15}{\meter\per\second}
  \begin{cases}
   -\tanh ((\phantom{1}4-t)(\phantom{-}4+t)/2) & t\in [\SI{0}{\day},\SI{4}{\day}]\\
   \phantom{-}\tanh ((12-t)(-4+t)/2) & t\in [\SI{4}{\day},\SI{8}{\day}]\\
  \end{cases},\\
  \alpha = \alpha(t) &= \begin{cases}
   72^\circ\quad t\in [\SI{0}{\day},\,\SI{4}{\day}]\\
   81^\circ\quad t\in [\SI{4}{\day},\,\SI{8}{\day}]\\
  \end{cases},\\
  m_x(t) = m_y(t) &=
  \begin{cases}
   \SI{256}{\kilo\meter} + \SI{51.2}{\kilo\meter\per\day}\, t & t\in [\SI{0}{\day},\SI{4}{\day}] \\
   \SI{665.6}{\kilo\meter} - \SI{51.2}{\kilo\meter\per\day}\, t & t\in [\SI{4}{\day},\SI{8}{\day}]\\
  \end{cases}\\
  \omega(\xt)=\omega(\xt,\,t) &= \frac{1}{50}\exp\left(-\frac{\sqrt{(x - m_x(t))^2 + (y - m_y(t))^2}}{\SI{100}{\kilo\meter}}\right).
\end{aligned}
\end{equation*}
\begin{figure}
 \centering
 \subcaptionbox{}{
 \includegraphics[height=.275\linewidth]{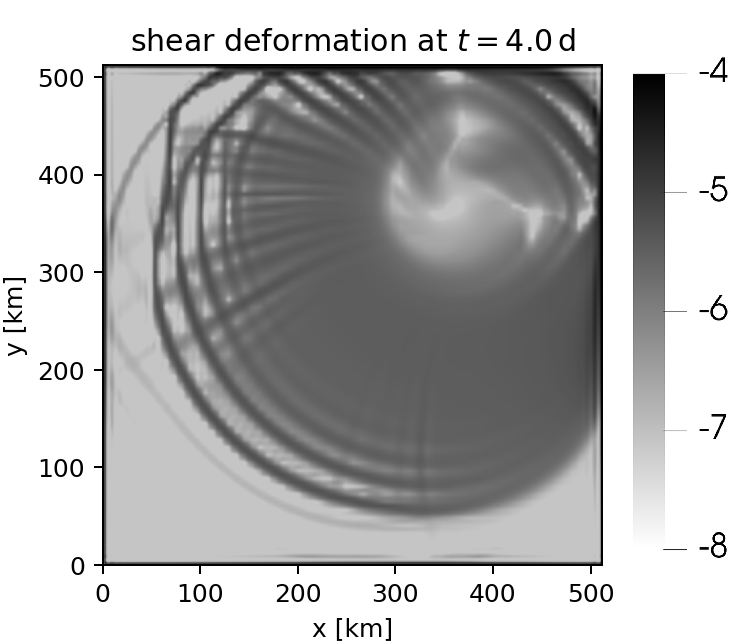}}
 \subcaptionbox{}{
 \includegraphics[height=.275\linewidth]{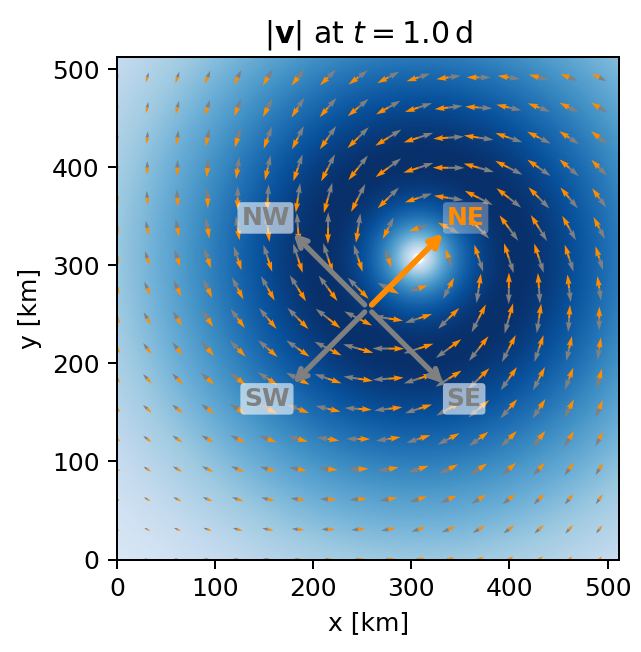}}
 \subcaptionbox{}{
 \includegraphics[height=.275\linewidth]{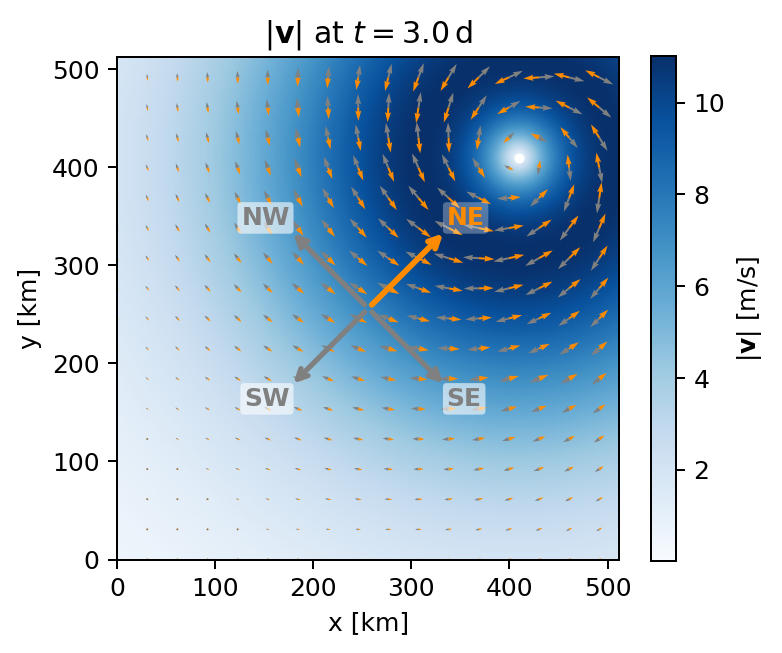}}
 \caption{\label{fig:bench-wind}\emph{
 Test and training setups. For training the model we consider the sea-ice benchmark problem defined in~\cite{Mehlmann2021}, where a cyclone is moving in northeast direction over an ice-covered domain. Panel (a) shows sea-ice shear deformation after 4 days of simulation. Panel (b) and (c) indicates
 the wind field of the cyclone $\vt_{\mathrm w}$ at
  $t=\SI{1}{\day}$ and $t=\SI{3}{\day}$, respectively. %
  For testing the network the wind field are either anti-cyclone or cyclone moving in north-west (NW),
  south-west (SW) or south-east (SE) direction. The different directions are indicated in gray in panel b and c. } }
\end{figure}
\subsection{\label{sec:nn-par-and-train}Neural Network Parameters, Training and Metrics for Evaluation}
For training and testing, we vary configuration parameters that, on the one hand, modify the benchmark problem itself, such as
\begin{enumerate}
\item \emph{Wind-field direction:} north-west, north-east, south-east, or south-west;
\item \emph{Direction of rotation:} cyclone or anticyclone.
\end{enumerate}
A visualization of these configurations is provided in Fig.~\ref{fig:bench-wind}.
{We note that, for the rotated wind forcing, the predicted fields align with the north-east
reference up to rotation, indicating approximate rotational equivariance despite
the local, patch-wise design.} On the other hand, we vary parameters that affect the network input and output dimensions, including
\begin{enumerate}
\item[3.] \emph{Auxiliary refinement depth $S$:} the number of uniform refinements used to define the auxiliary space $V_\text{fine}$ from the working space $V_\text{coarse}$ (cf.~Sec.~\ref{sec:dnnmg}); (a larger $S$ results in a finer $V_\text{fine}$);
\item[4.] \emph{Patch size $N_M$:} the number of uniform coarsening steps used
 to construct the patch from the working mesh.
\end{enumerate}
We collect these parameters in Tab.~\ref{tab:nn-parameters} and list the
corresponding input and output sizes $(N_{\mathrm{in}},N_{\mathrm{out}})$ per patch size $N_M$ and mesh refinement level $S$
in Tab.~\ref{tab:io-sizes-side}. We use the benchmark setting described in Section~\ref{sec:bench} to generate training data separately for each
pair $(N_M,S)$. We neglect the first 30 time steps and obtain 66 snapshots. The
training data are split into two sets: 75\% for training and 25\% for testing.
The remaining settings (other wind directions and anticyclones) are used for
validation.

The neural network architecture is a simple MLP with recurrent connections.
While more advanced architectures are supported, the MLP performed reliably in
previous works~\cite{margenberg2023hybrid}. Further, the focus of
this work is to formulate a hybrid NN-FE method, based on
DNN-MG~\cite{margenberg2023hybrid,margenberg2022nnmg,margenberg2021structure},
for coupled problems. In the experiments we conduct in this section, we
investigate how key parameters of the algorithm, namely the patch size, the
number of predicted levels, and the choice of network inputs, affect accuracy
and generalization. Rather than optimizing the network architecture, we conduct
a systematic study to understand the influence of the aforementioned parameters
on the performance of the NN-FEM approach.\@

\paragraph{Evaluation Metrics.}
Unless stated otherwise, all quality metrics are computed in the auxiliary fine
space \(\V_\text{fine}\). Given the corrected field \(\vt_\text{fine}\) and the fine
reference \(\vt_\text{fine}^{\mathrm{ref}}\),
\begin{enumerate}[label={(M\arabic*)},ref={(M\arabic*)},leftmargin=*,itemsep=0pt,topsep=0pt,parsep=0pt]
 \item\label{itm:metric-L2}\emph{Discrete \(\ell^2\) error:}
  \(E_{\ell^2} \coloneq \big\|\vt_\text{fine}-\vt_\text{fine}^{\mathrm{ref}}\big\|
  = {\big({(\vt_\text{fine}-\vt_\text{fine}^{\mathrm{ref}})}^\top(\vt_\text{fine}-\vt_\text{fine}^{\mathrm{ref}})\big)}^{1/2}\).
  Here $\vt_\text{fine}^{\mathrm{ref}}$ denotes the solution at
  \SI{2}{\kilo\meter} resolution.
 \item\label{itm:metric-LKF}\emph{LKF count:} We detect linear kinematic
  features (LKFs) on the coarse mesh as described in~\cite{Hutter.2019.code}.
  The detection algorithm is based on the shear deformation,
  \[
   \dot\epsilont_{\mathrm{shear}}
   = \sqrt{(\dot \epsilont_{11} - \dot \epsilont_{22})^2 + 4\,\dot\epsilont_{12}^2},
  \]
  where $\epsilont_{ij}$ denote the components of the symmetric strain-rate
  tensor \eqref{eq:strain}. The shear deformation, and thus the LKF detection, are computed on
  the coarse mesh. In accordance with the hybrid NN-FEM concept, only the
  neural network operates on the fine level, while the quantities of interest
  are evaluated on the coarse resolution.
  All parameters in the detection algorithm are kept constant across
  all experiments so that the number of LKFs, $N_{\mathrm{LKF}}$, is
  comparable between settings.
 \item\label{itm:metric-Newton}\emph{Solver effort:} total Newton iterations
 \(N_{\mathrm{Newton}} \coloneq \sum_{n} \kappa_n\), with \(\kappa_n\) the iterations in the momentum
 solve at time step \(n\).
 \item\label{itm:metric-time}\emph{Wall time:} Solution time for the momentum equation, \(T_{\mathrm{mom}}\), and its decomposition
 into Newton-Krylov time \(T_{\mathrm{NK}}\), CPU to GPU copy time
 \(T_{\mathrm{cpy}}\), and network inference time \(T_{\mathrm{nn}}\).
\end{enumerate}

Unless stated otherwise, the metrics~\ref{itm:metric-L2}--\ref{itm:metric-time}
are reported as arithmetic means over all test cases, all wind directions and
both cyclones and anticyclones. We do not average across jump levels
\(S\) or patch sizes \(N_M\); results for each \(S\) and \(N_M\) are always reported
separately. Accuracy is quantified by \(E_{\ell^2}\) and \(N_{\mathrm{LKF}}\);
efficiency by \(N_{\mathrm{Newton}}\) and \(T_{\mathrm{mom}}\).
\begin{table}[htb]
 \centering
 \caption{\label{tab:nn-parameters}\emph{Parameters of the hybrid neural network-finite element method. We
  use a multilayer perceptron with width \(w\) and number of network layers (depth)
  \(l\), and test different patch configurations by varying the additional refinements of the mesh
  \(S\) and patch size \(N_M\).}}
 \centering
 \setlength{\tabcolsep}{12pt}
 \renewcommand{\arraystretch}{1.15}
 \begin{tabular}{cccc}
 \toprule
\textbf{Width $w$} & \textbf{Layers $l$} & \textbf{Jump level $S$ (finer $V_\text{fine}$)} & \textbf{Patch size $N_M$} \\
 \midrule
 $\{128,\,256,\,512\}$ & $\{2,\,4,\,8\}$ & $\{1,\,2\}$ & $\{0,\,1,\,2\}$ \\
 \bottomrule
 \end{tabular}
\end{table}
\begin{table}[t]
 \centering
\caption{\label{tab:io-sizes-side} \emph{
Input and output sizes $N_{\mathrm{in}}$, $N_{\mathrm{out}}$ of the
  neural network for different patch sizes $N_M$ and additional mesh refinements $S$ (the larger $S$ the more degree of freedom contains 
  $\V_\text{fine}$).
  Here, $n_M$ denotes the number of DoFs per component and per patch on the auxiliary (refined) level. The input dimension is
  \(
  N_{\mathrm{in}} = 4n_M + 8,
  \)
  corresponding to 2 momentum components, 2 residual components, and 8 geometric features.
  The output dimension is \(
  N_{\mathrm{out}} = 2n_M,
  \) corresponding to the 2 momentum components.}}
 \setlength{\tabcolsep}{6pt}
 \begin{tabular}{c r r r r r r}
  \toprule
  & \multicolumn{3}{c}{$S=1$} & \multicolumn{3}{c}{$S=2$} \\
  \cmidrule(lr){2-4}\cmidrule(lr){5-7}
  $N_M$ & $n_M$ & $N_{\mathrm{in}}$ & $N_{\mathrm{out}}$ & $n_M$ & $N_{\mathrm{in}}$ & $N_{\mathrm{out}}$ \\
  \midrule
  0 &  25 &  108 &  50 &  81 &  332 &  162 \\
  1 &  81 &  332 &  162 &  289 &  1164 &  578 \\
  2 & 289 & 1164 &  578 & 1089 &  4364 &  2178 \\
  \bottomrule
 \end{tabular}
\end{table}

\subsection{Impact of Patch Size, Predicted Levels and Network
 Parameters}\label{sec:study-patchsize-jumplevels}
\begin{figure}[htbp]
 \sisetup{
  scientific-notation = true,
  round-mode = places,
  round-precision = 0
 }
 \centering
 \captionsetup[sub]{justification=centering}
  \textbf{Reference simulations}\\
 \begin{subfigure}[t]{0.32\linewidth}
  \includegraphics[width=\linewidth]{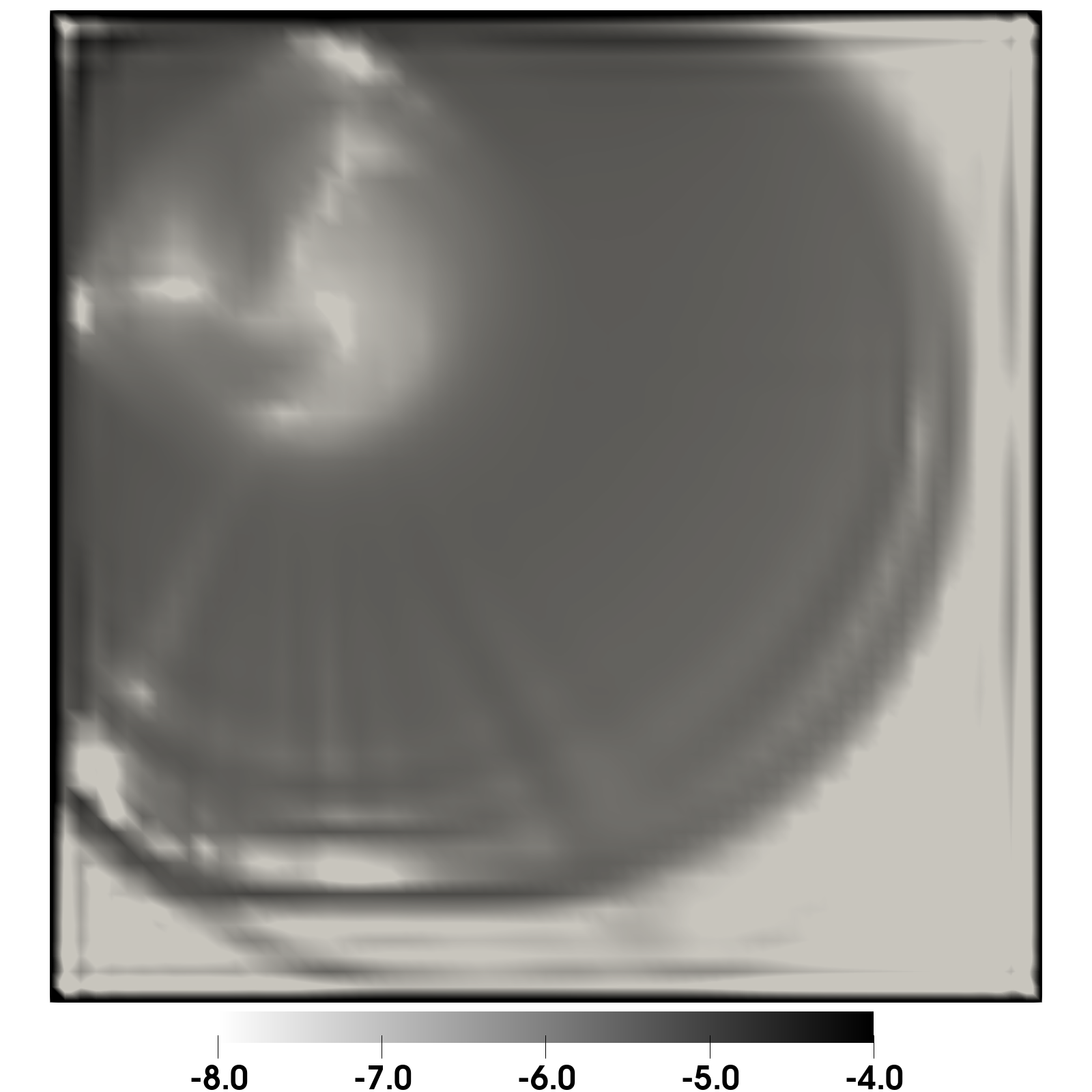}
  \subcaption{\SI{8}{\kilo\meter} resolution ($\V_\text{coarse}$, base)}\label{fig:grid-8k-base}
 \end{subfigure}\hfill
 \begin{subfigure}[t]{0.32\linewidth}
  \includegraphics[width=\linewidth]{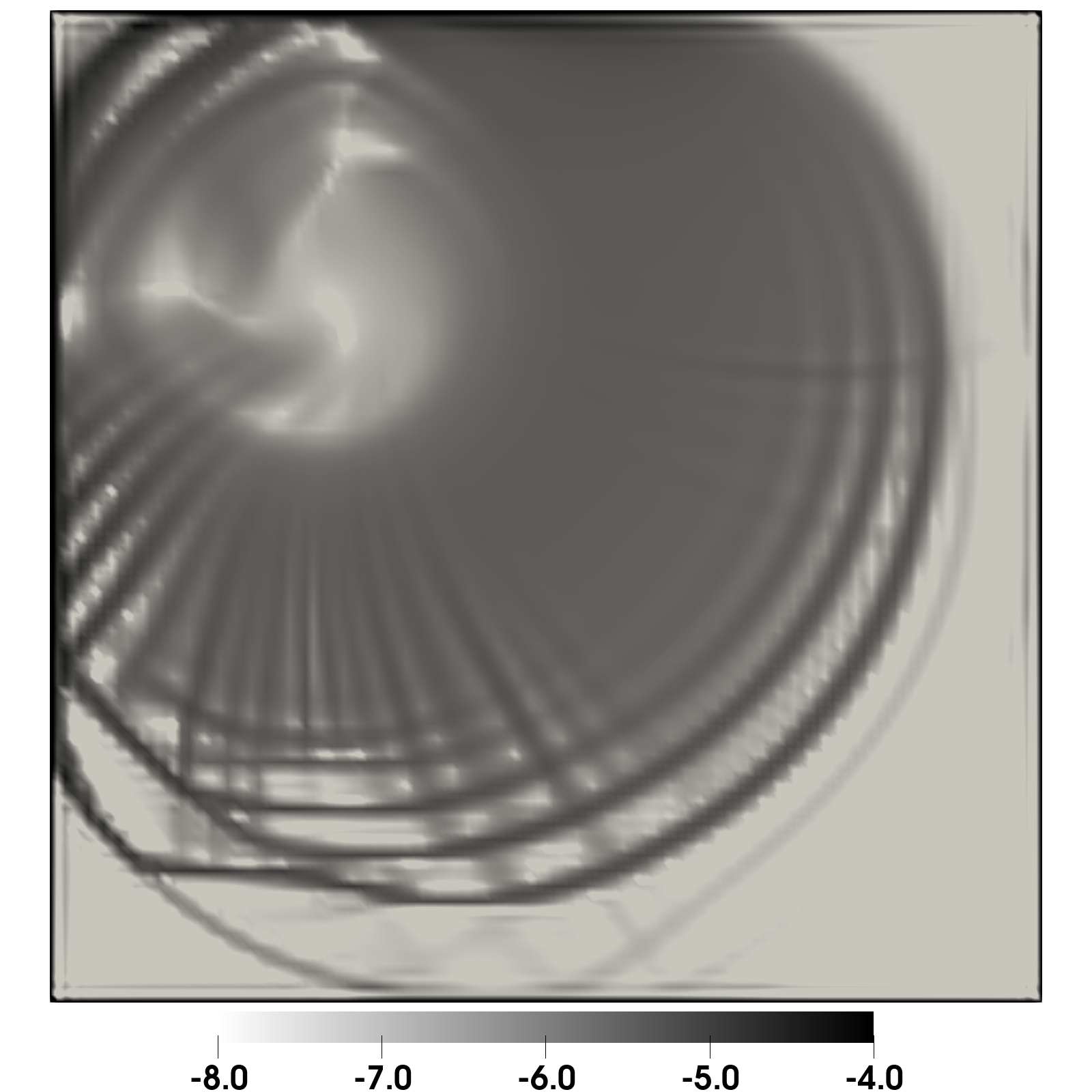}
  \subcaption{\SI{4}{\kilo\meter} resolution ($V_\text{fine}$, $S=1$)}\label{fig:grid-4k-ref}
 \end{subfigure}\hfill
 \begin{subfigure}[t]{0.32\linewidth}
  \includegraphics[width=\linewidth]{snapshot-final-r7-none-nw.png}
  \subcaption{\SI{2}{\kilo\meter} resolution ($V_\text{fine}$, $S=2$)}\label{fig:grid-2k-ref}
 \end{subfigure}\\[1ex]
 \textbf{NN-FEM with coarse base \SI{8}{\kilo\meter} trained on \SI{4}{\kilo\meter}}\\
 \begin{subfigure}[t]{0.32\linewidth}
  \includegraphics[width=\linewidth]{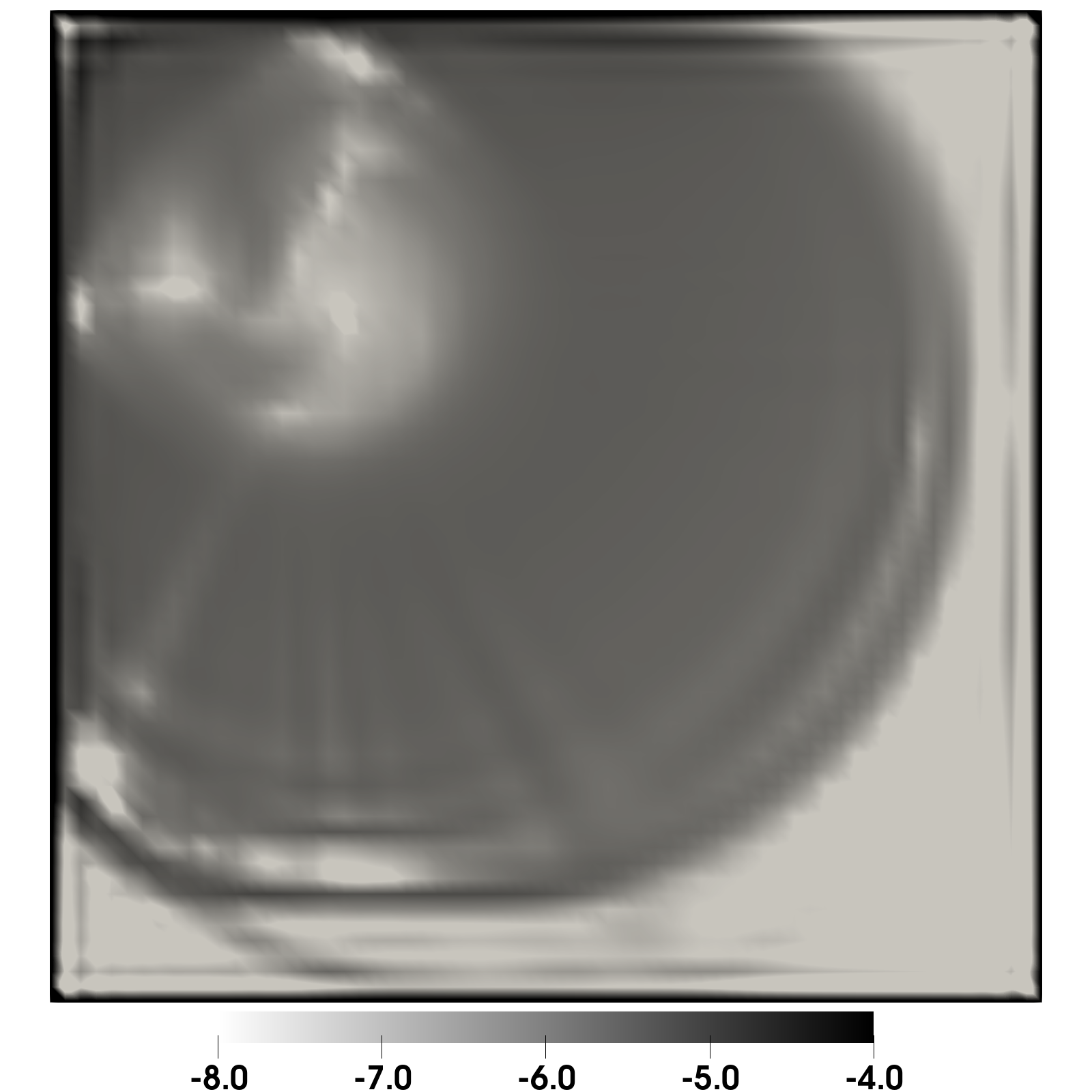}
  \subcaption{NN-FEM, $S=1$, $N_M=0$}\label{fig:grid-j1-p1}
 \end{subfigure}\hfill
 \begin{subfigure}[t]{0.32\linewidth}
  \includegraphics[width=\linewidth]{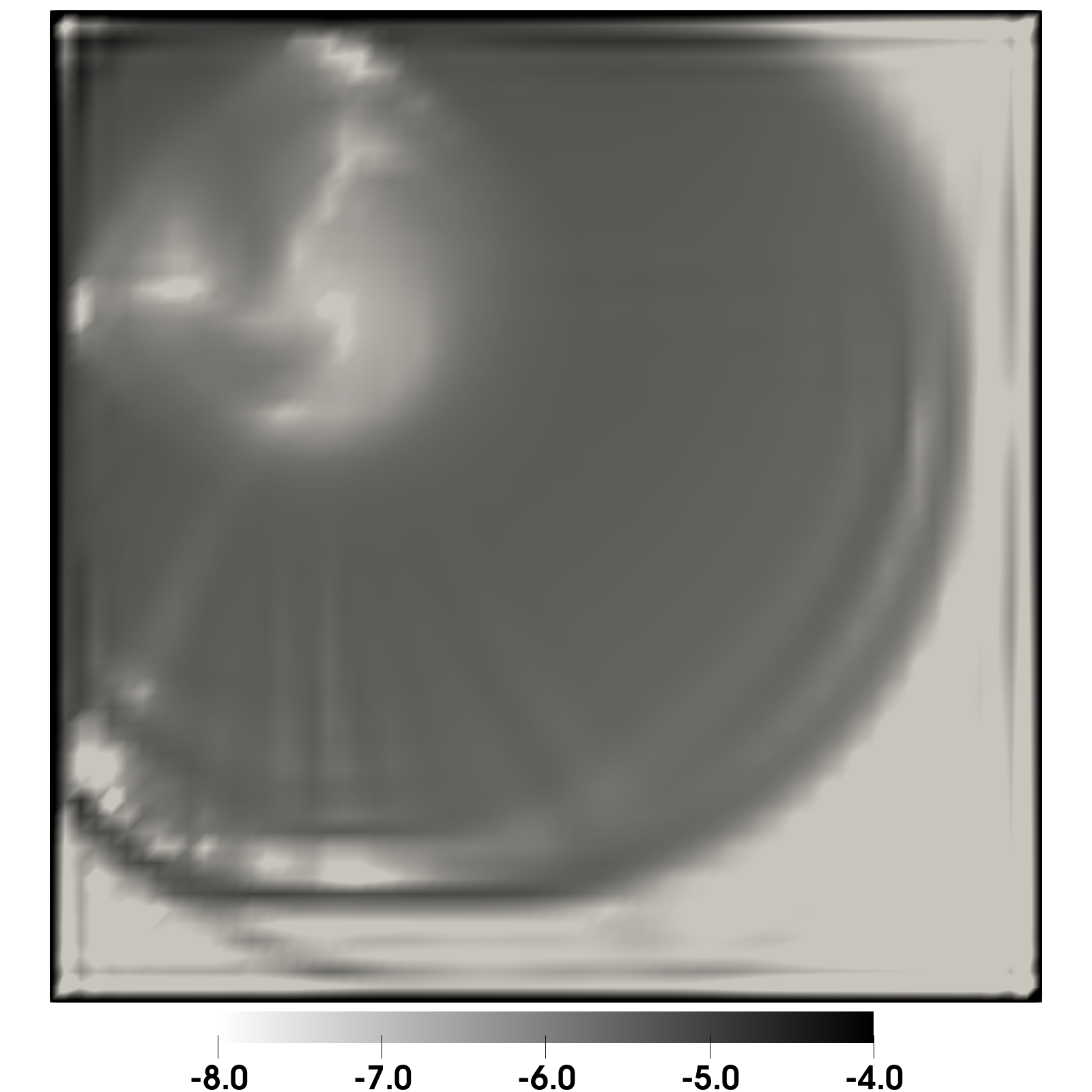}
  \subcaption{NN-FEM, $S=1$, $N_M=1$}\label{fig:grid-j1-p2}
 \end{subfigure}\hfill
 \begin{subfigure}[t]{0.32\linewidth}
  \includegraphics[width=\linewidth]{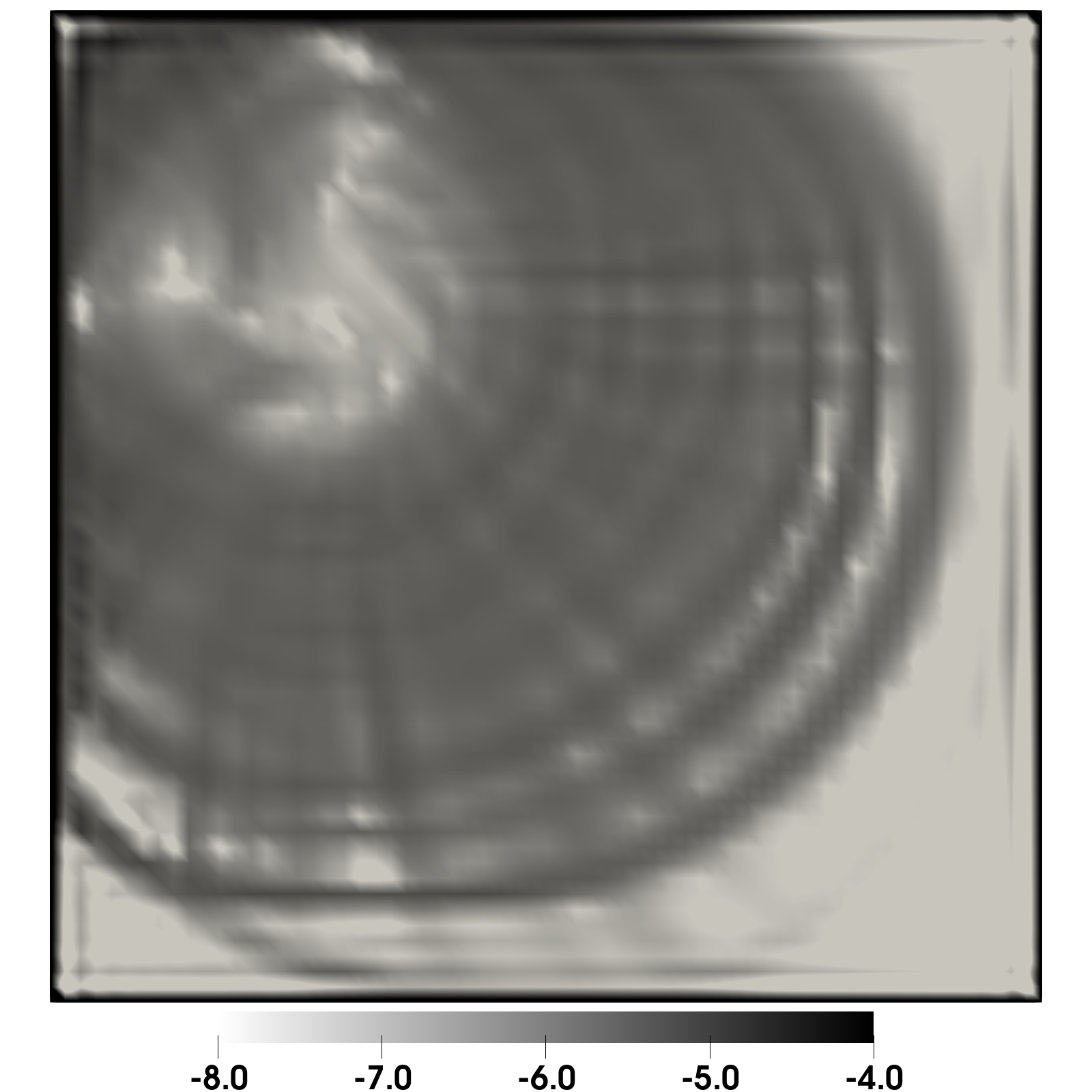}
  \subcaption{NN-FEM, $S=1$, $N_M=2$}\label{fig:grid-j1-p3}
 \end{subfigure}\\[1ex]
 \textbf{NN-FEM with coarse base \SI{8}{\kilo\meter} trained on \SI{2}{\kilo\meter}}\\
 \begin{subfigure}[t]{0.32\linewidth}
  \includegraphics[width=\linewidth]{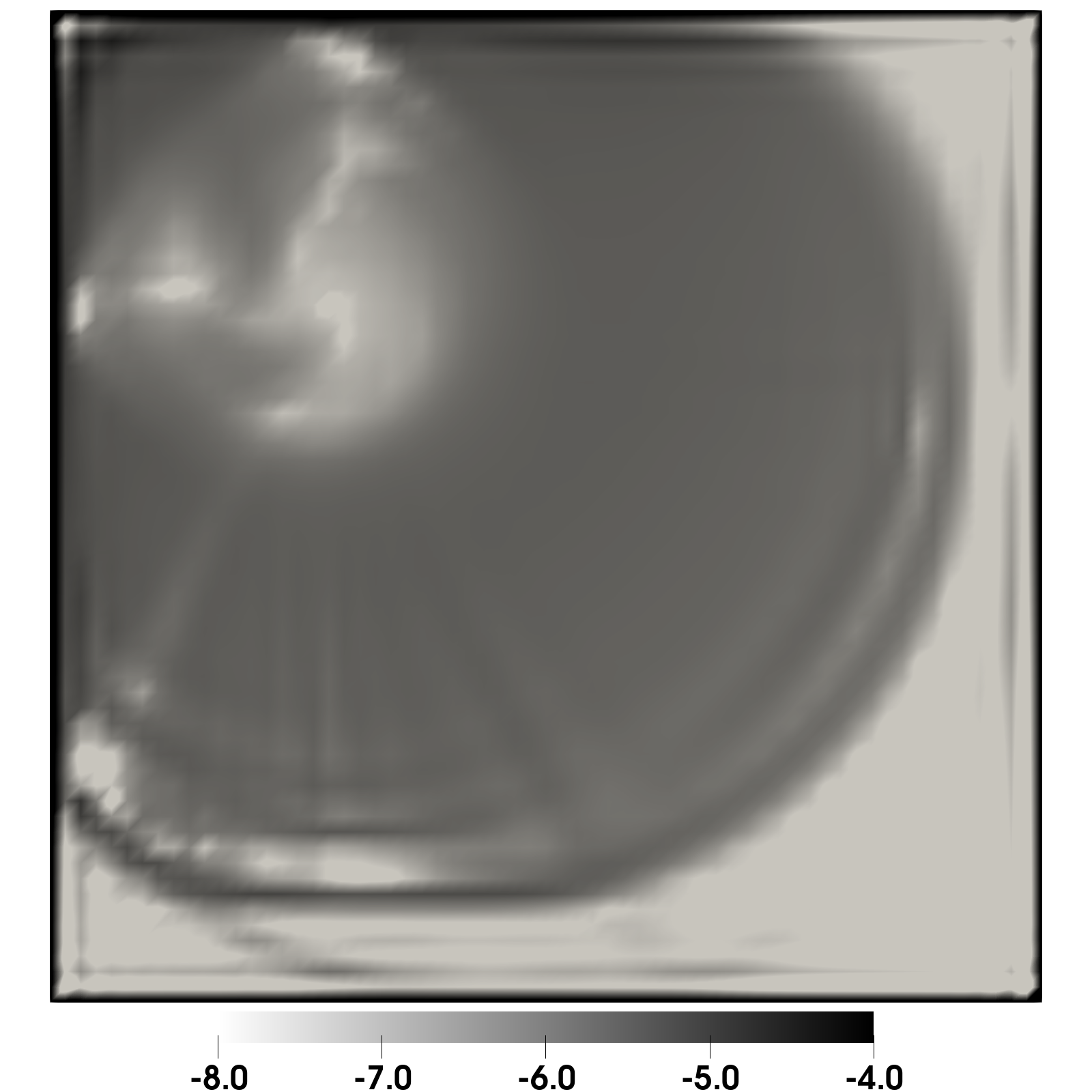}
  \subcaption{NN-FEM, $S=2$, $N_M=0$}\label{fig:grid-j2-p1}
 \end{subfigure}\hfill
 \begin{subfigure}[t]{0.32\linewidth}
  \includegraphics[width=\linewidth]{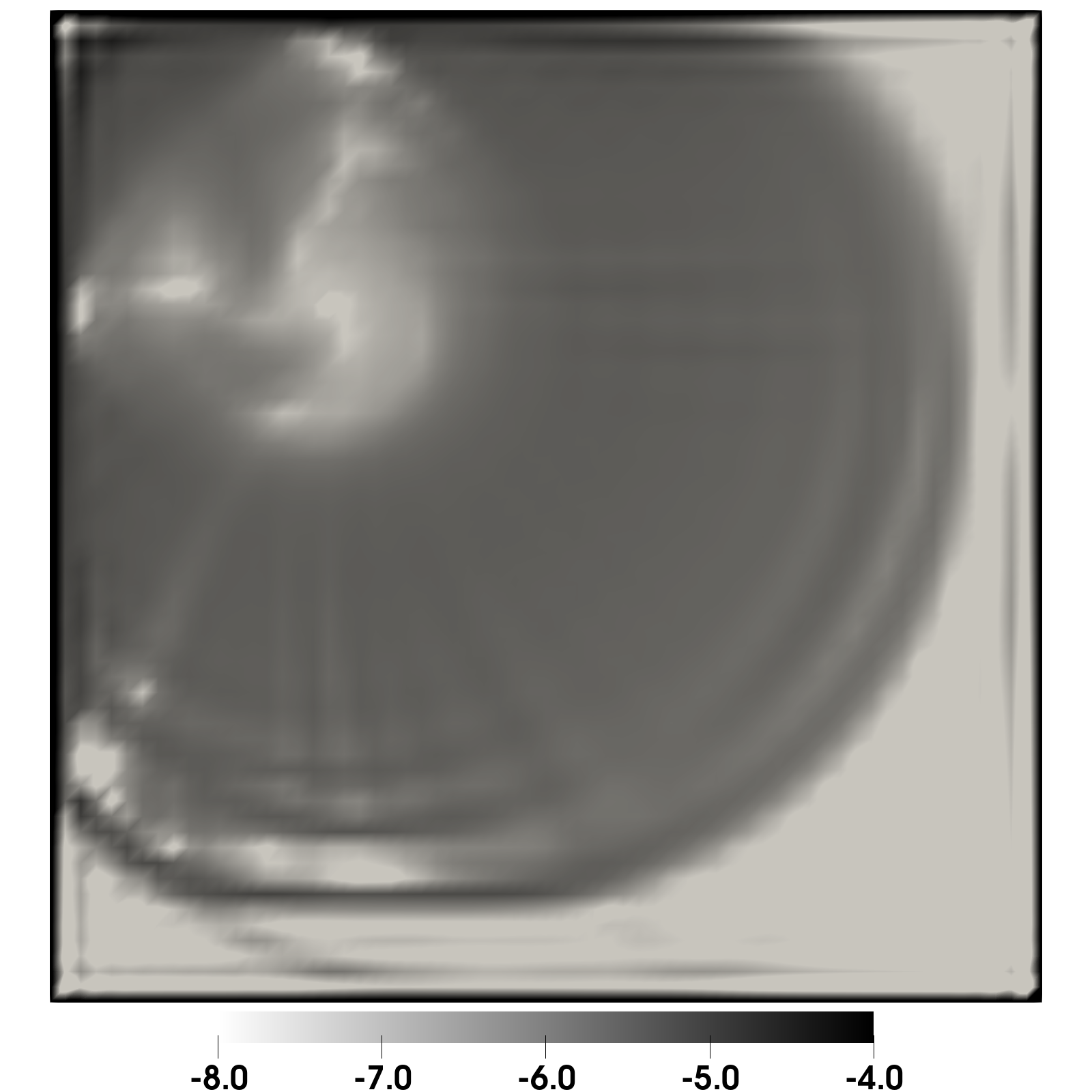}
  \subcaption{NN-FEM, $S=2$, $N_M=1$}\label{fig:grid-j2-p2}
 \end{subfigure}\hfill
 \begin{subfigure}[t]{0.32\linewidth}
  \includegraphics[width=\linewidth]{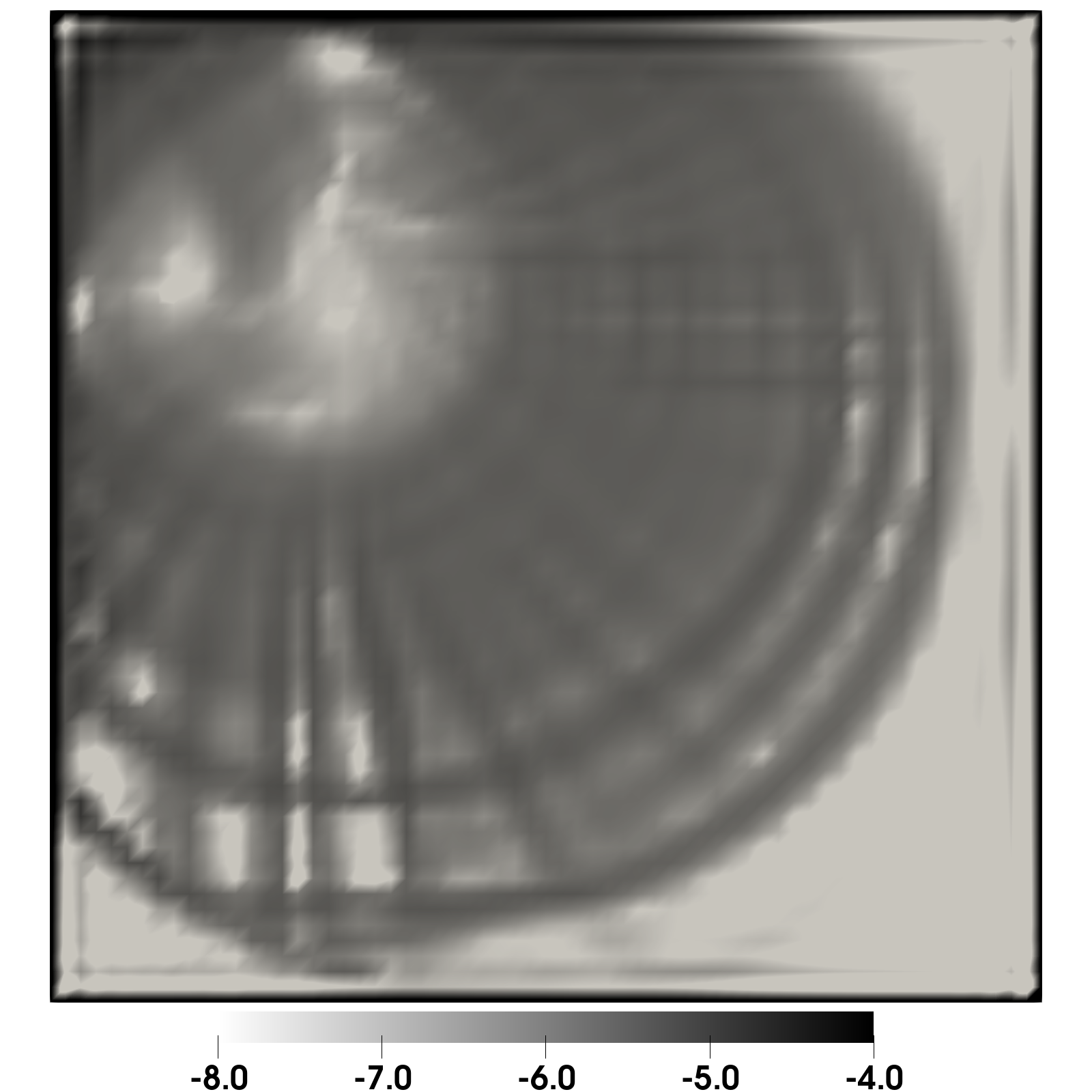}
  \subcaption{NN-FEM, $S=2$, $N_M=2$}\label{fig:grid-j2-p3}
 \end{subfigure}
 \caption{\label{fig:ss-3x3} \emph{Shear deformation simulated with the benchmark setup, where the direction of the cyclone is changed to north-west. The
  \emph{top row} shows the simulations performed without NN: (a) \SI{8}{\kilo\meter} coarse resolution simulation, (b) \SI{4}{\kilo\meter} reference simulation, the target for NN-FEM with $S=1$,
  (c) \SI{2}{\kilo\meter} reference simulation, the target for NN-FEM with $S=2$. The
  \emph{middle and bottom rows} show the simulated shear deformation derived from a NN-FEM simulation on $\V_\text{coarse}$ with a mesh size of \SI{8}{\kilo\meter}. The corresponding NN uses varying patch size $N_M$ and operates on $\V_\text{fine}$ with mesh sizes of \SI{4}{\kilo\meter}, i.e.\ $S=1$ (middle row) and \SI{2}{\kilo\meter}, i.e. $S=2$ (bottom row). }}
\end{figure}

\pgfplotstableread[row sep=\\, col sep=space]{
 l w val\\
 1 1 37 \\ 2 1 37 \\ 3 1 37 \\
 1 2 36 \\ 2 2 35 \\ 3 2 36 \\
 1 3 36 \\ 2 3 37 \\ 3 3 35 \\
}\lkfpsOneJ
\pgfplotstableread[row sep=\\, col sep=space]{
 l w val\\
 1 1 41 \\ 2 1 41 \\ 3 1 40 \\
 1 2 41 \\ 2 2 39 \\ 3 2 39 \\
 1 3 40 \\ 2 3 41 \\ 3 3 39 \\
}\lkfpsOne
\pgfplotstableread[row sep=\\, col sep=space]{
 l w val val2\\
 1 1 37 36 \\ 2 1 37 35 \\ 3 1 38 35 \\
 1 2 39 36 \\ 2 2 37 34 \\ 3 2 41 35 \\
 1 3 20 22 \\ 2 3 48 36 \\ 3 3 40 37 \\
}\lkfpsTwo
\pgfplotstableread[row sep=\\, col sep=space]{
 l w val val2\\
 1 1 37 41 \\ 2 1 43 40 \\ 3 1 60 41 \\
 1 2 49 42 \\ 2 2 42 40 \\ 3 2 44 39 \\
 1 3 44 42 \\ 2 3 43 41 \\ 3 3 56 42 \\
}\lkfpsThree
\pgfplotstableread[row sep=\\, col sep=space]{
 l w val val2\\
 1 1 65 57 \\ 2 1 63 54 \\ 3 1 73 56 \\
 1 2 53 49 \\ 2 2 64 55 \\ 3 2 70 61 \\
 1 3 59 55 \\ 2 3 68 58 \\ 3 3 71 60 \\
}\lkfpsFour
\begin{table}
 \caption{\label{tab:lkf-number} \emph{LKF counts evaluated from simulations on the 8
  km mesh corrected with the NN-FEM approach using patch sizes $N_M\in \{0, 1,
   2\}$ and trained on refined meshes corresponding to 4 km and 2 km
  resolution ($S \in \{1, 2\}$, respectively). Increasing $N_M$ increases the
  number of resolved LKFs across architectures. Network depth and width
  ($l,w$) shifts absolute counts only slightly compared to a change in the
  patch size.
  Solving the sea-ice momentum without NN results in an LKF count of 18,
  51, 171,
  at \SI{8}{\kilo\meter}, \SI{4}{\kilo\meter} and \SI{2}{\kilo\meter} resolution, respectively.
 }}
\begin{center}
 \begin{tikzpicture}[font=\small]
\begin{groupplot}[
 group style={group size=3 by 2, horizontal sep=1.5cm, vertical sep=1.5cm,group name=myplots},
 width=0.3\linewidth,height=3cm,
 xlabel={$l$}, ylabel={$w$},
 colormap/viridis, shader=flat,
 table/trim cells=true,
 xticklabel pos=right,
 y dir=reverse,
 xtick={1,2,3}, ytick={1,2,3},
 xticklabels={2,4,8},
 yticklabels={128,256,512},
 xmin=0.5, xmax=3.5, ymin=0.5, ymax=3.5,
 point meta min=20, point meta max=60,
 colorbar=false,
 title style={yshift=2ex},
 every axis x label/.style={at={(ticklabel* cs:0)},anchor=near ticklabel,yshift=0.25ex,xshift=-0.5ex},
 every axis y label/.style={at={(ticklabel* cs:1)},anchor=near ticklabel,xshift=-0.5ex,yshift=0.5ex}
]
\nextgroupplot[title={$N_M=0$, $S=1$}]
\addplot [matrix plot*, point meta=explicit, mesh/cols=3,
nodes near coords={
 \pgfmathfloattofixed{\pgfplotspointmeta}%
 \edef\MetaFixed{\pgfmathresult}%
 \pgfmathparse{\MetaFixed-50}%
 \ifdim\pgfmathresult pt>0pt \color{black}\else \color{white}\fi
 \pgfmathprintnumber[fixed,precision=1]{\pgfplotspointmeta}
},
nodes near coords align={center},
every node near coord/.append style={font=\footnotesize\bfseries}
] table[x=l, y=w, meta=val] {\lkfpsOneJ};
\coordinate (lu1) at (yticklabel* cs:1);
\nextgroupplot[title={$N_M=1$, $S=1$}]
\addplot [matrix plot*, point meta=explicit, mesh/cols=3,
 nodes near coords={
  \pgfmathfloattofixed{\pgfplotspointmeta}%
  \edef\MetaFixed{\pgfmathresult}%
  \pgfmathparse{\MetaFixed-50}%
  \ifdim\pgfmathresult pt>0pt \color{black}\else \color{white}\fi
  \pgfmathprintnumber[fixed,precision=1]{\pgfplotspointmeta}
 },
 nodes near coords align={center},
 every node near coord/.append style={font=\footnotesize\bfseries}
] table[x=l, y=w, meta=val2] {\lkfpsTwo};
\coordinate (lu3) at (yticklabel* cs:1);
\nextgroupplot[title={$N_M=2$, $S=1$}]
\addplot [matrix plot*, point meta=explicit, mesh/cols=3,
nodes near coords={
 \pgfmathfloattofixed{\pgfplotspointmeta}%
 \edef\MetaFixed{\pgfmathresult}%
 \pgfmathparse{\MetaFixed-50}%
 \ifdim\pgfmathresult pt>0pt \color{black}\else \color{white}\fi
 \pgfmathprintnumber[fixed,precision=1]{\pgfplotspointmeta}
},
nodes near coords align={center},
every node near coord/.append style={font=\footnotesize\bfseries}
] table[x=l, y=w, meta=val2] {\lkfpsThree};
\coordinate (lu5) at (yticklabel* cs:1);
\nextgroupplot[title={$N_M=0$, $S=2$}]
\addplot [matrix plot*, point meta=explicit, mesh/cols=3,
nodes near coords={
 \pgfmathfloattofixed{\pgfplotspointmeta}%
 \edef\MetaFixed{\pgfmathresult}%
 \pgfmathparse{\MetaFixed-50}%
 \ifdim\pgfmathresult pt>0pt \color{black}\else \color{white}\fi
 \pgfmathprintnumber[fixed,precision=1]{\pgfplotspointmeta}
},
nodes near coords align={center},
every node near coord/.append style={font=\footnotesize\bfseries}
] table[x=l, y=w, meta=val] {\lkfpsOne};
\coordinate (lu2) at (yticklabel* cs:1);
\nextgroupplot[title={$N_M=1$, $S=2$}]
\addplot [matrix plot*, point meta=explicit, mesh/cols=3,
nodes near coords={
 \pgfmathfloattofixed{\pgfplotspointmeta}%
 \edef\MetaFixed{\pgfmathresult}%
 \pgfmathparse{\MetaFixed-50}%
 \ifdim\pgfmathresult pt>0pt \color{black}\else \color{white}\fi
 \pgfmathprintnumber[fixed,precision=1]{\pgfplotspointmeta}
},
nodes near coords align={center},
every node near coord/.append style={font=\footnotesize\bfseries}
] table[x=l, y=w, meta=val] {\lkfpsTwo};
\coordinate (lu4) at (yticklabel* cs:1);
\nextgroupplot[title={$N_M=2$, $S=2$}]
\addplot [matrix plot*, point meta=explicit, mesh/cols=3,
 nodes near coords={
  \pgfmathfloattofixed{\pgfplotspointmeta}%
  \edef\MetaFixed{\pgfmathresult}%
  \pgfmathparse{\MetaFixed-50}%
  \ifdim\pgfmathresult pt>0pt \color{black}\else \color{white}\fi
  \pgfmathprintnumber[fixed,precision=1]{\pgfplotspointmeta}
 },
 nodes near coords align={center},
 every node near coord/.append style={font=\footnotesize\bfseries}
] table[x=l, y=w, meta=val] {\lkfpsThree};
\coordinate (lu6) at (yticklabel* cs:1);
\end{groupplot}
\draw[thick] (lu1) -- +(135:0.75) node[above,left,font=\bfseries] (A) {(a)};
\draw[thick] (lu2) -- +(135:0.75) node[above,left,font=\bfseries] (B) {(d)};
\draw[thick] (lu3) -- +(135:0.75) node[above,left,font=\bfseries] (C) {(b)};
\draw[thick] (lu4) -- +(135:0.75) node[above,left,font=\bfseries] (D) {(e)};
\draw[thick] (lu5) -- +(135:0.75) node[above,left,font=\bfseries] (E) {(c)};
\draw[thick] (lu6) -- +(135:0.75) node[above,left,font=\bfseries] (F) {(f)};
\begin{axis}[
 hide axis, scale only axis,
 height=4cm, width=0pt,
 at={($(myplots c3r1.east)!0.5!(myplots c3r2.east)$)}, anchor=west,
 colormap/viridis,
 colorbar,
 point meta min=20, point meta max=60,
 colorbar style={
  ytick distance=10,
  ylabel={LKFs},
  y label style={font=\footnotesize},
 },
]
 \addplot [draw=none] coordinates {(0,0)};
\end{axis}
\end{tikzpicture}
\end{center}
\end{table}

\pgfplotstableread[col sep=comma]{
conf,ps,epoch,max_abs,mean_abs,rmse,lkf_count
10,3,11,3.53903,0.124154433136095,0.285474626735255,43
02,2,38,1.45388,0.0268812118343195,0.0608195996963062,20
12,1,13,1.31848,0.0593552923076923,0.1215949145889409,41
11,3,7,1.93391,0.0842055147928994,0.164075866648408,42
21,2,9,2.60992,0.0922286390532544,0.251855977068608,41
12,3,16,2.54645,0.118464463905325,0.226632288072173,43
01,2,46,0.89124,0.0184953585798816,0.0419995778747454,39
20,3,3,3.24533,0.334932018934911,0.587822114617285,60
21,3,3,3.49369,0.295907633136095,0.525532845783115,44
02,1,25,1.31665,0.0443739195266272,0.1012410883677584,40
00,1,41,1.2603,0.0407588544378698,0.0877217885160604,41
11,1,46,1.20877,0.0440023786982248,0.1000015720933917,39
10,1,41,1.35857,0.0465493940828402,0.1002912223865782,41
11,2,47,1.37952,0.0237338106508875,0.0557714965029405,37
22,1,47,1.29944,0.0414871692307692,0.0990093514373297,39
20,2,7,2.60224,0.0629764118343195,0.167094468118165,38
21,1,47,1.07825,0.0402837775147929,0.0857327977248237,39
12,2,9,2.64628,0.0745961254437869,0.216763033464668,48
20,1,42,1.04635,0.0445708781065088,0.0919630868263299,40
01,1,47,1.24649,0.042589973964497,0.0915312164234732,41
10,2,47,0.984109999999999,0.0192159715976331,0.0439820326801133,37
22,3,3,3.632,0.283597197633136,0.495116668901047,56
00,3,39,1.43381,0.056137474556213,0.107088007688528,37
00,2,35,1.59125,0.0362286934911242,0.109943758193039,37
22,2,13,2.28873,0.0583030319526627,0.164997260263645,40
01,3,3,2.87449,0.0758820863905325,0.151764172781065,49
}\datatable

\def\ConfList{
 00/$2\times 128$,
 01/$2\times 256$,
 02/$2\times 512$,
 10/$4\times 128$,
 11/$4\times 256$,
 12/$4\times 512$,
 20/$8\times 128$,
 21/$8\times 256$,
 22/$8\times 512$
}

\pgfplotsset{
 myaxis/.style={
  width=.75\linewidth, height=5cm,
  grid=both, grid style={line width=.2pt},
  axis line style={semithick}, tick style={semithick},
  scaled ticks=false,
  legend columns=2, legend cell align=left,
  legend style={/tikz/every even column/.append style={column
    sep=0.2cm},at={(1.05,0.5)},draw=none,anchor=west},
  clip=false
 }
}
\pgfplotsset{
 mygroupaxis/.style={
  width=.5\linewidth, height=5cm,
  grid=both, grid style={line width=.2pt},
  axis line style={semithick}, tick style={semithick},
  scaled ticks=false,
  legend columns=9, legend cell align=left,
  legend style={/tikz/every even column/.append style={column
    sep=0.1cm},at={(-0.2,1.05)},draw=none,anchor=south west},
  clip=false
 }
}
\pgfplotscreateplotcyclelist{foo}{
{mark=*,    mark size=3pt},
{mark=square*, mark size=3pt},
{mark=triangle*,mark size=3pt},
{mark=diamond*, mark size=3pt},
{mark=pentagon*,   mark size=3pt},
{mark=halfcircle*,    mark size=3pt},
{mark=halfsquare left*,  mark size=3pt},
{mark=halfsquare right*, mark size=3pt},
{mark=10-pointed star, mark size=3pt},
}
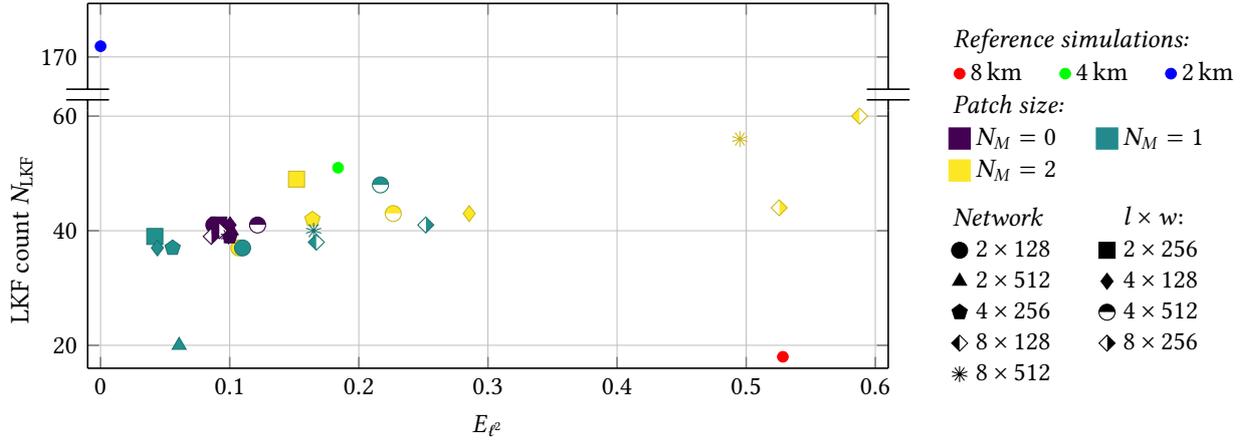
\begin{figure}\centering
 \pgfplotsset{
  every non boxed x axis/.style={}
 }
\begin{tikzpicture}[font=\small]
 \begin{groupplot}[
  group style={
   group size=1 by 2,
   xticklabels at=edge bottom,
   vertical sep=0pt
  },
 myaxis,xmin=-0.01,xmax=0.61,
 cycle list name=foo,
 colormap/viridis,
 point meta min=1, point meta max=3
 ]
 \nextgroupplot[myaxis,ymin=165,ymax=175,ytick={170},height=3cm,axis x line=top,
 axis y discontinuity=parallel,legend columns=3, legend cell align=left,
 legend style={at={(1.05,0.5)},draw=none,anchor=west}]
 \addlegendimage{empty legend}%
 \addlegendentry{}%
 \addlegendimage{empty legend}%
 \addlegendentry{}%
 \addlegendimage{empty legend}%
 \addlegendentry{\hspace{-3cm}\textit{Reference simulations:}}
 \addlegendimage{only marks, mark=*, draw=red,fill=red}
 \addlegendentry{\SI{8}{\kilo\meter}}%
 \addlegendimage{only marks, mark=*, draw=green,fill=green}
 \addlegendentry{\SI{4}{\kilo\meter}}%
 \addlegendimage{only marks, mark=*, draw=blue,fill=blue}
 \addlegendentry{\SI{2}{\kilo\meter}}%
 \addplot[
 only marks, blue
 ] coordinates {(0.0,171)};
 \nextgroupplot[myaxis,ymin=16,ymax=61,xlabel={$E_{\ell^2}$}, ylabel={LKF count
  $N_{\mathrm{LKF}}$},axis x line=bottom]
 \addlegendimage{empty legend}
 \addlegendentry{\hspace{-.3cm}\textit{Patch size:}}
 \addlegendimage{empty legend}%
 \addlegendentry{}%
 \addlegendimage{only marks, mark size=4pt, mark=square*, draw=mapped color, fill=mapped color, index of colormap=0}
 \addlegendentry{$N_M=0$}
 \addlegendimage{only marks, mark size=4pt, mark=square*, draw=mapped color, fill=mapped color, index of colormap=8}
 \addlegendentry{$N_M=1$}
 \addlegendimage{only marks, mark size=4pt, mark=square*, draw=mapped color, fill=mapped color, index of colormap=17}
 \addlegendentry{$N_M=2$}
 \addlegendimage{empty legend}%
 \addlegendentry{}%
 \addlegendimage{empty legend}%
 \addlegendentry{}%
 \addlegendimage{empty legend}%
 \addlegendentry{}%
 \addlegendimage{empty legend}%
 \addlegendentry{\hspace{-.3cm}\textit{Network}}%
 \addlegendimage{empty legend}%
 \addlegendentry{$l\times w$:}%
\foreach \y/\x in \ConfList {
 \addplot+[legend image post style={
  /pgfplots/scatter/use mapped color={draw=black,fill=black},
  draw=black, fill=black, mark options={draw=black,fill=black}
 },only marks,
 scatter, scatter src=explicit      %
 ]
 table[x=rmse, y=lkf_count, meta=ps,
 restrict expr to domain={\thisrow{conf}}{\y:\y}
 ]{\datatable};
 \addlegendentryexpanded{\x}%
}
\addplot[
only marks, red,
] coordinates {(0.52838852775,18)};
\addplot[
only marks, green
] coordinates {(0.1838852775,51)};
\end{groupplot}
\end{tikzpicture}
\caption{\label{fig:lkf-scatter} \emph{Scatter plot with the error $ E_{\ell^2}$ between the reference solution, $\vt^\text{ref}_\text{fine}$, and the corrected network solution, $\vt_\text{fine}$, over the
 number of LKFs (cf.~\ref{itm:metric-L2},~\ref{itm:metric-LKF}) under varying patch sizes ($N_M \in \{0,1,2\}$; in purple, green and yellow, respectively) and for different depth and width, $l \times w$, of the network indicated by markers. The jump level is set to $S=2$. The red, blue and green dot indicates the number of LKFs detected from the reference solutions (without NN) on the meshes with side length \SI{8}{\kilo\meter}, \SI{4}{\kilo\meter} and \SI{2}{\kilo\meter}, respectively.  All combinations are able to offer a reasonable trade-off
 between accuracy with respect to the reference solution and the number of LKFs.}}
\end{figure}

\pgfplotstableread[col sep=comma]{
conf,group_key,ps,n_files,copies_for_network_mean,copies_for_network_std,network_mean,network_std,solver_mean,solver_std,total_mean,total_std,newton_sum_mean,newton_sum_std,mean_abs,rmse,lkf
00,$2\times 128$,1,8,0.37124999999999997,0.01246423454758226,0.19,0.018516401995451025,238.37375,3.4060028416228456,249.33875,3.4127049912265814,546.0,4.956957592256421,0.0407588544378698,0.0877217885160604,41
01,$2\times 256$,1,8,0.37625,0.013024701806293204,0.2025,0.0158113883008419,241.76375000000002,4.0591270613273505,252.83375,4.230683624258779,541.5,8.053393251117374,0.042589973964497,0.0915312164234732,41
02,$2\times 512$,1,8,0.37625,0.010606601717798222,0.2,0.015118578920369089,240.63125,4.454523984509618,251.6425,4.680460142214355,547.0,7.050835816716038,0.0443739195266272,0.1012410883677584,40
10,$4\times 128$,1,8,0.37124999999999997,0.015526475085202983,0.20625000000000002,0.015979898086569355,241.6125,2.198855221635633,252.59625,2.342269580068506,548.5,4.44007721940573,0.0465493940828402,0.1002912223865782,41
11,$4\times 256$,1,8,0.37625,0.0140788595317336,0.2,0.011952286093343936,240.32874999999999,3.6016640499318417,251.24625,3.716864798263639,540.75,7.977647343851265,0.0440023786982248,0.1000015720933917,39
12,$4\times 512$,1,8,0.38,0.016903085094570332,0.22625,0.019226098333849674,240.47125,3.2421837327861045,251.43875,3.3405407496391986,542.0,4.7207747548166585,0.0593552923076923,0.1215949145889409,41
20,$8\times 128$,1,8,0.37625,0.013024701806293204,0.24875,0.018077215335491097,241.015,3.4005125664060887,252.15625,3.5550645463298514,546.75,5.4182232195751086,0.0445708781065088,0.0919630868263299,40
21,$8\times 256$,1,8,0.37625,0.011877349391654218,0.25375,0.018468119248354137,239.10625,4.858965020601455,250.12875,4.938841499495902,550.0,8.176621726430962,0.0402837775147929,0.0857327977248237,39
22,$8\times 512$,1,7,0.37285714285714283,0.012535663410560186,0.19428571428571428,0.01511857892036909,241.37,4.974635665051261,252.28142857142856,4.9631959749454015,545.8571428571429,3.8913824205360674,0.0414871692307692,0.0990093514373297,39
00,$2\times 128$,2,8,0.20625,0.00916125381312904,0.21375,0.024458419526091332,243.13875000000002,2.158302491046409,254.215,2.2334150660239733,549.5,6.301927142889365,0.0362286934911242,0.109943758193039,37
01,$2\times 256$,2,8,0.2025,0.008864052604279178,0.2125,0.03732100136460894,241.665,4.334855410342009,252.4475,4.517460569833459,551.25,5.4707011825333165,0.0184953585798816,0.0419995778747454,39
02,$2\times 512$,2,8,0.20375,0.007440238091428443,0.19625,0.015059406173077158,243.41375,10.323186021905116,254.11624999999998,10.53763857729587,554.5,18.228705776486553,0.0268812118343195,0.0608195996963062,20
10,$4\times 128$,2,8,0.20375,0.005175491695067647,0.215,0.014142135623730947,243.32125,3.424618591901878,254.2075,3.6142683511722984,546.25,3.7321001364608946,0.0192159715976331,0.0439820326801133,37
11,$4\times 256$,2,8,0.20375,0.005175491695067647,0.2175,0.02052872551885702,243.22875,5.295694848513404,254.025,5.6717973971471976,549.0,8.106434833777776,0.0237338106508875,0.0557714965029405,37
12,$4\times 512$,2,8,0.20375,0.007440238091428443,0.24125,0.009910312089651149,244.64499999999998,6.478727829266658,255.73875,6.883403305882106,544.25,3.882193783343198,0.0745961254437869,0.216763033464668,48
20,$8\times 128$,2,8,0.20625,0.007440238091428443,0.26375,0.03159452936370924,244.30125,4.671859334660788,255.26875,4.97594624884842,556.75,2.8660575211055543,0.0629764118343195,0.167094468118165,38
21,$8\times 256$,2,8,0.20625,0.007440238091428443,0.26,0.02138089935299395,240.23250000000002,5.280383508799337,251.21125,5.443555659151998,547.75,11.067971810589327,0.0922286390532544,0.251855977068608,41
22,$8\times 512$,2,8,0.2025,0.00886405260427918,0.24875,0.011259916264596043,241.69,4.891765675967262,252.50875000000002,5.178387607574946,545.25,5.257647491307985,0.0583030319526627,0.164997260263645,40
00,$2\times 128$,3,8,0.1675,0.010350983390135314,0.20125,0.02587745847533828,251.69,5.713867841863838,262.4725,5.878152649551443,568.75,9.452890714334048,0.056137474556213,0.107088007688528,37
01,$2\times 256$,3,8,0.1625,0.008864052604279192,0.2,0.013093073414159545,244.6775,3.3019766807171735,255.46875,3.4967431275402605,553.25,2.052872551885702,0.0758820863905325,0.151764172781065,49
02,$2\times 512$,3,8,0.16625,0.007440238091428449,0.21875,0.018077215335491087,248.14125,7.201000798500164,259.13,7.259866193178573,558.25,14.587176364386437,0.11193504852071,0.231829932454934,44
10,$4\times 128$,3,8,0.1625,0.007071067811865481,0.21375,0.015059406173077149,243.1,6.312364284255374,254.05625,6.4507583994885085,551.75,14.626297256263166,0.124154433136095,0.285474626735255,43
11,$4\times 256$,3,8,0.1625,0.004629100498862762,0.22375,0.02825268634509444,248.115,6.07578566348184,259.0675,5.948624690992501,563.0,17.3699247469379,0.0842055147928994,0.164075866648408,42
12,$4\times 512$,3,8,0.1625,0.004629100498862762,0.21625,0.011877349391654208,243.03625,6.702692289563309,253.96375,6.785002237077216,556.0,17.5662013130736,0.118464463905325,0.226632288072173,43
20,$8\times 128$,3,8,0.165,0.005345224838248492,0.25375000000000003,0.015059406173077166,247.58499999999998,8.363451099020919,258.5425,8.61042184133357,562.0,17.5662013130736,0.334932018934911,0.587822114617285,60
21,$8\times 256$,3,8,0.16625,0.005175491695067661,0.2475,0.013887301496588273,243.36125,4.365541817133679,254.33625,4.422988445125569,559.5,13.103979766249859,0.295907633136095,0.525532845783115,44
22,$8\times 512$,3,8,0.1625,0.007071067811865481,0.24125,0.01246423454758225,241.26125000000002,5.974468147279482,252.01999999999998,6.085174020284287,546.0,19.493588689617926,0.283597197633136,0.495116668901047,56
33,\SI{8}{\kilo\meter},5,9,0.0,0.0,0.0,0.0,268.,9,277.807,14.324357714488616,680,4,0.0,0.0,18
}\sumtab
\pgfplotscreateplotcyclelist{ps}{
 {mark=*, mark size=3pt, draw=mapped color, fill=mapped color, index of colormap=0},
 {mark=*, mark size=3pt, draw=mapped color, fill=mapped color, index of colormap=8},
 {mark=*, mark size=3pt, draw=mapped color, fill=mapped color, index of colormap=17},
 {mark=*, mark size=3pt, draw=red, fill=red},
}
\begin{figure}\centering
\begin{tikzpicture}[font=\small]
 \begin{groupplot}[
  group style={group name= myplot,group size=2 by 1, horizontal sep=1.2cm, vertical sep=1.8cm},
  width=.5\linewidth, height=4.5cm,
 cycle list name=foo,
 colormap/viridis,
 point meta min=1, point meta max=3
 ]
 \nextgroupplot[mygroupaxis,ylabel={$E_{\ell^2}$},xlabel={$T_{\mathrm{mom}}$},ymin=-0.01,ymax=0.61]
 \addlegendimage{empty legend}
 \addlegendentry{\hspace{-0.3cm}\textit{Patch size:}}
 \addlegendimage{only marks, mark size=4pt, mark=square*, draw=mapped color, fill=mapped color, index of colormap=0}
 \addlegendentry{$N_M=0$}
 \addlegendimage{only marks, mark size=4pt, mark=square*, draw=mapped color, fill=mapped color, index of colormap=8}
 \addlegendentry{$N_M=1$}
 \addlegendimage{only marks, mark size=4pt, mark=square*, draw=mapped color, fill=mapped color, index of colormap=17}
 \addlegendentry{$N_M=2$}
 \addlegendimage{empty legend}%
 \addlegendentry{}%
 \addlegendimage{empty legend}%
 \addlegendentry{}%
 \addlegendimage{empty legend}%
 \addlegendentry{}%
 \addlegendimage{empty legend}%
 \addlegendentry{}%
 \addlegendimage{only marks, mark=*, draw=red,fill=red}
 \addlegendentry{\SI{8}{\kilo\meter}}
 \addlegendimage{empty legend}%
 \addlegendentry{}%
 \addlegendimage{empty legend}%
 \addlegendentry{}%
 \addlegendimage{empty legend}%
 \addlegendentry{}%
 \addlegendimage{empty legend}%
 \addlegendentry{}%
 \addlegendimage{empty legend}%
 \addlegendentry{}%
 \addlegendimage{empty legend}%
 \addlegendentry{}%
 \addlegendimage{empty legend}%
 \addlegendentry{}%
 \addlegendimage{empty legend}%
 \addlegendentry{}%
 \addlegendimage{empty legend}%
 \addlegendentry{\hspace{-24.6cm}\textit{Neural Network $l\times w$}:}%
\foreach \y/\x in \ConfList {
 \addplot+[legend image post style={
  /pgfplots/scatter/use mapped color={draw=black,fill=black},
  draw=black, fill=black, mark options={draw=black,fill=black}
 },only marks,
 scatter, scatter src=explicit      %
 ]
 table[x=total_mean, y=rmse, meta=ps,
 restrict expr to domain={\thisrow{conf}}{\y:\y}
 ]{\sumtab};
 \addlegendentryexpanded{\x}%
 }
 \addplot+[
 only marks, red
 ] coordinates {(277.807,0.52838852775)};
 \nextgroupplot[mygroupaxis,ylabel={LKF},xlabel={$T_{\mathrm{mom}}$},ymin=16,ymax=61]
\foreach \y/\x in \ConfList {
 \addplot+[legend image post style={
  /pgfplots/scatter/use mapped color={draw=black,fill=black},
  draw=black, fill=black, mark options={draw=black,fill=black}
 },only marks,
 scatter, scatter src=explicit      %
 ]
 table[x=total_mean, y=lkf, meta=ps,
 restrict expr to domain={\thisrow{conf}}{\y:\y}
 ]{\sumtab};
}
\addplot+[
only marks, red
] coordinates {(277.807,18)};
\end{groupplot}
\node[font=\bfseries, anchor=north west] at (myplot c1r1.north west) {(a)};
\node[font=\bfseries, anchor=north west] at (myplot c2r1.north west) {(b)};
\end{tikzpicture}
\caption{\label{fig:walltime-scatter} \emph{Scatter plot with time to solution to solve the sea-ice momentum equation for varying patch sizes ($N_M \in\{0,1,2\}$; in purple, green and yellow, respectively) and for different network configurations analogous to Fig.~\ref{fig:lkf-scatter}. The time is presented with respect to (a) the error $E_{\ell^2}$ between the reference solution with $\SI{2}{\kilo\meter}$ resolution and the corrected network solution, and (b)
  LKF count, $N_{\mathrm{LKF}}$. Solving the sea-ice momentum without NN requires $T_{\mathrm{mom}}:$ \SI{277}{\second},
 \SI{2751}{\second},
 \SI{29656}{\second},
 at \SI{8}{\kilo\meter}, \SI{4}{\kilo\meter} and \SI{2}{\kilo\meter} resolution, respectively.
 }}
\end{figure}
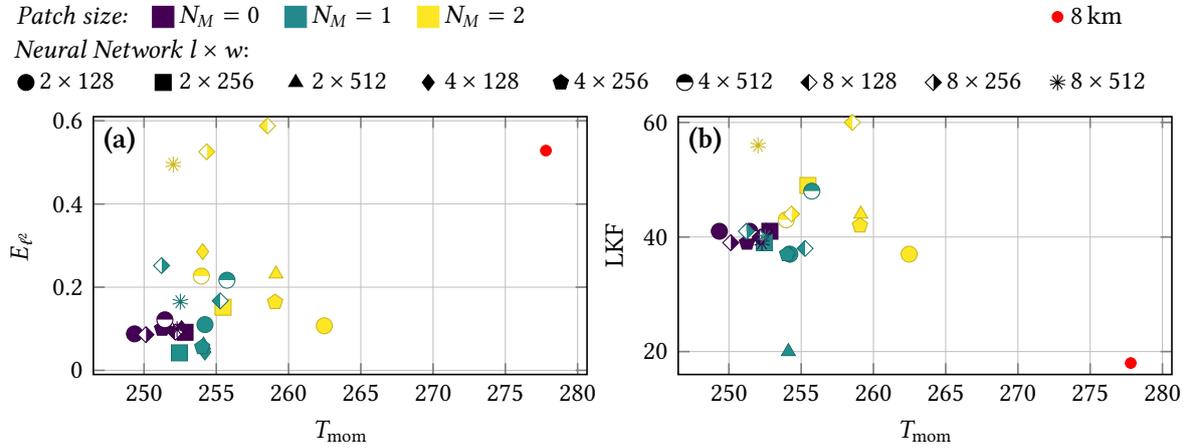

\pgfplotstableread[col sep=comma]{
conf,group_key,ps,n_files,copies_for_network_mean,copies_for_network_std,network_mean,network_std,solver_mean,solver_std,total_mean,total_std,newton_sum_mean,newton_sum_std,mean_abs,rmse,lkf
01,$2\times 256$,1,8,0.37625,0.013024701806293204,0.2025,0.0158113883008419,241.76375000000002,4.0591270613273505,252.83375,4.230683624258779,541.5,8.053393251117374,0.042589973964497,0.0915312164234732,41
11,$4\times 256$,1,8,0.37625,0.0140788595317336,0.2,0.011952286093343936,240.32874999999999,3.6016640499318417,251.24625,3.716864798263639,540.75,7.977647343851265,0.0440023786982248,0.1000015720933917,39
21,$8\times 256$,1,8,0.37625,0.011877349391654218,0.25375,0.018468119248354137,239.10625,4.858965020601455,250.12875,4.938841499495902,550.0,8.176621726430962,0.0402837775147929,0.0857327977248237,39
01,$2\times 256$,2,8,0.2025,0.008864052604279178,0.2125,0.03732100136460894,241.665,4.334855410342009,252.4475,4.517460569833459,551.25,5.4707011825333165,0.0184953585798816,0.0419995778747454,39
11,$4\times 256$,2,8,0.20375,0.005175491695067647,0.2175,0.02052872551885702,243.22875,5.295694848513404,254.025,5.6717973971471976,549.0,8.106434833777776,0.0237338106508875,0.0557714965029405,37
21,$8\times 256$,2,8,0.20625,0.007440238091428443,0.26,0.02138089935299395,240.23250000000002,5.280383508799337,251.21125,5.443555659151998,547.75,11.067971810589327,0.0922286390532544,0.251855977068608,41
01,$2\times 256$,3,8,0.1625,0.008864052604279192,0.2,0.013093073414159545,244.6775,3.3019766807171735,255.46875,3.4967431275402605,553.25,2.052872551885702,0.0758820863905325,0.151764172781065,49
11,$4\times 256$,3,8,0.1625,0.004629100498862762,0.22375,0.02825268634509444,248.115,6.07578566348184,259.0675,5.948624690992501,563.0,17.3699247469379,0.0842055147928994,0.164075866648408,42
21,$8\times 256$,3,8,0.16625,0.005175491695067661,0.2475,0.013887301496588273,243.36125,4.365541817133679,254.33625,4.422988445125569,559.5,13.103979766249859,0.295907633136095,0.525532845783115,44
33,\SI{8}{\kilo\meter},5,9,0.0,0.0,0.0,0.0,268.,9,277.807,14.324357714488616,680,4,0.0,0.0,18
}\sumtabII
\pgfplotscreateplotcyclelist{psII}{
 {mark=*, mark size=3pt, draw=mapped color, fill=mapped color, index of colormap=0},
 {mark=*, mark size=3pt, draw=mapped color, fill=mapped color, index of colormap=8},
 {mark=*, mark size=3pt, draw=mapped color, fill=mapped color, index of colormap=17},
 {mark=*, mark size=3pt, draw=red, fill=red},
 {mark=*, mark size=3pt, draw=green, fill=green},
}
\begin{figure}
\begin{tikzpicture}[font=\small]
\begin{groupplot}[
 group style={group name=myplot,group size=2 by 2, horizontal sep=2cm, vertical sep=1.5cm},
 width=.48\linewidth, height=4.5cm,
 ymajorgrids, ymin=0,
 xtick=data, title style={at={(0.5,0.95)}},
 x tick label style={rotate=45, yshift=4pt, xshift=2pt, font=\footnotesize,
  anchor=north east},
 xticklabels from table={\sumtabII}{group_key},
 enlarge x limits=0.06,
 colormap/viridis,
 cycle list name=ps,
 legend columns=7, legend cell align=left,
 legend style={/tikz/every even column/.append style={column
   sep=0.2cm},at={(0.0,1.05)},draw=none,anchor=south west},
]
\nextgroupplot[ylabel={\(T_{\mathrm{mom}}\)}, title style={at={(0.5,1.2)}},ymax=450]
\addlegendimage{empty legend}
\addlegendentry{\textit{Numbers above bars: }Total Newton iters}
\addlegendimage{empty legend}
\addlegendentry{$\quad$}
\addlegendimage{empty legend}
\addlegendentry{\textit{Patch size}}
\addlegendimage{only marks, mark size=4pt, mark=square*, draw=mapped color, fill=mapped color, index of colormap=0}
\addlegendentry{$N_M=0$}
\addlegendimage{only marks, mark size=4pt, mark=square*, draw=mapped color, fill=mapped color, index of colormap=8}
\addlegendentry{$N_M=1$}
\addlegendimage{only marks, mark size=4pt, mark=square*, draw=mapped color, fill=mapped color, index of colormap=17}
\addlegendentry{$N_M=2$}
\addplot+[forget plot,
ybar, bar width=0pt,fill,
error bars/.cd, y dir=both, y explicit
] table[
x expr=\coordindex,
y=total_mean,
]{\sumtabII};
\foreach \P in {1,2,3,5,6}{%
\addplot+[forget plot,
ybar, bar width=0.8pt,fill,
error bars/.cd, y dir=both, y explicit
] table[
x expr=\coordindex,
y=total_mean,
y error=total_std,
restrict expr to domain={\thisrow{ps}}{\P:\P}
]{\sumtabII};
\addplot+[forget plot,only marks, nodes near coords,
visualization depends on=\thisrow{newton_sum_mean}\as\newtonlabel,
nodes near coords style={font=\footnotesize,text=gray,yshift=12pt, xshift=6pt,rotate=90},
nodes near coords={\pgfmathprintnumber[precision=0]{\newtonlabel}}
] table[
x expr=\coordindex,
y=total_mean,
restrict expr to domain={\thisrow{ps}}{\P:\P}
]{\sumtabII};
\addplot+[only marks, nodes near coords,
nodes near coords style={font=\footnotesize,text=black,fill=white, yshift=-1cm, xshift=6pt,rotate=90},
] table[
x expr=\coordindex,
y=total_mean,
restrict expr to domain={\thisrow{ps}}{\P:\P}
]{\sumtabII};
}
\nextgroupplot[ylabel={\(T_{\mathrm{NK}}\)}]
\addplot+[forget plot,
ybar, bar width=0pt,fill,
error bars/.cd, y dir=both, y explicit
] table[
x expr=\coordindex,
y=solver_mean,
]{\sumtabII};
\foreach \P in {1,2,3,5,6}{%
\addplot+[forget plot,
ybar, bar width=0.8pt,fill,
error bars/.cd, y dir=both, y explicit
] table[
x expr=\coordindex,
y=solver_mean,
y error=solver_std,
restrict expr to domain={\thisrow{ps}}{\P:\P}
]{\sumtabII};
\addplot+[only marks, nodes near coords,
nodes near coords style={font=\footnotesize,text=black,yshift=-1cm,fill=white,xshift=6pt,rotate=90},
] table[
x expr=\coordindex,
y=solver_mean,
restrict expr to domain={\thisrow{ps}}{\P:\P}
]{\sumtabII};
}
\nextgroupplot[ylabel={\(T_{\mathrm{cpy}}\)}]
\addplot+[forget plot,
ybar, bar width=0pt,fill,
error bars/.cd, y dir=both, y explicit
] table[
x expr=\coordindex,
y=copies_for_network_mean,
restrict expr to domain={\thisrow{ps}}{1:3}
]{\sumtabII};
\foreach \P in {1,2,3}{%
 \addplot+[forget plot,
 ybar, bar width=0.8pt,fill,
 error bars/.cd, y dir=both, y explicit
] table[
 x expr=\coordindex,
 y=copies_for_network_mean,
 y error=copies_for_network_std,
 restrict expr to domain={\thisrow{ps}}{\P:\P}
]{\sumtabII};
\addplot+[
 only marks, nodes near coords,
 point meta=y,
 nodes near coords style={font=\footnotesize,text=black, yshift=-0.66cm,fill=white, xshift=6pt,rotate=90},
] table[
 x expr=\coordindex,
 y=copies_for_network_mean,
 restrict expr to domain={\thisrow{ps}}{\P:\P}
]{\sumtabII};
}

\nextgroupplot[ylabel={\(T_{\mathrm{nn}}\)}]
\addplot+[forget plot,
ybar, bar width=0pt,fill,
error bars/.cd, y dir=both, y explicit
] table[
x expr=\coordindex,
y=network_mean,
restrict expr to domain={\thisrow{ps}}{1:3}
]{\sumtabII};
\foreach \P in {1,2,3}{%
\addplot+[forget plot,
 ybar, bar width=0.8pt,fill,
 error bars/.cd, y dir=both, y explicit
] table[
 x expr=\coordindex,
 y=network_mean,
 y error=network_std,
 restrict expr to domain={\thisrow{ps}}{\P:\P}
]{\sumtabII};
\addplot+[only marks, nodes near coords,
nodes near coords style={font=\footnotesize,text=black, yshift=-1cm,fill=white, xshift=6pt,rotate=90},
] table[
 x expr=\coordindex,
 y=network_mean,
 restrict expr to domain={\thisrow{ps}}{\P:\P}
 ]{\sumtabII};
}
\end{groupplot}
\node[font=\bfseries, anchor=north west] at (myplot c1r1.north west) {(a)};
\node[font=\bfseries, anchor=north west] at (myplot c2r1.north west) {(b)};
\node[font=\bfseries, anchor=north west,fill=white,fill opacity=0.5,text opacity=1] at (myplot c1r2.north west) {(c)};
\node[font=\bfseries, anchor=north west] at (myplot c2r2.north west) {(d)};
\end{tikzpicture}
\caption{\label{fig:timings}\emph{Timing of NN-FEM approach across patch sizes with respect to the sea-ice momentum solve, Newton-Krylov
 solver time, host-device copy time, and network inference time. The run times are presented for varying patch sizes
 \(N_M\in\{0,1,2\}\) (colored purple, green and yellow, respectively), different network widths and depths (presented on the horizontal axis). Numbers above
 bars in panel (a) (colored in gray) denote total Newton iterations \(N_{\text{Newton}}\). The red
\SI{8}{\kilo\meter} label in panel (a) and (b) denotes the solution at \SI{8}{\kilo\meter} grid
 resolution without NN-FEM corrections. On a mesh with side length of \SI{4}{\kilo\meter} and \SI{2}{\kilo\meter}, the run time without NN-FEM is
 \SI{2751}{\second},
 \SI{29656}{\second},
  respectively.
 }}
\end{figure}
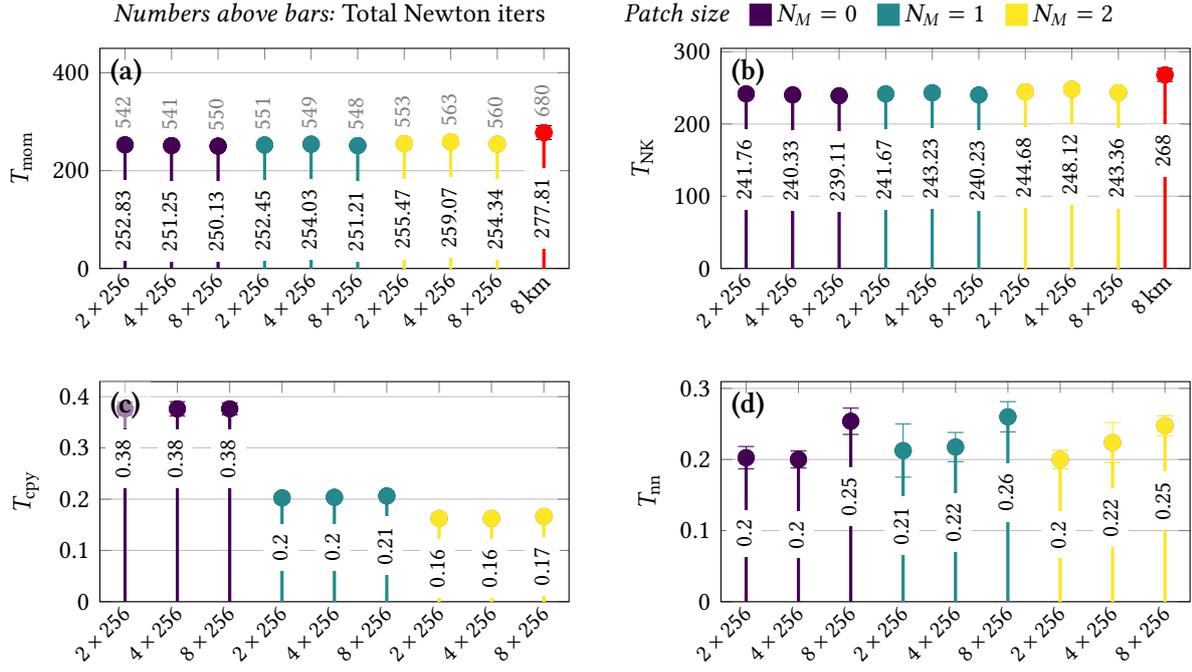
We study how the \emph{patch size} (collection of neighboring mesh cells), \(N_M\), the \emph{prediction jump} \(S\) (which
controls the number of degrees of freedom in the auxiliary space \(\V_\text{fine}\)),
and the \emph{network depth and width} \((l,w)\) affect the performance of the hybrid
NN-FEM method with respect to simulate LKFs. We vary \(N_M\in\{0,1,2\}\), \(S\in\{1,2\}\), and
\(l\in\{2,4,8\}\), \(w\in\{128,256,512\}\), see Tab.~\ref{tab:nn-parameters}.
Accuracy and efficiency are quantified by the metrics~\ref{itm:metric-L2}--\ref{itm:metric-time}. In the following we evaluate the simulated shear deformation qualitatively by visual comparison (Section \ref{sec:M1}). Accuracy of the NN-FEM is addressed by deriving the error \(E_{\ell^2}\) between the NN-corrected
solution and the high-resolution reference~\ref{itm:metric-L2} combined with an analysis of the number of
LKFs \(N_{\mathrm{LKF}}\)~\ref{itm:metric-LKF} (Section \ref{sec:M2}). The Newton-iterations \ref{itm:metric-Newton}  and the run time \ref{itm:metric-time} required for simulations with the NN-FEM are evaluated in Section \ref{sec:M4}. The main findings are summarized in Section \ref{sec:Msum}.

\subsubsection{Comparison of shear deformation fields.}\label{sec:M1}
\emph{Single-level prediction (\(S=1\)).}
We first assess accuracy for single-level prediction, where the NN is trained on a
\SI{4}{\kilo\meter} mesh (\(S=1\)) and used to correct simulations on the
\SI{8}{\kilo\meter} working mesh. We exemplarily visualize the NN-FEM simulations for an anti-cyclone setting with northwest wind direction. All other setups perform qualitatively similar. In
Fig.~\ref{fig:grid-j1-p1}--\ref{fig:grid-j1-p3} we show the NN-corrected fields on
the \SI{8}{\kilo\meter} resolution grid for
\(N_M \in \{0,1,2\}\). These are compared with the \SI{8}{\kilo\meter} baseline
without NN (Fig.~\ref{fig:grid-8k-base}) and the \SI{4}{\kilo\meter} target
(Fig.~\ref{fig:grid-4k-ref}). The visual comparison shows that increasing \(N_M\)
systematically improves the approximation relative to the baseline. We summarize
the observations for each patch size:
\begin{description}[leftmargin=*,labelsep=0.5em,itemsep=0pt,parsep=0pt,topsep=0pt]
\item[$N_M=0$:] The smallest patch setting
 (Fig.~\ref{fig:grid-j1-p1}) stays close to the \SI{8}{\kilo\meter} baseline
 (Fig.~\ref{fig:grid-8k-base}); the influence of the network correction on the simulated shear deformation is barley visible.
\item[$N_M=1$:] Compared to \(N_M=0\), the larger patch size visually improves the agreement
 with the \SI{4}{\kilo\meter} reference (Fig.~\ref{fig:grid-j1-p2}).
\item[$N_M=2$:] The setting
 (Fig.~\ref{fig:grid-j1-p3}) approaches the \SI{4}{\kilo\meter} target
 (Fig.~\ref{fig:grid-4k-ref}) even better than \(N_M=1\). Larger patches provide
 more sharply defined features and enable richer corrections compared to the baseline.
\end{description}

\emph{Two-level prediction (\(S=2\)).}
We now consider two-level prediction with \(S=2\), where the NN is trained on a
\SI{2}{\kilo\meter} mesh and used to correct the \SI{8}{\kilo\meter} working mesh.
Fig.~\ref{fig:grid-j2-p1}--\ref{fig:grid-j2-p3} shows snapshots of the NN-FEM
solution for \(N_M\in\{0,1,2\}\). We compare these to the \SI{8}{\kilo\meter}
baseline without NN (Fig.~\ref{fig:grid-8k-base}) and the \SI{2}{\kilo\meter}
target (Fig.~\ref{fig:grid-2k-ref}). For larger \(N_M\), we obtain approximations
that are qualitatively similar to a \SI{4}{\kilo\meter} resolution
(cf.~Fig.~\ref{fig:grid-4k-ref}), even though the corrections are formed on the
\SI{2}{\kilo\meter} auxiliary grid. For fixed \(N_M\), raising the jump level from
\(S=1\) to \(S=2\) yields  higher LKF counts
across the tested setups (Fig.~\ref{fig:grid-j2-p3},
Tab.~\ref{tab:lkf-number}(d)--\ref{tab:lkf-number}(f). We did not observe
patch-boundary artifacts in the analyzed cases. The auxiliary space \(\V_\text{fine}\)
constructed on a \SI{2}{\kilo\meter} grid allows the patch corrections to represent
higher spatial frequencies than those obtained from an auxiliary space defined on
a \SI{4}{\kilo\meter} grid. Comparing the second and third rows of
Fig.~\ref{fig:ss-3x3} and Tab.~\ref{tab:lkf-number}(a)--\ref{tab:lkf-number}(c) with
Tab.~\ref{tab:lkf-number}(d)--\ref{tab:lkf-number}(f) shows that \(S=2\) yields a
qualitative improvement over \(S=1\). Therefore, we focus on \(S=2\) in what follows.

\begin{remark}[NN-FEM approach on large patches] We exclude the results for \(N_M=3\) combined with \(S=1,2\), as the network corrections
produce numerous spurious high-frequency artifacts. A plausible explanation is
data scarcity relative to the input and parameter dimension: as \(N_M\) increases,
the feature dimension grows while the number of statistically independent patches
per snapshot decreases, worsening the sample-to-parameter ratio. With
substantially larger training datasets, \(N_M=3\) could become advantageous;
however, under the current data volume it leads to noisy predictions.
\end{remark}

\subsubsection{Metric~\ref{itm:metric-L2} and \ref{itm:metric-LKF}: Accurancy and LKF Counts}\label{sec:M2}
Tab.~\ref{tab:lkf-number} summarizes the LKF counts obtained from simulations on
the \SI{8}{\kilo\meter} mesh corrected with the NN-FEM approach for different patch
sizes \(N_M\in\{0,1,2\}\) and jump levels \(S\in\{1,2\}\). Increasing \(N_M\)
increases the number of resolved LKFs across architectures, while changes in
network depth and width \((l,w)\) shift the absolute counts only slightly. In
Fig.~\ref{fig:lkf-scatter} we plot the error \(E_{\ell^2}\) versus the number of
LKFs \(N_{\mathrm{LKF}}\) for \(S=2\). All combinations achieve a reasonable
trade-off between accuracy and additional LKF features under varying patch sizes
and architectures. The patch size is the dominant factor, whereas the architecture
modulates the absolute counts without changing the ranking
(Tab.~\ref{tab:lkf-number}a--\ref{tab:lkf-number}f).

For \(S=2\) we can resolve the effect of \(N_M\) and \((l,w)\) in more detail:
\begin{description}[leftmargin=*,itemsep=0pt,topsep=0pt,parsep=0pt]
\item[$N_M=0$:]
 In terms of LKFs the NN-FEM approach yields a moderate increase over the
 baseline count of 18 LKFs (Tab.~\ref{tab:lkf-number}). The model size cannot
 compensate for the limited local context, but this configuration provides a
 low-cost baseline when memory is constrained.
\item[$N_M=1$:] The \(l=2\), \(w=256\) model is competitive, but higher-capacity
 networks are advantageous. A model with \(l=4\) and \(w=512\) performs best in
 terms of LKFs (Tab.~\ref{tab:lkf-number}), and a model with \(l=8\) and
 \(w=256\) exhibits similar LKF counts. The larger number of training patches
 for \(N_M=1\) supports stable training of the larger networks and compensates
 for the reduced context. This choice offers a good runtime/accuracy trade-off
 and was consistent across the tested scenarios.
\item[$N_M=2$:]
 A compact model
 with \(l=2\), \(w=256\) strikes a favorable error-LKF balance
 (Fig.~\ref{fig:lkf-scatter}). Deeper or wider variants show no systematic
 gains (Tab.~\ref{tab:lkf-number}f).
 \end{description}

 Among the tested settings, \(N_M=2\) was
 most favorable in our configuration. 
Overall, hybrid NN-FEM on the \SI{8}{\kilo\meter} grid substantially increases
\(N_{\mathrm{LKF}}\) compared to the \SI{8}{\kilo\meter} baseline, approaching the
\SI{4}{\kilo\meter} reference with 51 detected LKFs.
(Fig.~\ref{fig:lkf-scatter},
Tab.~\ref{tab:lkf-number}). At \SI{2}{\kilo\meter} resolution, the reference field contains \(171\) LKFs, i.e.,
more than three times the best hybrid NN-FEM counts and roughly an order of
magnitude more than the \SI{8}{\kilo\meter} baseline with \(18\) LKFs
(Tab.~\ref{tab:lkf-number}). 
This discrepancy is expected: although hybrid NN-FEM forms corrections on the
auxiliary grid \(\mathcal T_\text{fine}\), the accepted update is constrained to
the working space \(\V_\text{coarse}\).
Some of the fine-scale features cannot be represented on the coarse mesh as the cells are too wide. Reproducing the \SI{2}{\kilo\meter} LKF density
is therefore not a meaningful target; the appropriate benchmark is the improvement
relative to the \SI{4}{\kilo\meter} resolution. %

\subsubsection{Metrics~\ref{itm:metric-Newton} and~\ref{itm:metric-time}: Newton Iterations and Runtime.}\label{sec:M4}

Fig.~\ref{fig:walltime-scatter} relates the time-to-solution \(T_{\mathrm{mom}}\)
for the sea-ice momentum equation to (a) the error \(E_{\ell^2}\)
(Fig.~\ref{fig:walltime-scatter}a) and (b) the LKF count \(N_{\mathrm{LKF}}\)
(Fig.~\ref{fig:walltime-scatter}b). Runtime scales primarily with the total number
of Newton iterations \(N_{\mathrm{Newton}}\). The favorable region (low error,
low time) is dominated by \(N_M=1\) and \(N_M=2\) configurations: for \(N_M=2\),
networks with two layers of width 256 are particularly efficient. For \(N_M=1\),
four layers of width 512 (or eight layers of width 256) achieve similar errors at
comparable or slightly lower times, but with slightly fewer LKFs. While
\(N_M=0\) leads to fast runtimes, it reduces the error only modestly. Hence, the best
accuracy-efficiency trade-offs are \(N_M=1\) with a medium-large network and
\(N_M=2\) with a compact network.

Fig.~\ref{fig:timings} decomposes the runtime \(T_{\mathrm{mom}}\)
(Fig.~\ref{fig:timings}a) into the Newton-Krylov solver
(Fig.~\ref{fig:timings}b), memory copies for the network
(Fig.~\ref{fig:timings}c), and network evaluation
(Fig.~\ref{fig:timings}d), again as a function of \(N_M\) and \((l,w)\).
Across all configurations, the Newton-Krylov solver dominates the wall time
(\(\sim 200\text{--}350\)\,s), whereas the DNN overhead (copies and inference)
is \(\mathcal{O}(10^{-1})\)\,s in total, i.e., \(<\SI{1}{\percent}\) of the runtime.
Increasing the patch size \(N_M\) does not noticeably raise the cost of copies and
inference. Although the network input and output sizes \((N_{\mathrm{in}},N_{\mathrm{out}})\)
grow with \(N_M\) and \(S\) (Tab.~\ref{tab:io-sizes-side}), their contribution
to the runtime remains minor compared with the Newton-Krylov solver
(cf.~Fig.~\ref{fig:timings}). The problem and data sizes do not saturate the
available hardware. Total time scales with the total Newton iterations
\(N_{\text{Newton}}\) (shown above the bars in Fig.~\ref{fig:timings}), indicating
an approximately constant per-iteration cost.

Network depth (2-8 layers) and width (128-512 units) have only a minor effect on
runtime. For this reason, we only show results for \(w=256\) in
Fig.~\ref{fig:timings}. The inference cost \(T_{\mathrm{nn}}\) grows slightly
with network size but remains below \SI{1}{\percent} of the numerical solver
cost. Thus, further improvements of hybrid NN-FEM should primarily target a
reduction of \(N_{\text{Newton}}\). The copy and inference times each account for
less than \SI{1}{\percent} of the total runtime and have a minor impact on overall
performance. We note that the timings reflect our hardware configuration: the
numerical solver runs on the CPU, whereas the neural network is evaluated on the
GPU. On systems with a different CPU/GPU balance or interconnect bandwidth the
relative shares may change, but on typical hardware setups the qualitative
conclusion, that the wall time is governed primarily by the total number of
Newton iterations \(N_{\text{Newton}}\), should remain valid.

The Newton and linear solvers are able to solve the nonlinear systems of
equations, but are not fully \(h\)-robust for these equations: iteration counts
increase with refinement and the per-iteration cost grows with the number of
degrees of freedom. Consequently, the time to solution for the momentum equation
at \SI{2}{\kilo\meter} resolution exceeds that at \SI{8}{\kilo\meter} by a factor
\(\approx 100\). In contrast, hybrid NN-FEM applied on the \SI{8}{\kilo\meter}
grid achieves both lower error and lower time-to-solution than the
\SI{8}{\kilo\meter} baseline (Fig.~\ref{fig:walltime-scatter},
Fig.~\ref{fig:timings}). The cost of data generation in the experiments is
essentially a single simulation of the benchmark at \SI{2}{\kilo\meter}
resolution. On our hardware, training a single network took less than a minute
across the tested configurations, as training is not compute-bound for these
small models.

\subsubsection{Overall Performance Summary.}\label{sec:Msum}
We summarize the discussed performance of the NN-FEM approach in Section \ref{sec:M1}-\ref{sec:M4}. Across all runs, the \emph{patch size} \(N_M\) is the dominant driver of prediction
quality, while the \emph{jump level} \(S\) has a secondary influence and the \emph{architecture}
\((l,w)\) primarily shifts absolute LKF counts without changing the overall ranking
(Tab.~\ref{tab:lkf-number}a--\ref{tab:lkf-number}c versus
Tab.~\ref{tab:lkf-number}d--\ref{tab:lkf-number}f)

Compared to the coarse baseline simulation performed on the
\SI{8}{\kilo\meter}, which has
\(T_{\mathrm{mom}} = \SI{277.8}{\second}\),
\(E_{\ell^2} = 0.528\),
\(N_{\mathrm{LKF}} = 18\),
and \(N_{\mathrm{Newton}} = 680\), the compact hybrid NN-FEM configuration at \(S=2\), \(N_M=1\) with
\((l=4,\, w=256)\) achieves
\(T_{\mathrm{mom}} = \SI{252.45}{\second}\)
(\SI{-9.1}{\percent} compared to \SI{8}{\kilo\meter}),
\(E_{\ell^2} = 0.042\),
\(N_{\mathrm{LKF}} = 39\) (about \(2.2\times\) the \SI{8}{\kilo\meter} baseline),
and \(N_{\mathrm{Newton}} = 551\) (\SI{-19}{\percent}).
The \SI{4}{\kilo\meter} reference requires
\(T_{\mathrm{mom}} = \SI{2751.16}{\second}\)
(\(\approx 10.9\times\) the hybrid NN-FEM runtime), with
\(N_{\mathrm{LKF}} \approx 51\).
Across all tested hybrid NN-FEM settings,
\(T_{\mathrm{mom}}\) lies in the narrow range
\SI{250}{}--\SI{262}{\second}, and the network overhead (copies and inference)
remains \(<\SI{1}{\percent}\) of \(T_{\mathrm{mom}}\)
(Fig.~\ref{fig:timings}). In summary, hybrid NN-FEM on the
\SI{8}{\kilo\meter} grid attains accuracy and LKF statistics comparable to the
\SI{4}{\kilo\meter} reference at roughly one-tenth of the Newton run-time of the high-resolution
simulations, while keeping the additional neural network overhead negligible.

\subsection{Ablation Studies on Neural Network Inputs}\label{sec:study-ablation}
We test the \emph{causal} role of network inputs by masking them while keeping
architecture, weights, data, and the solver as well as finite element pipeline fixed. In
particular, we assess the importance of the state and residual inputs by masking
them in the batched input
\[
 \mathbf X = \big[ \vt_{\mathrm{batched}} \ \big|\ \rt_{\mathrm{batched}}
 \ \big|\ \boldsymbol\omega \big].
\]
Geometry $\boldsymbol\omega$ is kept unchanged, as it is constant across the
whole batch (uniform Cartesian refinement). A masked input is \emph{dispensable}
if stability is unchanged (no divergence, no increase in
\(N_{\mathrm{Newton}}\)~\ref{itm:metric-Newton}) and accuracy degrades at most
marginally; otherwise it is \emph{essential}. Removing $\vt_{\mathrm{batched}}$ tests
\emph{state-awareness} (local Jacobian scaling), removing $\rt_{\mathrm{batched}}$ tests
\emph{residual-awareness} (magnitude/direction), cf.~Remark~\ref{rem:nn-design}.
Poor neural network predictions usually leads to residual growth and Newton divergence.

So far, we considered $\mathbf X_{\text{full}} = [\vt_{\mathrm{batched}} \ \big|\ \rt_{\mathrm{batched}}
\ \big|\ \boldsymbol\omega]$. The goal is to compare the performance of hybrid NN-FEM when
one of the inputs is masked:
\[
 \mathbf X_{(\vt_{\mathrm{batched}=0)}} \coloneq [\mathbf 0 \mid \rt_{\mathrm{batched}} \mid \boldsymbol\omega],\qquad
 \mathbf X_{(\rt_{\mathrm{batched}=0)}} \coloneq [ \vt_{\mathrm{batched}} \mid \mathbf 0 \mid \boldsymbol\omega].
\]
For each $\alpha\in\{\vt_{\mathrm{batched}}=0,\rt_{\mathrm{batched}}=0\}$, we evaluate the network and assemble the
correction as in Section~\ref{sec:patch-ops},
\[
 \mathbf D_{(\alpha)}=\mathcal N(\mathbf X_{(\alpha)}),\qquad
 \delta\vt_{\text{fine},(\alpha)} = BW\sum_M S_M\mathbf D_{(\alpha);\,M,:},\qquad
 {\vt}_{\text{fine},(\alpha)} = \vt_\text{fine} + \delta\vt_{\text{fine},(\alpha)}\,.
\]
\paragraph{Inter-step Effect of the Hybrid NN-FEM Correction.}
The hybrid NN-FEM update is applied \emph{after} Newton has converged at $t^{n+1}$ and
thus perturbs only the right-hand side of the \emph{next} step. As discussed in
Remark~\ref{rem:dnnmg-stab}, standard inexact Newton considerations suggest that
this perturbation must be appropriately scaled relative to the unperturbed
residual to retain local convergence. We therefore regard the hybrid NN-FEM correction
as a controlled perturbation and \emph{test} this expectation via input
ablations.

\paragraph{Ablation $(\vt_{\mathrm{batched}}=0)$: Residual-only Input.}
We evaluate the residual-only variant by zeroing the state block in the batched
network input. In all experiments, the correction constructed within hybrid
NN-FEM increased the nonlinear residual and the subsequent Newton iteration
diverged, yielding no usable solutions. On a patch $N_M$, an ideal stabilizing
update solves (approximately) the local linearized equation
$\mathcal J\restrict{M}(\vt)\delta\vt\restrict{M}\approx-\mathbf
R\restrict{M}(\vt)$ with
$\mathcal J\restrict{M}(\vt)=\partial_\vt \mathcal R\restrict{M}(\vt)$, (cf.~\eqref{eq:newton-iterate})
Using the splitting time approach to approximate the sea-ice dynamics the
Jacobian of the momentum equation~\eqref{eq:Newvar}, $\mathcal J\restrict{M}(\vt)$, depends on the strain-rate invariants
through $\eta(\vt),\zeta(\vt)$ and is therefore state dependent and
anisotropic. Masking $\vt_{\text{batched}}$ forces the network to produce
$\delta\vt\restrict{M}$ without
access to this local scaling/rotation. While this may still reduce a
\emph{same-step} fine-grid residual in expectation, it provides no control of
the forcing in the next time step. Empirically, the correction
$\delta\vt_{\text{fine}}$, and, in turn, the perturbation of
$\mathbf f_{\text{fine}}$ becomes too large (or misaligned), such
that condition~\eqref{eq:newton-preservation} presented in the Appendix is not fulfilled anymore and the next Newton
solve diverges. Thus, the state input is necessary to calibrate
$\delta\vt_{\text{fine}}^{n+1}$ to the local Jacobian and keep the perturbation induced by
$\delta\vt_{\text{fine}}$ within admissible bounds.

\paragraph{Ablation $(\rt_{\mathrm{batched}}=0)$: State-only Input.}
We evaluate the state-only variant by zeroing the residual block in the batched
network input. Similar to the state ablation study, the hybrid NN-FEM correction
leads to diverging Newton iterations and unusable solutions. If the residual
block is masked, the network outputs $\delta\vt_{\text{fine}}$ with no
dependence on the local residual. Consequently, $\delta\vt_{\text{fine}}^{n}$
does not scale with the \emph{current} residual and, in particular, does not
vanish when the residual is small. The perturbation of $\mathbf f_{\text{fine}}$
induced by $\delta\vt_{\text{fine}}^{n}$ is then uncontrolled in magnitude and
direction relative to the residual, again causing Newton divergence at
$t^{n+1}$ due to the violation of condition~\eqref{eq:newton-preservation} stated in the Appendix.

Because the hybrid NN-FEM correction acts \emph{between} time steps as an
explicit forcing, we conclude that stability requires \emph{both} (i) state
information $\vt_{\mathrm{batched}}$ to approximate the local action of
$\mathcal J(\vt)^{-1}$ (scaling/rotation), and (ii) residual information
$\rt_{\mathrm{batched}}$ to set the correct direction and scale. Masking either input
destroys this coupling, yields a destructive perturbation, and destabilizes the
Newton-Krylov solve in the next time step. This confirms our hypothesis in
Remark~\ref{rem:nn-design}. However, at this point, we do not prove a general
guarantee beyond the empirical observations made here.

\section{\label{sec:conc}Conclusion}
In this paper, we developed a hybrid NN–FEM approach to efficiently simulate sea-ice deformation, in particular linear kinematic features (LKFs), on coarse horizontal grids. This is especially advantageous given that, in the viscous–plastic model, resolving LKFs typically requires grids with high spatial resolution. The hybrid NN–FEM approach for the viscous–plastic sea-ice model augments an under-resolved finite element solution with patch-local neural corrections, evaluated in a richer auxiliary space and subsequently assembled into a global correction. The construction preserves Dirichlet conditions.

With regard to the simulation of LKFs, our parameter study shows that patch size (the union of neighboring mesh elements of the fine grid), $N_M$, is the most important influencing factor. Larger patches provide a richer state context and better alignment with dominant deformation modes. We observe that an increasing number of refinement levels used in the training of the NN leads to a more accurate representation of LKFs, although this effect is less dominant than an increased patch size. A larger jump ($S=2$) enables the correction to account for finer structures. The network depth and width mainly shift the absolute values of the LKFs and have a weaker influence on their representation quality compared to increasing $S$. In
terms of accuracy-efficiency, $N_M=2$, $S=2$ with a compact MLP (two hidden
layers, width 256) performs best. When training data are scarce, $N_M=1$, $S=2$
offers a stable trade-off with comparable runtime and slightly fewer LKFs. In our analysis, we observe that very large patches ($N_M=3$) are data-sensitive. This is likely due to the fact that the feature dimension grows while the
number of statistically independent patches per snapshot decreases, which can
induce noisy corrections unless substantially more data (or stronger inductive
bias) are available.

Across all settings, the nonlinear and linear solves dominate
the runtime of solving the sea-ice momentum equation as typical for the sea-ice dynamics. The NN overhead (copies and inference) remains below
\SI{1}{\percent} on standard CPU/GPU nodes.
The learned correction even accelerates the Newton method compared to the reference Newton solve on the coarse grid. In this setup, the solution of the momentum equation is up to \SI{10}{\percent} faster with the NN–FEM approach than without the NN correction.
Moreover, the hybrid NN–FEM elevates a \SI{8}{\kilo\meter} baseline approximation to a level of accuracy qualitatively comparable to a \SI{4}{\kilo\meter} discretization, while the \SI{4}{\kilo\meter} run is nearly $11\times$ more expensive than the hybrid NN–FEM.

The data requirements of very large patches suggest benefits from additional
inductive bias (e.g., equivariance to rotations/reflections) or adaptive
patching. A possible direction is \emph{online}, semi-supervised training with
physics-aware objectives that penalize fine-space momentum residuals after
correction and incorporate stability cues from local sensitivity information.
Such training has the potential to reduce labeled-data needs, improve
generalization across different forcing terms and meshes, and further decrease Newton
iterations without increasing network capacity.

\appendix
\section{\label{sec:technical}Technical Details of the Hybrid Neural Network-Finite Element Method}
In Alg.~\ref{alg:dnnmg} the correction at $t^{n}$ modifies the
state $\vt_\text{coarse}$ obtained in Alg.~\ref{alg:dnnmg}, l.~\ref{alg:dnnmg:S1} at $t^n$.
Thereby, the hybrid NN-FEM modifies the time integrator, as the previous state is used in
the backward Euler at $t^{n+1}$. Hence, the hybrid NN-FEM correction in the current time
step alters the right-hand side of the next step. Implementation-wise, we
assemble all \emph{load} terms on an auxiliary fine space using the corrected
field and obtain the right-hand side in $\V_\text{coarse}$ by \emph{restriction}. The
notation below makes this restriction/prolongation precise. The subsequent
remarks state the modified right hand side $\mathbf f_\mathrm{fine}$ assembly and quantify the perturbation through
the correction, together with a simple smallness condition that guarantees
admissibility for the inexact Newton solve at $t^{n+1}$. Beforehand, we need to
introduce some additional notation.

To define the restricted right-hand side, we use the restriction and prolongation
operators between the fine and coarse spaces. Let
$f_{fine}\in \V_{\text{fine}}$. Its \emph{restriction} to the coarse space $\V_{\text{coarse}}$ is defined by the
Petrov-Galerkin duality
\[
 ( \mathbf R f_{\text{fine}},\, w_{\text{coarse}})
 \coloneq
 ( f_{\text{fine}},\, \mathbf P w_{\text{coarse}})
 \qquad \forall\, w_{\text{coarse}}\in \V_{\text{coarse}}\,.
\]
In discrete form, if $\mathbf f_{\text{fine}}$ denotes the assembled fine-grid
right-hand side with entries
$(\mathbf f_{\text{fine}})_i = (f_{\text{fine}}, {\varphi_{\text{fine}}}_i)$,
then the restricted coarse right-hand side is obtained by
\[
 \mathbf f_{\text{coarse}} = \mathbf R \mathbf f_{\text{fine}},
 \qquad\text{with}\qquad \mathbf R = \mathbf P^{\top}.
\]

\begin{remark}[Impact of the correction on inexact Newton admissibility]\label{rem:dnnmg-stab}
Let $\mathcal R_{L,\,\mathrm{base}}(\cdot)$ denote the residual
built in $\V_\text{coarse}$ with the \emph{uncorrected} history state $\vt_\text{coarse}$,
and let $\mathcal R$ be the residual built with the
\emph{corrected} state. The modification is an additive
functional
\[
 \mathcal R(\vt_\text{coarse},\,\wt_\text{coarse};\,\delta\vt_\text{fine})
 = \mathcal R_{\mathrm{base}}(\vt_\text{coarse},\,\wt_\text{coarse})
 + \gt_\text{fine}(\delta\vt_\text{fine},\,\wt_\text{coarse}),
\]
where
\[
 \gt_\text{fine}(\delta\vt_\text{fine},\,\wt_\text{coarse})
 = \Big(\tfrac{\rho}{k}\,\delta\vt_\text{fine},\,\mathbf P\wt_\text{coarse}\Big).
\]
Let $\vt_\text{coarse}^{(0)}\in \V_\text{coarse}$ be the initial Newton guess at $t^{n+1}$
(typically the state from $t^{n}$). By the triangle inequality,
\[
\big|\,\|\mathcal R(\vt_\text{coarse}^{(0)},\,\wt_\text{coarse};\,\delta\vt_\text{fine})\|
-\|\mathcal R_{\mathrm{base}}(\vt_\text{coarse}^{(0)},\,\wt_\text{coarse})\|\,\big|
\le \|\gt_\text{fine}\|.
\]
Thus, if
\begin{equation}\label{eq:smallness-assumption-i}
 \|\gt_\text{fine}\|
 \le \beta\|\mathcal R_{\mathrm{base}}(\vt_\text{coarse}^{(0)},\,\wt_\text{coarse})\|,
 \qquad \beta\in(0,1),
\end{equation}
then
\begin{equation}\label{eq:newton-preservation}
 (1-\beta)\|\mathcal R_{\mathrm{base}}(\vt_\text{coarse}^{(0)},\,\wt_\text{coarse})\|
 \le \|\mathcal R(\vt_\text{coarse}^{(0)},\,\wt_\text{coarse};\,\delta\vt_\text{fine})\|
 \le (1+\beta)\|\mathcal R_{\mathrm{base}}(\vt_\text{coarse}^{(0)},\,\wt_\text{coarse})\|.
\end{equation}
Let \(\eta_0 \in [0,1)\) denote the inexact Newton forcing parameter controlling the relative accuracy of the linear solve in the first Newton step.
The inexact Newton linear solve at the first step is imposed on the
\emph{perturbed} residual:
\[
 \big\|\mathcal R(\vt_\text{coarse}^{(0)},\,\wt_\text{coarse};\,\delta\vt_\text{fine})
 + \mathcal J_\text{coarse}(\vt_\text{coarse}^{(0)})[\delta\vt_\text{coarse}^{(0)},\,\wt_\text{coarse}]\big\|
 \le \eta_0\big\|\mathcal R(\vt_\text{coarse}^{(0)},\,\wt_\text{coarse};\,\delta\vt_\text{fine})\big\|,
\qquad \eta_0<1,
\]
and the above bound~\eqref{eq:newton-preservation} ensures that, provided the unperturbed guess lies in the
local Newton convergence region, the perturbed one remains in the convergence
region (with a radius reduced at most by the factor $1+\beta$). Thus, one should enforce the bound~\eqref{eq:smallness-assumption-i}.
\end{remark}

\section{Data Statistics}
In this section, we summarize the main characteristics of the data sets used in
our study. We report descriptive statistics for the training and test data as described in Sec.~\ref{sec:nn-par-and-train}.
This includes ranges and distributions of the input parameters and relevant output
quantities (see Tab.~\ref{tab:groupwise-stats-final} for grouped statistics of inputs
\((R_x,R_y,V_x,V_y)\) and outputs \((u_1,u_2)\) by patch size \(N_M\) and jump level \(S\)).

\sisetup{
 scientific-notation = true,
 round-mode = places,
 round-precision = 3
}
\begin{table}[t]
\centering
\caption{Grouped statistics (mean, std, min, max) for inputs ($R_x,R_y,V_x,V_y$) and outputs ($u_1,u_2$) by patch size $N_M$ and jump level $S$. Geometry rows are omitted.}
\label{tab:groupwise-stats-final}
\begin{tabular}{ccllS[table-format=+1.3e-2]
 S[table-format=+1.3e-2]
 S[table-format=+1.3e-2]
 S[table-format=+1.3e-2]}
\toprule
$N_M$ & $S$ & Kind & Group & {mean} & {std} & {min} & {max} \\
\midrule
0 & 1 & In & $R_x$   & 1.8347873861e-11 & 2.15933829793e-09 & -4.51514298678e-08 & 4.11465608129e-08 \\
0 & 1 & In & $R_y$   & -4.05724491448e-11 & 2.07449227366e-09 & -4.65420648674e-08 & 4.11968744825e-08 \\
0 & 1 & In & $V_x$ & 1.39630205837e-06 & 6.35819332314e-05 & -1.61367661691e-04 & 1.79014554894e-04 \\
0 & 1 & In & $V_y$ & 8.12178631694e-07 & 6.35378248048e-05 & -1.59044202183e-04 & 1.78911044233e-04 \\
0 & 1 & Out & $u_1$   & 9.4556332856e-09 & 1.83386855416e-06 & -4.08221607045e-05 & 5.7146642637e-05 \\
0 & 1 & Out & $u_2$   & -1.15764643151e-07 & 2.95573789183e-06 & -9.29173729574e-05 & 5.7146642637e-05 \\
\midrule
0 & 2 & In & $R_x$   & -1.44903600828e-11 & 9.30551376062e-10 & -2.58851541834e-08 & 2.27764784082e-08 \\
0 & 2 & In & $R_y$   & 6.60390862705e-12 & 9.73042220553e-10 & -2.25118033960e-08 & 2.23310733976e-08 \\
0 & 2 & In & $V_x$ & 1.18493397840e-06 & 6.36260512984e-05 & -1.63699263261e-04 & 1.79564441510e-04 \\
0 & 2 & In & $V_y$ & 1.00453730013e-06 & 6.35370389972e-05 & -1.59047801039e-04 & 1.79645803947e-04 \\
0 & 2 & Out & $u_1$   & 2.52409393595e-08 & 2.89806793897e-06 & -7.45123199686e-05 & 1.01293884001e-04 \\
0 & 2 & Out & $u_2$   & -2.25698884481e-07 & 4.43796236537e-06 & -1.04717258327e-04 & 1.01631729579e-04 \\
\midrule
1 & 1 & In & $R_x$   & -1.32757684834e-11 & 2.21793176923e-09 & -4.65420648674e-08 & 4.11465608129e-08 \\
1 & 1 & In & $R_y$   & 3.78419200331e-12 & 2.15200716448e-09 & -4.53045577487e-08 & 4.11968744825e-08 \\
1 & 1 & In & $V_x$ & 1.34481129386e-06 & 6.34259105150e-05 & -1.59044202183e-04 & 1.79014554894e-04 \\
1 & 1 & In & $V_y$ & 8.70321149280e-07 & 6.35536740049e-05 & -1.61367661691e-04 & 1.77594188274e-04 \\
1 & 1 & Out & $u_1$   & 4.49137825136e-08 & 1.82935881580e-06 & -4.21330510720e-05 & 5.71466426370e-05 \\
1 & 1 & Out & $u_2$   & -1.59827082538e-07 & 3.13779848213e-06 & -9.29173729574e-05 & 5.71466426370e-05 \\
\midrule
1 & 2 & In & $R_x$   & 2.18192793878e-12 & 9.27436205370e-10 & -2.56575027209e-08 & 2.23304472604e-08 \\
1 & 2 & In & $R_y$   & -4.70096333123e-12 & 9.44881475921e-10 & -2.58851541834e-08 & 2.27764488445e-08 \\
1 & 2 & In & $V_x$ & 1.29377721499e-06 & 6.34817791884e-05 & -1.59047966976e-04 & 1.79645728384e-04 \\
1 & 2 & In & $V_y$ & 8.98536633452e-07 & 6.36080124906e-05 & -1.63699430786e-04 & 1.78800904449e-04 \\
1 & 2 & Out & $u_1$   & 5.85021824347e-08 & 2.94199365155e-06 & -7.46877228327e-05 & 1.01293861713e-04 \\
1 & 2 & Out & $u_2$   & -2.65214644759e-07 & 4.56114423157e-06 & -1.04717287025e-04 & 1.01631631518e-04 \\
\midrule
2 & 1 & In & $R_x$   & 1.72666890247e-12 & 2.27514016550e-09 & -4.65420648674e-08 & 4.06122729915e-08 \\
2 & 1 & In & $R_y$   & -9.10145587994e-13 & 2.19157039844e-09 & -4.53045577487e-08 & 4.11968744825e-08 \\
2 & 1 & In & $V_x$ & 1.81294631766e-06 & 6.30617232728e-05 & -1.48834354466e-04 & 1.74077158302e-04 \\
2 & 1 & In & $V_y$ & 4.03378161998e-07 & 6.38037789113e-05 & -1.61367661691e-04 & 1.79014554894e-04 \\
2 & 1 & Out & $u_1$   & 4.59290486835e-08 & 1.89297822780e-06 & -3.29728878026e-05 & 5.71466426370e-05 \\
2 & 1 & Out & $u_2$   & -1.63774932320e-07 & 3.20311280589e-06 & -9.29173729574e-05 & 5.71466426370e-05 \\
\midrule
2 & 2 & In & $R_x$   & -1.98097681309e-13 & 9.12178856240e-10 & -2.56575027209e-08 & 2.27757087970e-08 \\
2 & 2 & In & $R_y$   & 8.69119089784e-13 & 9.51481224019e-10 & -2.58851541834e-08 & 2.11364453591e-08 \\
2 & 2 & In & $V_x$ & 1.82425554622e-06 & 6.31424076642e-05 & -1.48838826583e-04 & 1.76489774143e-04 \\
2 & 2 & In & $V_y$ & 3.68085095744e-07 & 6.38858252469e-05 & -1.63699676006e-04 & 1.79644540484e-04 \\
2 & 2 & Out & $u_1$   & 5.53189449569e-08 & 3.00326931560e-06 & -6.88438812191e-05 & 1.01292846291e-04 \\
2 & 2 & Out & $u_2$   & -2.62993947989e-07 & 4.59769652850e-06 & -1.04717264157e-04 & 1.01629868620e-04 \\
\bottomrule
\end{tabular}
\end{table}
\section{Extended Results from the Numerical Experiments}
\begin{figure}
\begin{tikzpicture}[font=\small]
\begin{groupplot}[
 group style={group size=1 by 4, horizontal sep=1cm, vertical sep=1.5cm},
 width=\linewidth, height=4.5cm,
 ymajorgrids, ymin=0,
 xtick=data, title style={at={(0.5,0.95)}},
 x tick label style={rotate=45, yshift=4pt, xshift=2pt, font=\footnotesize,
  anchor=north east},
 xticklabels from table={\sumtab}{group_key},
 enlarge x limits=0.03,
 colormap/viridis,
 cycle list name=ps,
 legend columns=7, legend cell align=left,
 legend style={/tikz/every even column/.append style={column sep=0.2cm},at={(0.5,1.05)},draw=none,anchor=south},
]
\nextgroupplot[ylabel={Seconds}, title style={at={(0.5,1.2)}}, title={Total hybrid NN-FEM timing},ymax=450]
\addlegendimage{empty legend}
\addlegendentry{\textit{Numbers above bars: }Total Newton iters}
\addlegendimage{empty legend}
\addlegendentry{$\quad$}
\addlegendimage{empty legend}
\addlegendentry{\textit{Patch size}}
\addlegendimage{only marks, mark size=4pt, mark=square*, draw=mapped color, fill=mapped color, index of colormap=0}
\addlegendentry{$N_M=0$}
\addlegendimage{only marks, mark size=4pt, mark=square*, draw=mapped color, fill=mapped color, index of colormap=8}
\addlegendentry{$N_M=1$}
\addlegendimage{only marks, mark size=4pt, mark=square*, draw=mapped color, fill=mapped color, index of colormap=17}
\addlegendentry{$N_M=2$}
\addplot+[forget plot,
ybar, bar width=0pt,fill,
error bars/.cd, y dir=both, y explicit
] table[
x expr=\coordindex,
y=total_mean,
]{\sumtab};
\foreach \P in {1,2,3,5}{%
\addplot+[forget plot,
ybar, bar width=0.8pt,fill,
error bars/.cd, y dir=both, y explicit
] table[
x expr=\coordindex,
y=total_mean,
y error=total_std,
restrict expr to domain={\thisrow{ps}}{\P:\P}
]{\sumtab};
\addplot+[forget plot,only marks, nodes near coords,
visualization depends on=\thisrow{newton_sum_mean}\as\newtonlabel,
nodes near coords style={font=\scriptsize,text=black,yshift=12pt, xshift=6pt,rotate=90},
nodes near coords={\pgfmathprintnumber[precision=0]{\newtonlabel}}
] table[
x expr=\coordindex,
y=total_mean,
restrict expr to domain={\thisrow{ps}}{\P:\P}
]{\sumtab};
\addplot+[only marks, nodes near coords,
nodes near coords style={font=\scriptsize,text=black,fill=white, yshift=-1cm, xshift=6pt,rotate=90},
] table[
x expr=\coordindex,
y=total_mean,
restrict expr to domain={\thisrow{ps}}{\P:\P}
]{\sumtab};
}
\nextgroupplot[ylabel={Seconds}, title={Timing of the Newton-Krylov solver}]
\addplot+[forget plot,
ybar, bar width=0pt,fill,
error bars/.cd, y dir=both, y explicit
] table[
x expr=\coordindex,
y=solver_mean,
]{\sumtab};
\foreach \P in {1,2,3,5}{%
\addplot+[forget plot,
ybar, bar width=0.8pt,fill,
error bars/.cd, y dir=both, y explicit
] table[
x expr=\coordindex,
y=solver_mean,
y error=solver_std,
restrict expr to domain={\thisrow{ps}}{\P:\P}
]{\sumtab};
\addplot+[only marks, nodes near coords,
nodes near coords style={font=\scriptsize,text=black,yshift=-1cm,fill=white,xshift=6pt,rotate=90},
] table[
x expr=\coordindex,
y=solver_mean,
restrict expr to domain={\thisrow{ps}}{\P:\P}
]{\sumtab};
}
\nextgroupplot[ylabel={Seconds}, title={Timing of copies for network}]
\addplot+[forget plot,
ybar, bar width=0pt,fill,
error bars/.cd, y dir=both, y explicit
] table[
x expr=\coordindex,
y=copies_for_network_mean,
restrict expr to domain={\thisrow{ps}}{1:3}
]{\sumtab};
\foreach \P in {1,2,3}{%
 \addplot+[forget plot,
 ybar, bar width=0.8pt,fill,
 error bars/.cd, y dir=both, y explicit
] table[
 x expr=\coordindex,
 y=copies_for_network_mean,
 y error=copies_for_network_std,
 restrict expr to domain={\thisrow{ps}}{\P:\P}
]{\sumtab};
\addplot+[
 only marks, nodes near coords,
 point meta=y,
 nodes near coords style={font=\scriptsize,text=black, yshift=-0.66cm,fill=white, xshift=6pt,rotate=90},
] table[
 x expr=\coordindex,
 y=copies_for_network_mean,
 restrict expr to domain={\thisrow{ps}}{\P:\P}
]{\sumtab};
}

\nextgroupplot[ylabel={Seconds}, title={Timing of network evaluation}]
\addplot+[forget plot,
ybar, bar width=0pt,fill,
error bars/.cd, y dir=both, y explicit
] table[
x expr=\coordindex,
y=network_mean,
restrict expr to domain={\thisrow{ps}}{1:3}
]{\sumtab};
\foreach \P in {1,2,3}{%
\addplot+[forget plot,
 ybar, bar width=0.8pt,fill,
 error bars/.cd, y dir=both, y explicit
] table[
 x expr=\coordindex,
 y=network_mean,
 y error=network_std,
 restrict expr to domain={\thisrow{ps}}{\P:\P}
]{\sumtab};
\addplot+[only marks, nodes near coords,
nodes near coords style={font=\scriptsize,text=black, yshift=-1cm,fill=white, xshift=6pt,rotate=90},
] table[
 x expr=\coordindex,
 y=network_mean,
 restrict expr to domain={\thisrow{ps}}{\P:\P}
 ]{\sumtab};
}
\end{groupplot}
\end{tikzpicture}
\caption{\label{fig:timings-detail}Timing of hybrid NN-FEM for $S=2$ across patch sizes
 \(N_M\in\{0,1,2\}\) (colors) and network widths/depths: total runtime, Newton-Krylov
 solver time, host-device copy time, and network inference time; numbers above
 bars denote total Newton iterations \(N_{\text{Newton}}\). }
\end{figure}
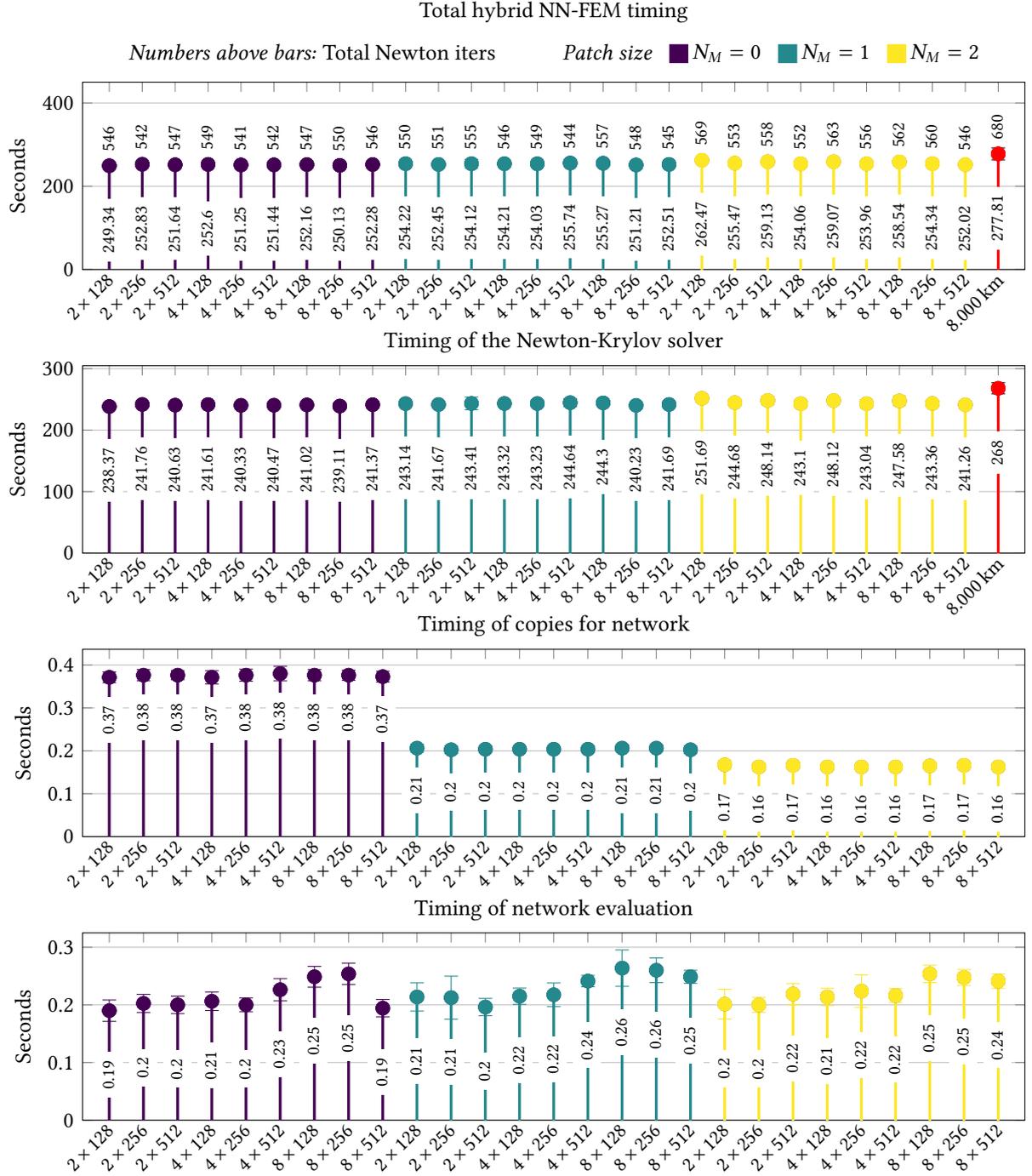

In this section we report extended timing results for the hybrid NN-FEM solver.
Fig.~\ref{fig:timings-detail} provides a detailed breakdown for $S=2$ across
patch sizes $N_M\in\{0,1,2\}$ and all tested network widths and depths, showing
total runtime, Newton-Krylov solver time, host-device copy time, and network
inference time; the numbers above the bars indicate the total number of Newton
iterations $N_{\text{Newton}}$. These data complement the timing results shown
in Fig.~\ref{fig:timings} by including all considered network architectures and
patch-size/jump-level configurations.

\printbibliography
\end{document}